\documentclass{amsart}
 \usepackage{amsmath}
 \usepackage{epsfig} 
 \usepackage{graphics} 

\usepackage{amsthm}
\usepackage{amssymb}
\usepackage{amsfonts,mathrsfs}
\usepackage{stmaryrd}
\usepackage{amsxtra}  
\usepackage{verbatim}

\newtheorem{theorem}{Theorem}
\newtheorem{corollary}{Corollary}

\newtheorem{lemma}{Lemma}
\newtheorem{proposition}{Proposition}

\newtheorem{problem 1}{Problem 1}
\newtheorem{problem 2}{Problem 2}
\newtheorem{problem 3}{Problem 3}
\theoremstyle{definition}
\newtheorem{definition}{Definition}
\newtheorem{remark}{Remark}

\newcommand{\bea}{\begin{eqnarray*}}
\newcommand{\eea}{\end{eqnarray*}}

\title[ERGODICITY  OF HARMONIC CURRENTS ]
      {UNIQUE ERGODICITY  OF HARMONIC CURRENTS ON SINGULAR FOLIATIONS
OF $\mathbb{P}^2$}

\author[John Erik Forn\ae ss, Nessim Sibony]{}

\email{fornaess@umich.edu, nessim.sibony@math.u-psud.fr}

\begin{document}

\begin{abstract}
Let ${\mathcal F}$ be a holomorphic foliation of $\mathbb{P}^2$ by Riemann
surfaces. Assume all the singular points of ${\mathcal F}$ are hyperbolic. If ${\mathcal F}$
has no algebraic leaf, then there is a unique
positive harmonic $(1,1)$ current $T$ of mass one, directed by ${\mathcal F}.$ This implies strong
ergodic properties for the foliation ${\mathcal F}.$ We also study the harmonic flow associated to
the current $T.$
  \end{abstract}

\maketitle



\centerline{John Erik Forn\ae ss\footnote{The first author is 
supported by an NSF grant.
} and Nessim Sibony
}

\medskip


\medskip

 \medskip


\section{Introduction}

Let ${\mathcal F}$ be a holomorphic foliation of the complex projective space
$\mathbb{P}^2$. Our purpose is to study the ergodic properties of
$\mathcal F,$ using the theory of harmonic currents as developed by the
authors in \cite{FS2005}. 

\medskip

A holomorphic foliation can be seen as a rational vector field in $\mathbb{C}^2.$ Our goal is to
develop an ergodic theory for the dynamics of such vector fields. The two main
difficulties are: the presence of singularities (they always exist) and the
absence (generically) of algebraic leaves. And hence it is not clear where to start the
analysis. Our method is geometric but requires difficult estimates. To our knowledge, this
is the first paper where global dynamical results for rational
vector fields are obtained. The subject is classical and related to polynomial vector
fields in $\mathbb{R}^2.$

\medskip

We first recall a few facts. Let $\pi: \mathbb{C}^3 \setminus \{0\}\rightarrow \mathbb{P}^2$
denote the canonical projection. The foliation $\pi^*{\mathcal F}$ can be defined
in $\mathbb{C}^3$ by a global $1-$form  $\omega_0=a_1(x)dx_1+ a_2(x)dx_2+a_3(x)dx_3$
where  the $a_j(x)$ are homogeneous polynomials of 
the same degree $\delta \geq 1$ without common factors. 
Moreover since every line through the
origin is in the kernel of $\omega_0,$ they satisfy the condition $\sum x_i a_i(x)=0$ 

\medskip

The degree of ${\mathcal F}$ is by definition deg ${\mathcal F}=d:=$deg$\;\delta-1.$
It represents the number of tangencies of a generic line $L,$ with ${\mathcal F}.$
Let $Fol(d)$ denote the space of foliations of degree $d.$  The space of coefficients of $1$
forms of degree $\delta$ is a projective space. The subspace given by $\sum x_ia_i=0$
is  a linear subspace, so also a projective space. The subspace of $1$ forms of degree $\delta$
of the form $H\lambda$ where $H$ is a homogenous polynomial of degree $0<\delta'<\delta$ and $\lambda $
is a $1-$ form of degree $\delta-\delta'$ is an algebraic subvariety. So together this gives that
$Fol(d)$ is the complement of an algebraic subvariety of some $\mathbb{P}^N.$
 It follows from the B\'ezout
theorem that the foliation $\mathcal F$ has a finite number of singularities bounded
uniformly by some function of the degree. If in a coordinate chart $U,\mathcal F$
is defined by $\omega_1=\alpha(z,w)dz+\beta(z,w)dw,$ then sing$(\mathcal F)\cap U=
\{\alpha=\beta=0\}.$ We can assume that all the singular points are in the same
$\mathbb{C}^2$, $\{p_j=(\alpha_j,\beta_j)\}_{j\leq N}.$

\medskip

\begin{definition} Suppose there is a 
change of coordinates around $p_j$ sending $p_j$ to $0$ and such that
$\omega_0(z,w)= zdw-\lambda wdz +{\mathcal O}(z,w)^2$ where
$\lambda=a+ib$ and $b$ is a nonzero number. We say in this case that 
the singularity is hyperbolic and that we are in the Poincar\'e domain. 

\end{definition}

The following is a classical fact due to Poincar\'e, see \cite{Ch1986}.

\begin{theorem}
Suppose that the singular point is hyperbolic. Then there is a local biholomorphic
change of coordinates so that the form $\omega_0$ in these coordinates can be written
$\omega_0=zdw-\lambda wdz$ (with the same $\lambda$).
\end{theorem}

We remark that the form $\omega_0$ is invariant under scaling except for multiplication by a 
constant which of course does not affect the zero set. Hence we can assume
that the linearization is valid in a fixed large ball, in particular in a neighborhood
of the unit bidisc.\\

The following result is due to Lins Neto, Soares \cite{LNSo1996} (we give only the two
dimensional version, their result is also valid in $\mathbb P^k$). The result uses Jouanolou's
example of a foliation in $\mathbb P^2$ without algebraic leaves.

\begin{theorem}
There exists a real Zariski dense open subset ${\mathcal H} (d)\subset Fol(d)$
such that any ${\mathcal F}\subset {\mathcal H}(d)$ satisfies:\\
i) ${\mathcal F} $ has only
hyperbolic singularities and no other singular points.\\
ii) ${\mathcal  F}$ has no invariant algebraic curve.
\end{theorem}

The global behavior of foliations is not well understood. It is unknown
whether every leaf of a given foliation ${\mathcal F}$, clusters
at a singular point. This problem, known as the  problem of existence
of a minimal exceptional set is discussed in \cite{CLS1988}
and \cite{BLM1992} for example.
It is conjectured in \cite{Il1987} that a generic holomorphic foliation by Riemann surfaces
in $\mathbb P^k$ has dense leaves. Recently Loray and Rebelo \cite{LR2003}
have constructed non empty open sets of holomorphic foliations by Riemann
surfaces in $\mathbb P^k$ such that every leaf is dense.

\medskip

L. Garnett  \cite{G1983} has introduced the notion of harmonic measure
for smooth foliations (without singularities) of a compact Riemannian manifold.
She studied their ergodic properties. The article by Candel \cite{C2003}
contains a recent approach to that theory. In \cite{FS2005} the authors have shown that
a ${\mathcal C}^1$ laminated set  in  $\mathbb P^2$, without singularities carries a unique harmonic current
of mass $1$ directed by the lamination. Very recently Deroin and Klepsyn \cite{DK}
 developed the theory of diffusion on transversally conformal foliations and they showed
 that there are only finitely many harmonic measures.

\medskip

For holomorphic foliations (with singularities) of $\mathbb{P}^2$ the following analogue
was proved in \cite{BS2002}. It is valid for laminations by Riemann surfaces with a small set
of singularities, see \cite{BS2002} and \cite{FS2005}.

\begin{theorem}
Let ${\mathcal F}$ be a holomorphic foliation of  $\mathbb{P}^2$. There exists
a positive current $T$ on  ${\mathbb P^2}$, of bidimension $(1,1)$ and mass $1$ which is harmonic, i.e.
$i \partial \overline{\partial} T=0.$ Moreover in any flow box $B,$ (without singular
points) the current can be
expressed as 
$$
T=\int h_\alpha [V_\alpha] d\mu(\alpha).
$$
The functions $h_\alpha$ are positive harmonic on the local leaves $V_\alpha$
and $\mu$ is a  Borel measure on the transversal.  The function $H: B \rightarrow
\mathbb R^+$, $H_{|V_\alpha}=h_\alpha$ is Borel measurable.
\end{theorem}

Observe that if ${\mathcal F}$ is defined in $B$ by a  smooth form $\omega_0$, then
$T \wedge \omega_0=0.$ We will say that the current is directed by ${\mathcal F}.$

\medskip

A theory of intersection of 
positive harmonic currents of bidegree $(1,1)$ is developed in \cite{FS2005}.
The main purpose of the present article is, using that intersection theory,  to prove:

\begin{theorem}
Let ${\mathcal F}$ be a holomorphic foliation in $\mathbb{P}^2$ without algebraic leaves. Assume that all
singular points of $\mathcal F$ are hyperbolic. Then there is a unique positive harmonic
current $T$ of mass one, directed by $\mathcal F.$ 
\end{theorem}

A consequence of Theorem 4 and of results from \cite{FS2005}
is that the foliations ${\mathcal F}$ with only hyperbolic singular points
are uniquely ergodic in a very strong sense, see Corollary 2, i.e. the current $T$ can be obtained by an averaging process on the leaves, whose limit is independent of the leaf.  We will show, in Remark 2 section 26, a similar uniqueness
result for some classes of foliations with non hyperbolic singularities.

\medskip

Observe that under the assumption of Theorem 4 there is no non zero positive closed current directed by ${\mathcal F}$, see \cite{FS2005} and Brunella  \cite{B200x}
for a general discussion of closed cycles on foliations by Riemann surfaces.

\medskip

The intersection theory of positive harmonic currents in \cite{FS2005} is valid on
compact K\"{a}hler manifolds. We just recall a few facts restricting to $\mathbb{P}^2.$

\medskip

Let $T$ be a positive harmonic current of bidegree $(1,1)$ in $\mathbb{P}^2,$
i.e. $i\partial \overline{\partial}T=0.$ Let
$\omega$ denote the standard K\"{a}hler form on $\mathbb{P}^2.$ Then $T$ can
be written as

$$
T=c\omega+\partial S+\overline{\partial S}+i \partial \overline{\partial}u$$

\noindent with $c \geq 0$ and $S$ is a $(0,1)$ form such that $S,\partial S, \overline{\partial} S$
are in $L^2$ and $u\in L^1.$ The current $\overline{\partial} S$ depends only on $T$ and is zero
only if $T$ is closed. So the quantity $\int \overline{\partial} S \wedge  \partial \overline{S} $
which we called energy, measures how far $T$ is from being closed.
The expression

$$
\int T \wedge T:= \int(c\omega+\partial S+\overline{\partial S})
\wedge (c\omega+\partial S+\overline{\partial S})
$$

\noindent makes sense and is finite. It is independent on the choice of $S.$ Moreover if $T_1$ and $T_2$ are $2$ positive
harmonic currents such that $\int T_1 \wedge T_2=0$, then $T_1$ and $T_2$ are proportional mod
$(\partial \overline{\partial} u)$. For currents directed by foliations and whose support does carry a positive closed current, then $\int T_1 \wedge T_2=0$ implies that $T_1, T_2$ are proportional, see 
\cite{FS2005} Lemma 3.10.
On the other hand the currents directed by holomorphic foliations can be expressed
in a flow box $B$ as

$$
T=\int h_\alpha [V_\alpha] d \mu(\alpha)
$$

\noindent as described in Theorem 3. It is hence possible to consider the
geometric self intersection of such currents. More precisely consider suitable automorphisms
$\Phi_\epsilon$
of $\mathbb{P}^2$ 
which are close to the identity.
For  a current
$T$ directed by a foliation ${\mathcal F}$, it is possible to define the geometric intersection
$T \wedge_g {\Phi_\epsilon}_* (T)$ as the measure on the complement of the singular points
given locally by the expression

$$
 \int \left[\sum_{p\in J_{\alpha,\beta}^\epsilon}h_\alpha(p)h^\epsilon_\beta(p) \delta_p\right] d\mu(\alpha)d\mu(\beta)\;\;\;\;\;\;\;\;\;\;\;\;\;\;\;\;\;\;\;\;\;\;(1)
$$

\noindent where $J_{\alpha,\beta}^\epsilon$ denotes the points of intersection of the plaque $L_\alpha$
and the plaque $(\Phi_\epsilon)_*L_\beta$ and $\delta_p$ denotes the Dirac mass at $p.$ It is shown in \cite{FS2005}  that $\int T_1 \wedge T_2=\lim_{\epsilon \rightarrow 0}
\int T_1 \wedge_g T_{2,\epsilon}$ (\cite{FS2005}, Lemma 19). To show that $\int T_1 \wedge T_2=0$
it is enough to count the number of points of intersection of a given plaque with perturbed plaques and
estimate the harmonic functions.
This is done in \cite{FS2005} (Theorem 6.2) when we assume that the currents $T_1,T_2$ are
supported on a minimal laminated compact set, which is transversally of class ${\mathcal C}^1.$

\medskip

Indeed the minimality hypothesis is not used and the argument there gives the following
stronger result.

\begin{theorem} 
Let ${\mathcal F}$ be a ${\mathcal C}^1$ lamination with singularities
by Riemann surfaces in $\mathbb P^2$. Assume that there is a laminated compact set $X$ without singularities.
Then there is a unique positive harmonic current $T$, of mass $1$, directed by ${\mathcal F}.$
\end{theorem}

\begin{proof}
We know there is a harmonic current $T_1$ of mass $1$, supported on $X.$ Let $T_2$ be another
such current directed by ${\mathcal F}$, but not necessarily supported by $X.$ The argument
in \cite{FS2005} Theorem 6.2 shows that $\lim_{\epsilon \rightarrow 0}\int T_1 \wedge_g T_{2,\epsilon}
=0$. Hence $\int T_1 \wedge T_2=0.$ Therefore $T_1$ and $T_2$ are proportional.
\end{proof}

We now deal with the case where the foliation is holomorphic and the current $T$ contains
in its support singular points (which are all hyperbolic).

We will prove the following more general result than Theorem 4.

\begin{theorem} (MAIN THEOREM)
Let ${\mathcal F}$ be a holomorphic foliation of
$\mathbb{P}^2$ without algebraic leaves. Let $X$ be a closed invariant set for ${\mathcal F}.$ Assume that all singular points of $X$ are hyperbolic.
Then there is a unique positive harmonic current $T$ of mass $1$, directed
by $X.$
\end{theorem}

The result is valid for a laminated set $(X,{\mathcal L},E)$ where $X \setminus E$ is a
${\mathcal C}^1$ lamination by Riemann surfaces. The set $E=\{p_1,\dots,p_\ell\}$ is a finite set
and in a neighborhood $U_j$ of every singular point $p_j$ we assume that $X \cap U_j$ is 
holomorphically equivalent
to a lamination contained in $z=Cw^{\lambda_j}, \lambda_j=a_j+ib_j, b_j \neq 0.$ 
One of  the consequences of the main theorem is Corollary 2 (section 26) which says that appropriate
weighted averages of the leaves always converge to the current $T.$ This is a strong ergodic theorem.
The uniqueness of $T$ also permits to show that $\lambda \rightarrow T_\lambda$ is continuous when $\lambda$ varies in a holomorphic family of foliations as considered in the main theorem.

\medskip

It is easy to see that $\overline{\partial}T= \overline{\tau}\wedge T,$ $\tau$ is a $(1,0)$ form along
leaves. We consider in Section 27 a metric $g_T:=\frac{i}{2} \tau \wedge \overline{\tau}$ and we show that the curvature $\kappa$ of that metric satisfies $\kappa(g_T)=-1.$ We also define a $\underline{{\mbox{finite}}}$ measure $\mu_T:=
i \tau \wedge \overline{\tau}\wedge T$. 
 We have that the measures vary continuously with the foliation. The metric  $g_T$ and the measure $\mu_T$ were introduced by
S. Frankel \cite{F1995} in the nonsingular case.

\medskip

\section{Proof of the Main Theorem}
 
Let $T$ be a harmonic current of mass $1$ supported on $X$ and directed by ${\mathcal F}.$
In a flow box

$$
T=\int h_\alpha [V_\alpha] d\mu(\alpha). \;\;\;\;\;\;\;\;\;\;(2)
$$

We have to estimate the number of intersection points of a plaque with perturbed plaques near a singularity and also to study the behaviour of the harmonic continuation $\tilde{h}_\alpha$ 
of $h_\alpha$ along a leaf near a hyperbolic singularity.

This will give us that the geometric intersection is zero and hence $\int T \wedge T=0.$
Since $T$ is arbitrary, the intersection theory of positive harmonic currents implies that $T$
is unique.

\medskip

After a change of coordinates
we do the analysis for the form $\omega_0=zdw-\lambda wdz, \lambda=a+ib, b \neq 0,$ near $(0,0).$

\medskip

In order to study positive harmonic currents near $0,$ we cover a deleted neighborhood of $0$ 
by finitely many "flow boxes" $(B_i)_{1 \leq i \leq N},$ with $0\in \overline{B_i}$ for every $i.$ Each $B_i=S_i \times \Delta$, where $S_i$ is a sector in $\mathbb C$ such that the map $\zeta \rightarrow e^\zeta$ is injective in a strip in the $\zeta-$plane $\gamma_1 < \Im \zeta <\gamma_2,$ with values in
$S_i$, $\Delta$ is a disc in $\mathbb C$, centered at $0.$ So the leaves in $B_i$ are graphs over 
all or part of $S_i$. We will consider them as the local plaques. For the sake of argument we will
use the sector $S$ given by $0<u<2 \pi.$

\medskip

The strategy for the proof is to choose a family of automorphisms $(\Phi_\epsilon)$ of 
$\mathbb P^2,$  close to the identity and to estimate the integral (1) in the flow boxes $(B_i)_{i \leq N}.$ For
that purpose we need to estimate the growth of the harmonic continuation of $h_\alpha$
along the leaves  and also the number of intersection points of a plaque
$L_\alpha,$ with perturbed plaques $L^\epsilon_\beta$. 

 Away from singularities this is just the proof given in \cite{FS2005} for a
lamination. In the present case we have to divide the phase space in many regions where
the estimates are technically different. The estimates are different close
to separatrices and in other regions. This requires a precise subdivision of a polydisc near a
singular point. We will describe the subdivision in more details after stating Theorem 7.

Consider again the foliation $zdw-\lambda wdz=0, \lambda=a+ib, b \neq 0.$
Notice that if we flip $z$ and $w$, we replace $\lambda$ by $1/\lambda=\overline{\lambda}/|\lambda|^2=a/(a^2+b^2)-ib/(a^2+b^2).$ We will assume below that the axes are chosen so that
$b>0.$ However, it is important to note that the estimates are also valid if $b<0.$
The point is that it will be seen that the case $a=1$ is a degeneracy that complicates the estimates. However if we flip coordinates, the constant $a=1$ becomes $a/(a^2+b^2)=1/(1+b^2)<1.$
We now describe a general leaf.

There are two separatrices, $(w=0), (z=0).$ Other than that a leaf $L_\alpha$
can be parametrized by

\bea
(z,w) & = & \psi_\alpha(\zeta)\\
z & = & e^{i(\zeta+(\log |\alpha|)/b)}, \zeta = u+iv\\
w & = & \alpha e^{i\lambda(\zeta+(\log |\alpha|)/b)}\\
|z| & = & e^{-v}\\
|w| & = & e^{-bu-av}\\
\eea

Notice that as we follow $z$ once counterclockwise around the origin, $u$ increases by $2\pi,$ so the absolute value
of $|w|$ decreases by the multiplicative factor of $e^{-2\pi b}$. Hence we cover all
leaves by restricting the values of $\alpha$ so that
$e^{-2\pi b} \leq |\alpha| <1.$ We observe that with the above parametrization, the intersection with the unit bidisc
of the leaf is given by $v>0, u>-av/b$ independently of $\alpha.$ In the $(u,v)$-plane this corresponds to
a sector $S=S_\lambda$ with corner at $0$ and given by $0<\theta< \arctan (-b/a)$ where the arctan
is chosen to have values in $(0,\pi).$ Let $\gamma:= \frac{\pi}{\arctan (-b/a)}$. Then
the map $\phi: \tau \rightarrow \tau^\gamma$ maps this sector to the upper half plane with
coordinates $(U,V).$\\ The fact that $\gamma >1 $ will be crucial, this is where the hyperbolicity of singularities is used.

 Let $h_\alpha$ denote the harmonic function associated to the current $T$ on the leaf $L_\alpha.$
 The local leaf clusters on both separatrices. To investigate the clustering on the $z-$ axis,
 we use a transversal $D_{z_0}:=\{(z_0,w);|w|<1\}$ for some $|z_0|=1.$ We can normalize so that
 $h_\alpha(z_0,w)=1$ where $(z_0,w)$ is the point on the local leaf with
 $e^{-2\pi b}\leq |w|<1.$ So $(z_0,w)=\psi_\alpha(\zeta_0)=\psi_\alpha(u_0+iv_0)$ with
 $v_0=0$ and $0<u_0\leq 2\pi$ determined by the equations $|z_0|=e^{-v_0}=1$
 and $e^{-2\pi b}\leq |w|=e^{-bu_0-av_0}<1.$ Let $\tilde{h}_\alpha$ denote the harmonic continuation along $L_\alpha.$
 Define $H_\alpha(\zeta):=\tilde{h}_\alpha(
 e^{i(\zeta+(\log |\alpha|)/b)},\alpha e^{i\lambda (\zeta+(\log |\alpha|)/b)})$ on $S_\lambda$.

\begin{proposition}
The harmonic function $\tilde{H}_\alpha:=H_\alpha \circ \phi^{-1}$ is the Poisson integral of its
boundary values. So in the upper half plane $\{U+iV;V>0\},$
$$
\tilde{H}_\alpha (U+iV)=\frac{1}{\pi} \int_{-\infty}^\infty \tilde{H}_\alpha(x) \frac{V}{V^2+(x-U)^2}dx
$$
[for a.e. $\alpha, d \mu$].
Moreover, 
$$
\int_{e^{-2\pi b}\leq |\alpha|<1} \int_{-\infty}^\infty \tilde{H}_\alpha(x) |x|^{\frac{1}{\gamma}-1}dx d\mu(\alpha)<\infty.
$$
\end{proposition}

\begin{proof}
Let $A_n:=\{(z_0,w); e^{-2\pi b(n+1)}\leq |w|<e^{-2\pi bn},  \; n=0,1,\dots\}. $ The holonomy map around $(z=0)$ as described above gives a map

$$
A_n \rightarrow A_{n+1}.
$$

\medskip

The transverse  masses of these sets are $\int_{A_0}
H_\alpha (\zeta_0+2\pi  n) d\mu(\alpha)=B_n(\zeta_0).$ The functions $B_n(\zeta)$ are harmonic on
$\{v>0,u>-av/b-2\pi n\}.$ 
Since the transverse mass is finite on $(z=z_0)$ and since the annuli $A_n$
are disjoint we get,

$$
\sum_{n=0}^\infty B_n(\zeta_0) <\infty. \;\;\;\;\;\;\;\;(*)
$$

 We get a similar estimate along the other separatrix.
It follows that

$$
\int_{A_0}\left(\int_{\partial S_\lambda}H_\alpha \right)d\mu(\alpha) <\infty.\;\;\;\;\;\;\;\;\;\;\;\;\;\;\;(3)
$$

We show now that for almost every $\alpha$, $\tilde{H}_\alpha(x,y)$ is equal to the Poisson integral of its restriction to $y=0.$ Every positive harmonic function on the upper half plane can
be written as a sum of a Poisson integral and $cy, c \geq 0.$ The problem is to show that $c=0.$

We consider the restriction $L'_\alpha$ of $L_\alpha$ to the bidisc $\Delta^2(0,e^{-1}).$ The leaf
$L'_\alpha$ equals $\psi_\alpha(S'_\lambda)$ where
$S'_\lambda:=\{v>1,u>-av/b+1/b\}.$ The image of this sector under $\phi$ is a domain of the
form 
$\Delta_{\lambda,\alpha}'= \{x+iy; y>\gamma_\alpha (x)\}$ where
$\gamma_\alpha$ is a continuous strictly positive function so that
$\gamma_\alpha \rightarrow +\infty$ when $|x| \rightarrow \infty.$ 
The function $B_1$ is bounded on the edges of $S'_\lambda.$ So $\tilde{B}_1:=
B_1 \circ \phi^{-1}$  is bounded on the graph
of $\gamma_\alpha$ and hence there is no term $cy, c>0$ in the canonical representation
of $\tilde{B}_1.$ The same argument is valid for the functions $\tilde{H}_\alpha$ at least
for $\mu$ almost every $\alpha.$

It follows that the representation as a Poisson integral is valid. On the other hand
estimate (3) can be read as

$$
\int_{e^{-2\pi b}\leq |\alpha|<1} \int_0^\infty H_\alpha(u) du d\mu(\alpha)<\infty\; {\mbox{and}}
$$
$$\int_{e^{-2\pi b}\leq |\alpha|<1} \int_0^\infty H_\alpha(u e^{i\arctan(-b/a)}) du d\mu(\alpha)<\infty,
$$

\noindent which gives after a change of variables the estimate on the growth of $\tilde{H}_\alpha.$

\end{proof}

\begin{corollary} Let ${\mathcal F}$ be a foliation as in Theorem 6. Then for any positive harmonic current $T$, directed by $\mathcal F,$ the transverse measure $\mu$ is diffuse.
\end{corollary}

\begin{proof}
Assume $\mu$ has an atomic part i.e. a Dirac mass at $p.$ Let $L$ be the leaf through $p.$ The restriction $T$ to $L$ is a non zero positive harmonic current. We can normalize so that the transverse measure is one. Then we have  a positive harmonic function $h$ defined on $L.$ 

If there is a flow box $B$, away from the singularities, such that $L$ crosses $B$ on infinitely many plaques on which $h$ is bounded below by a strictly positive constant, then we get a contradiction because the mass of $T$ should be finite. In any flow box the leaf must intersect in infinitely many plaques $P_j$ and the harmonic functions $h_j=h_{|P_j}$ must go uniformly to zero as $j \rightarrow \infty.$



Let $f$ denote the lifting of the harmonic function to the unit disc, so $f=h \circ \phi$ where $\phi:\Delta \rightarrow L$ is a universal covering map. 
Since $f>0$
 there is a set $S \subset \partial \Delta$ of full measure on which $f$ has nontangential limits $f(e^{i\theta}).$

\begin{lemma}
The function $f(e^{i\theta})=0$ a.e. on $S.$ 
\end{lemma}

\begin{proof}Suppose that  $f(e^{i\theta_0})>0, \theta_0\in S.$
We consider the curve $\phi(re^{i \theta_0}).$ By the above argument, it follows that this curve can only intersect finitely many plaques in any flow box away from the singular points. But if some plaque is visited infinitely many times as $r\rightarrow 1$,  we see that $h$ must be constant on this plaque, hence constant on the leaf, a contradiction. It follows that the curve converges to a singular point. 

Then it follows from \cite{FS2005}, p 991 that this only happens on a set of measure $0$ because almost every radius leaves some ball around the singularity.

\end{proof}

A consequence of the Lemma is that the function $f$ is given by the convolution of the Poisson kernel with a singular measure. This implies that the function $f$ is unbounded. Outside any given neighborhood of the singular set the function must be uniformly bounded. But then the Poisson integrals of Proposition 1 are also uniformly bounded. Hence by Proposition 1 the function is uniformly bounded everywhere, a contradiction.

\end{proof}

\begin{remark}
It is convenient in some later calculations to replace $|x|^{1/\gamma-1}$
by $(|x|+1)^{1/\gamma-1}$ in the integral of Proposition 1. By Harnack, this doesn't effect the order
of magnitude of the integral.
\end{remark}

We  decompose a leaf $L_\alpha$ into 
plaques $L_{\alpha,n}$ where $2n\pi <u< 2(n+1)\pi.$ Here $n$ is an integer. [ Note that if $a\leq 0$,
these $n$ must be positive to have a nonempty intersection with the bidisc.]
In this way $L_{\alpha,n}$ is a graph over some part of the $z-$ axis.

We let $(z,w)$ be a point in $L_{\alpha}$ parametrized by a point $(u,v).$ We write in polar coordinates, $u+iv=\rho e^{i \theta}$ with
$\rho=\sqrt{u^2+v^2}, \theta = \arctan (v/u).$ Then in the $(U,V)$ plane this point
corresponds to $U+iV=\phi(u+iv)=(u+iv)^\gamma,$

$$
U+iV= \rho^\gamma e^{i \gamma \theta}
=\rho^\gamma \cos (\gamma \theta)+i\rho^\gamma \sin (\gamma \theta).
$$

We hence get the following formula for the function $H_\alpha(u+iv):$

\begin{lemma}
$$
H_{\alpha}(u+iv)= \frac{1}{\pi} \int \tilde{H}_\alpha(x) \frac{\rho^\gamma \sin (\gamma \theta)}
{(\rho^\gamma \sin (\gamma \theta))^2+(x-\rho^\gamma \cos (\gamma \theta))^2} dx
$$
\end{lemma}

Now we write the formula for the perturbed foliation ${\mathcal F}_\epsilon
=(\Phi_\epsilon)_*\mathcal F$ where $\Phi_\epsilon$ is a family of automorphisms of $\mathbb P^2.$
We will need as in \cite{FS2005} that all our estimates stay valid when composing
$\Phi_\epsilon$ with $\Psi$ in a neighborhood of the identity in $U(3)$
(depending on $\epsilon$). We will need that $\Phi_\epsilon$ moves the singular point in a direction away from the separatrices
near all the hyperbolic points. We also need the $\Phi_\epsilon$ to have a common fixed point $p$ in the
support of $T$ and that the tangent space of the leaf through $p$ moves to first order with
$\epsilon.$ So we write in $\mathbb C^2$ 

$$
\Phi_\epsilon (z,w) =(\alpha(\epsilon),\beta(\epsilon))+(z,w)+\epsilon {\mathcal O}(z,w)
$$

\noindent with $\alpha'(0), \beta'(0) \neq 0. $ We will also need that $\lambda \neq \beta'(0)/\alpha'(0).$

Suppose that $(z,w)$ is
a point in the perturbed bidisc $\Phi_\epsilon(\Delta^2)$, not on an indicatrix of ${\mathcal F}_\epsilon.$ 
Then $\Phi_\epsilon^{-1}(z,w)$ is on some plaque $L_{\beta,m}$ with parameters $(u',v')$. 
We ignore the problem that we need $u' \neq 2\pi$ because we can also use other flow boxes. 
The original point $(z,w)$
is on a plaque $L_{\beta,m}^\epsilon$ and we get:

\begin{lemma}
$$
H_{\beta}^\epsilon(u'+iv')= \frac{1}{\pi} \int \tilde{H}_\beta^\epsilon(y) \frac{(r')^\gamma \sin (\gamma \theta')}
{((r')^\gamma \sin (\gamma \theta'))^2+(y-(r')^\gamma \cos (\gamma \theta'))^2} dy
$$
\end{lemma}

Next, for each $(\alpha, \beta, m,n,\epsilon),$ let $I_{\alpha,\beta,m,n,\epsilon}$ denote the set of points $p$
in a slightly smaller bidisc which belong to $L_{\alpha,n} \cap L^\epsilon_{\beta,m}.$ Our main technical result is the following Theorem,which says that the geometrical intersection is zero, so that the current T is unique, see section 26.

\begin{theorem}
$$\lim_{\epsilon \rightarrow 0} \int  \sum_{m,n} \sum_{p\in I_{\alpha, \beta,m,n,\epsilon}}
\tilde{h}_{\alpha,n}(p)\tilde{h}^\epsilon_{\beta,m}(p) d\mu(\alpha) d\mu(\beta)=0.$$
\end{theorem}

\begin{proof}
During the proof it will be convenient to divide up the region of integration into several pieces.
For constants $0<c< C$ and $\delta>0$, we consider the regions around one of the finitely many singular points.The regions are defined as follows:

\bea
D_1 & = & \{(z,w); |z|\leq c\epsilon, |w|\leq c \epsilon\}\\
D_2 & = & \{(z,w); |z|\leq C\epsilon, |w|\leq C\epsilon\}\setminus D_1 \\ 
D_3 & = & \{(z,w); |z|\leq \delta, |w|\leq \delta\}\setminus D_1\cup D_2.\\
\eea

By \cite{FS2005}, for any given $\delta>0$, the contribution to the integral from outside these
regions goes to zero when $\epsilon \rightarrow 0,$ this uses that the measure is diffuse.
We will subdivide the regions $D_1,D_2,D_3$ further. For most of these new subregions
the contributions go to zero with $\epsilon$. But for some of the subregions, we need $\delta$ to go to
zero for the contribution  to  go to zero. Hence in the following arguments, $\delta$ will be an unspecified small number which will later go to zero. So the way the argument works  is, in order to show that the integral becomes smaller than some given $\tau>0$ when $\epsilon\rightarrow 0,$ we first fix a very small $\delta$ and then let $\epsilon \rightarrow 0$. This is the case at the end of Sections 8, 9. We will constantly use the finiteness of the integral in Proposition 1, in order to show that the limits are zero.

\medskip

The constants $c$ and $C$ are determined by the geometry of the leaves near the singularity.
We choose $c>0$ small enough that the region $D_1$ does not contain the singular point of the
perturbed foliation. In fact we will make $c>0$ so small that the slopes of the leaves of the
perturbed foliation are almost constant on $D_1.$ The precise estimate is done in Lemma 3.

For the constant $C$ we want to make sure that the singular point of the perturbed
foliation is inside $\Delta^2(0,C|\epsilon|/2).$  So for example the choice
$C=3 \max\{|\alpha'(0)|,|\beta'(0)|\}$ will work.

\section{Proof of theorem 7 for the intersection points in $D_1$ (close to the singularity)}

\begin{lemma}
Let $\delta>0$. Then for all small enough $c$, $|\epsilon|$, the slopes of the leaves
of ${\mathcal F}_\epsilon,$ 
$dw/dz\in \Delta(\lambda \frac{\beta'(0)}{\alpha'(0)}; \delta)$ at all points in $D_1.$
\end{lemma}

\begin{proof}
Recall that $$
\Phi_\epsilon (z,w) =(\alpha(\epsilon),\beta(\epsilon))+(z,w)+\epsilon {\mathcal O}(z,w).
$$
We estimate $\omega_\epsilon$ in $D_1$.

\bea
\omega_\epsilon & = & (\Phi_\epsilon)_*(\omega_0)\\
& = &  {\mathcal O}(\epsilon^2)+
\left[(z-\alpha(\epsilon))(1+A\epsilon)+B\epsilon(w-\beta(\epsilon))\right]dw\\
& 
+ & \left[(z-\beta(\epsilon))(-\lambda+C\epsilon)+D\epsilon (z-\alpha(\epsilon))\right]dz\\
& = &  {\mathcal O}(\epsilon^2)+(z-\alpha(\epsilon))dw+(z-\beta(\epsilon))(-\lambda)dz\\
& = & (z-\alpha'(0)\epsilon+ {\mathcal O}(\epsilon^2))dw-
\lambda (w-\beta'(0)\epsilon+{\mathcal O}(\epsilon^2))
dz.\\
\eea

So,
\bea
dw/dz & = & \frac{\lambda (w-\beta'(0)\epsilon+{\mathcal O}(\epsilon^2))}
{z-\alpha'(0)\epsilon+ {\mathcal O}(\epsilon^2)}\\
& = & \lambda \frac{-\beta'(0)\epsilon + \cdots}{-\alpha'(0)\epsilon +\cdots}\\
& = & \lambda \frac{\beta'(0)}{\alpha'(0)} +\cdots\\
\eea

The Lemma follows immediately.

\end{proof}

The following lemma decribes the lamination associated to  $\omega_\epsilon$ near $D_1$
after possibly shrinking $c$ further and is an immediate consequence of Lemma 3.

\begin{lemma}
The plaques of ${\mathcal F}_\epsilon$ near $D_1$ are of the form
$w=f_\eta(z)$ where $f_\eta(\eta)=0$ and $f'_\eta\in \Delta(\lambda \frac{\beta'(0)}{\alpha'(0)}; \delta).$
\end{lemma}

To estimate the geometric wedge product we will consider three types of points in a plaque
$L^\epsilon_{\beta,m}$, namely if they are close to where the plaque crosses
the $z-$ axis (Case 1) or $w-$ axis or otherwise (Case 2). The estimates for $\tilde{h}^\epsilon_\beta$ are
fairly independent of which case we are in because of the choice of $c$, but $h_\alpha$ is very sensitive to the cases.

\medskip

We estimate the function $\tilde{h}_\beta^\epsilon$ on these plaques. First observe that the points
in $B_2:= \Delta^2((-\alpha'(0)\epsilon,-\beta'(0)\epsilon); 2c |\epsilon|)$ are mapped by $\Phi_\epsilon$ to a
region covering $D_1.$

\begin{lemma}
There is a constant $C>0$  so that if some leaf $L^\epsilon_\beta$ intersects $D_1$ for a
parameter value $u+iv$ then 
$$
\frac{1-a}{b} \log (1/|\epsilon|)-C<u<\frac{1-a}{b} \log (1/|\epsilon|)+C,
\log (1/|\epsilon|)-C<v<\log(1/|\epsilon|)+C.
$$
\end{lemma}

\begin{proof}
First recall that $z=e^{i(u+iv+(\log |\beta|/b))}.$ Hence
$|z|=e^{-v}.$ But $(z,w)\in B_2$. Hence
$$
(|\alpha'(0)|-2c)||\epsilon |< |z|=e^{-v}<(|\alpha'(0)|+2c)||\epsilon|.
$$
So $$\log| \epsilon |-C< -v< \log |\epsilon|+C$$ which gives the estimate on $v.$
To get the estimate on $u,$ we consider:
\bea
|w| & = & e^{-bu-av}\\
(|\beta'(0)|-2c)|\epsilon| & < & |w|=e^{-bu-av}<(|\beta'(0)|+2c)|\epsilon|\\
\log |\epsilon|-C' & < & -bu-av < \log |\epsilon|+C'\\
\log(1/ |\epsilon|)-C' & < & bu+av < \log (1/|\epsilon|)+C'\\
\frac{1}{b}\log(1/ |\epsilon|)-C''& < & u+av/b < \frac{1}{b}\log (1/|\epsilon|)+C''.\\
\eea

So

\bea
\frac{1}{b}\log(1/ |\epsilon|)-av/b-C'' & < & u< \frac{1}{b}\log (1/|\epsilon|)-av/b +C'',\\
\eea

\noindent hence
\bea
u & < & \frac{1}{b}\log (1/|\epsilon|)+\frac{a}{b} [-\log (1/|\epsilon|)+C] +C'' \;[a>0]\\
& < & \frac{1-a}{b} \log (1/|\epsilon|)+C'''\\
u & < & \frac{1}{b}\log (1/|\epsilon|)-\frac{a}{b} [\log (1/|\epsilon|)+C] +C'' \;[a\leq 0]\\
& < & \frac{1-a}{b} \log (1/|\epsilon|)+C'''\\
\eea

For the other estimate 

\bea
u & > & \frac{1}{b}\log (1/|\epsilon|)+\frac{a}{b} [-\log (1/|\epsilon|)-C] -C'' \;[a>0]\\
& > & \frac{1-a}{b} \log (1/|\epsilon|)-C'''\\
u & > & \frac{1}{b}\log (1/|\epsilon|)-\frac{a}{b} [\log (1/|\epsilon|)-C] -C'' \;[a\leq 0]\\
& > & \frac{1-a}{b} \log (1/|\epsilon|)-C'''.\\
\eea

\end{proof}

In what follows we use the notation $A \sim B$ to mean that there is a constant $L$ so that
$\frac{1}{L}A\leq B \leq LA$ and $L$ is chosen independent of the other parameters.
Also $A \lesssim   B$ means similarly that there is a constant $L$ so that $A \leq LB.$

Next we estimate the value of $\tilde{h}_\beta^\epsilon$ for a point $(u,v)$ as in the previous Lemma.
Let $\theta,  \tan \theta =v/u$  be the argument. By  Lemma  5, it follows that for
all small $\epsilon$, $\tan \theta\sim b/(1-a)\neq b/(-a)$ so that the angle $\theta$ is uniformly inside the
sector $S_\lambda$ for all small $\epsilon.$ It follows that $\gamma \theta$ is strictly
inside a sector $0<s<\gamma \theta<\pi-s<\pi$ for some fixed $s>0$ which only depends on $\lambda$ and is independent of all other choices, and for all small enough $\epsilon.$
This implies that $\sin \gamma \theta > k>0$ uniformly. As in Lemma 2, for a point in $D_1$, this allows us to estimate the kernel
for $H^\epsilon_{\beta}(u+iv)$:

\begin{lemma}
Suppose $(u+iv)$ is such that the corresponding point on the leaf $L^\epsilon_\beta$
is in $D_1,$ then if $|y|<2 (\log (1/|\epsilon|))^\gamma,$
$$
 \frac{(r)^\gamma \sin (\gamma \theta)}
{((r)^\gamma \sin (\gamma \theta))^2+(y-(r)^\gamma \cos (\gamma \theta))^2} \sim\frac{1}{(
\log (1/|\epsilon|))^\gamma}
$$
On the other hand if $|y|\geq 2 (\log (1/|\epsilon|))^\gamma$ then
$$
 \frac{(r)^\gamma \sin (\gamma \theta)}
{((r)^\gamma \sin (\gamma \theta))^2+(y-(r)^\gamma \cos (\gamma \theta))^2} \sim\frac{(\log (1/|\epsilon|))^\gamma}{y^2}
$$
\end{lemma}

Hence we get using Lemma 2:

\begin{lemma}
We have the following  estimate of $H_\beta^\epsilon$ for points in $D_1$:
$$
H_\beta^\epsilon \sim  \frac{1}{(\log (1/|\epsilon|))^\gamma}
\int_{|y|<2 (\log (1/|\epsilon|))^\gamma} \tilde{H}_\beta(y)dy+
(\log (1/|\epsilon|))^\gamma \int_{|y|\geq 2 (\log (1/|\epsilon|))^\gamma}\frac{\tilde{H}_\beta(y)}{y^2}dy
$$
\end{lemma}

Next we fix $\alpha, \beta$ and plaques $L_{\alpha,n}, L^\epsilon_{\beta,m}$ and assume they
intersect in $D_1.$ By Lemma 5, there are conditions on $m$ for this to happen, namely:

\bea
2m\pi & < & u' < 2(m+1)\pi\\
\frac{1-a}{b} \log (1/|\epsilon|)-C & < & u' < \frac{1-a}{b} \log (1/|\epsilon|)+C \\
\eea

So 

\bea
\frac{1-a}{b} \log (1/|\epsilon|)-C & < & 2(m+1) \pi < \frac{1-a}{b} \log (1/|\epsilon|)+C+2\pi\\
\frac{1-a}{b} \log (1/|\epsilon|)-C-2\pi & < & 2m \pi < \frac{1-a}{b} \log (1/|\epsilon|)+C.\\
\eea

We pick a plaque $L^\epsilon_{\beta,m}$ with an intersection point in $D_1$.
Then this plaque is of the form $w=f(z)=f_\eta(z)$ where $f_\eta(\eta)=0$ and
$f'$ is as in Lemma 4. i.e. close to $\lambda \frac{\beta'(0)}{\alpha'(0)}.$
Next consider a plaque $L_{\alpha,n}$

\bea
z & = & e^{i(u+(\log |\alpha|/b))-v}\\
w & = & \alpha e^{i \lambda (\zeta+ (\log |\alpha|)/b)}\\
2n\pi & <  & u < 2(n+1)\pi.\\
|w| & = & e^{-bu-av}.\\
\eea

We estimate the location of the intersection points.

 Case 1: $|z-\eta|<d|\eta|, 0<d<<1.$ The constant d will be chosen small enough, in order to satisfy an inequality at the end of the proof of Lemma 8.
 
 \medskip
 
 We estimate the parameter values $(u,v)$ for $L_{\alpha,n}.$\\
 
 Since $|\eta|(1-d)<|z|=e^{-v} < |\eta|(1+d),$ 
 $\log (1/|\eta|)-2d<v< \log (1/|\eta|)+2d.$ Note that also, for the point $(z,w)$ to be on
 $L_{\beta,m}^\epsilon$ with $|z-\eta|<d|\eta|$ we must have that $|w|< 2|\lambda| \frac{|\beta'(0)|}{|\alpha'(0)|}d |\eta|.$

 \begin{lemma} For $(z,w)$ to be an intersection point between $L_{\alpha,n}$ and $L_\beta^\epsilon$
 in $D_1$ with $|z-\eta|<d|\eta|$, we must have\\
 (i) $2n\pi< u< 2(n+1)\pi$\\
 (ii) $2n\pi > \frac{1-a}{b} \log (1/|\eta|)-C$\\
 (iii)   $\log (1/|\eta|)-2d<v< \log (1/|\eta|)+2d.$\\
 Moreover there is at most one such intersection point.
  \end{lemma}
  
  \begin{proof}
 We have already proved (iii) and (i) is given. To prove (ii):
 \bea
 |w| & = & e^{-bu-av}\\
 & < & 2 |\lambda|\frac{|\beta'(0)|}{|\alpha'(0)|}d |\eta|.\\
 \eea
 
 So,
 \bea
 -bu-av & < & \log |\eta|+C.\\
 \eea
 
 Using the estimate on $v$ we get
 \bea
 u & > & (-a/b) v- (\log |\eta| )/b-C/b\\
 u & > & (-a/b) \log (1/|\eta|) - (\log |\eta| )/b-C'/b\\
  u & > & ((1-a)/b) \log (1/|\eta|)-C''\\
  \eea
  \noindent where $C',C''$ are absolute constants.
  
  To prove the last part, notice that the slope of $L_\beta^\epsilon$ is about $\lambda$
  while the slope of $L_\alpha$ is $\lambda w/z$ so is at most $ |\lambda| (2|\lambda|
  \frac{|\beta'(0)|}{|\alpha'(0)|}d||\eta|)
  /(|\eta|(1-d))<< |\lambda|$ if we just make $d$ small enough.
  
  \end{proof}
  
  \begin{lemma}
  We estimate the value of $H_{\alpha}$ at intersection points between
  $L_{\alpha,n}$ and $L_{\beta}^\epsilon$ in $D_1$ with $|z-\eta|<d|\eta|.$ We have two cases:\\
  (i) $\frac{1-a}{b} \log (1/|\eta|)-C <2\pi n< C \log (1/|\eta|).$ Then we have
    \bea
   H_{\alpha}(u+iv)  & \sim &
   \int_{|x|< 2(\log (1/|\eta|))^\gamma} \frac{\tilde{H}_\alpha(x)}{(\log (1/|\eta|))^{\gamma}} +
    \int_{|x|> 2(\log (1/|\eta|))^\gamma} \tilde{H}_\alpha(x) \frac{(\log (1/|\eta|))^\gamma}{ x^2}\\
   \eea
   In the next case,
    (ii) $2\pi n \geq C \log (1/|\eta|),$  we then have: $U+iV\sim n^\gamma+in^{\gamma-1}\log (1/|\eta|)$ and
     \bea
   H_{\alpha}(u+iv) & \sim & \int_{|x-U|<n^{\gamma-1}\log (1/|\eta|)} \frac{\tilde{H}_\alpha(x)}
   {n^{\gamma-1} \log (1/|\eta|)}dx\\
   & + & \int_{|x-U|>n^{\gamma-1}\log (1/|\eta|)} \frac{\tilde{H}_\alpha(x) n^{\gamma-1} \log (1/|\eta|)}
   {(x-U)^2}dx\\
   \eea
  \end{lemma}
  
  \begin{proof}
  Case (i): We use that  $\sin (\gamma \theta)$ is bounded below by a strictly positive constant.
      Case (ii) is clear.

  \end{proof}
  
  Case 2: 
  Our next step is to discuss intersection points of $L_{\alpha,n}$ and $L^\epsilon_{\beta,m}$ in $D_1$
  for which $|z-\eta|>d|\eta|.$ Note that $L^\epsilon_{\beta,m}$ intersects the $w-$ axis close to
  $(0,-\lambda \frac{\beta'(0)}{\alpha'(0)}\eta)$ and the above argument applies as well to the region
  $|w+\lambda\frac{\beta'(0)}{\alpha'(0)} \eta|<d|\eta|.$ Hence we only need to consider intersections of
  $L_{\alpha,n}$ and $L^\epsilon_{\beta,m}$ when $|w+\lambda \frac{\beta'(0)}{\alpha'(0)}\eta|>d|\eta|$ and also $|z-\eta|>d|\eta|$,
  call this set $S'$.
  
  \bigskip
  
  Note: This is the place in the argument where we will assume that $a\neq 1.$

  \bigskip
  
  Since we are excluding the points near where $L^\epsilon_{\beta,m}$ crosses the two axes, we have the following
  estimate on points in $L^\epsilon_{\beta,m}:$
  For some fixed constant $R>1$ we have that
  
  $$\frac{1}{R} |w|<|z|<R |w|$$
  
  \noindent for points in $S'.$
  
  Hence
  
  \bea
  \frac{1}{R}e^{-av-bu} & < & e^{-v} <R e^{-av-bu}\\
  -av-bu-\log R & < & -v<-av-bu+\log R\\
  bu-\log R & < & (1-a)v < bu+\log R\\
    2n\pi b-\log R &  < &  (1-a)v< 2(n+1)\pi b+\log R,\\
  2n\pi -C &  < &  \frac{1-a}{b}v< 2n\pi+C,\\
  2nb\pi /(1-a)-C' & < & v < 2nb\pi/(1-a)+C.'\\
   \eea
  
\end{proof}

\begin{lemma} (Intersection Lemma)
There is a constant $N>1$ so that if we cover
the rectangle $2n\pi<u<(2n+1)\pi, 2nb\pi /(1-a)-C'  < v < 2nb\pi/(1-a)+C'$
with $N$ equal squares, then there are at most two intersection points
in each square.
\end{lemma}

\begin{proof}
We choose $N$ so that in each square, the slope of $L_{\alpha,n}$ is almost constant and will produce
at most one intersection point. The exception is when the slope is close to $\lambda \frac{\beta'(0)}{\alpha'(0}.$
Then there might be a tangency between $L_{\alpha,n}$ and $L_\beta.$ Hence there
might be two or more intersection points counted with multiplicity. We will show there
are at most $2.$ Note that the slope $S$ of $L_{\alpha,n}$ is given by the quotient
$\lambda w/z$.

\bea
dw/dz & = & \lambda w/z\\
& = & \lambda \frac{\alpha e^{i \lambda (\zeta+(\log |\alpha|/b))}}{e^{i(u+(\log |\alpha|/b))-v}}\\
& = & \frac{\lambda \alpha e^{((\log |\alpha|)/b)(-b+ia)}}{e^{i (\log |\alpha|)/b}}
 \frac{e^{i\lambda \zeta}}{e^{i\zeta}}.\\
 \eea
 So,
 \bea
 S & = & C e^{i(\lambda-1)\zeta}\;\\
 \eea
 and
 \bea
\frac { \partial S}{\partial \zeta} & = & i (\lambda-1)S\\
& \sim & i(\lambda-1)\lambda \frac{\beta'(0)}{\alpha'(0)}\\
& \sim & 1.\\
 \eea

This says that the slope of $L_{\alpha,n}$ near intersection points vary very rapidly, while
we also see from Lemma 4 that the slope of $L^\epsilon_{\beta,m}$ varies slowly. This implies that near tangential intersection
points there are at most two of them.

\end{proof}

We estimate the value of $H_{\alpha}$ at points $p$ where $L_{\alpha,n}$ and $L^\epsilon_{\beta,m}$ intersect
in $D_1$ away from the axes ($|z-\eta|>d|\eta|, |w+\lambda \frac{\beta'(0)}{\alpha'(0)}\eta|>d|\eta|$).

\begin{lemma} For the intersection point to be in $D_1$ we need \\
$|n|> [|1-a| \log (1/|\epsilon|)]/[2\pi b]  -C_1.$ Then
\bea
H_{\alpha}(p) & \sim & \int_{|x|<C_2|n|^\gamma} \frac{\tilde{H}_\alpha(x)dx}{|n|^\gamma}+
\int_{|x|>C_2|n|^\gamma} \frac{\tilde{H}_\alpha(x) |n|^\gamma}{x^2}dx\\
\eea
\end{lemma}

\begin{proof}
For the first estimate, recall that $|z|=e^{-v}<c|\epsilon|$ and that 
$  2nb\pi /(1-a)-C'  <  v < 2nb\pi/(1-a)+C'$.
For the integral estimate we see that $(u+iv)^\gamma=U+iV$ with $V \sim |n|^\gamma$ and
$|U| < \sim  |n|^\gamma.$ Then the estimate is immediate from the Poisson kernel.

\end{proof}

We finish the estimate for $D_1.$

\begin{theorem}
The contribution to the geometric wedge product of $T$ and $T_\epsilon$ from
intersection points in $D_1$ goes to zero when $\epsilon \rightarrow 0.$
\end{theorem}

\begin{proof}
Let $I=I_\epsilon$ consist of all intersection points $p$ in $D_1$. They
are labeled $p=p_{\alpha,\beta,n,m,\ell}$ if they belong to the plaques
$L_{\alpha,n}, L^\epsilon_{\beta,m}$ and $\ell$ lists them (with multiplicity) if there are more than one.
By Lemma 5, 
$$
\frac{(1-a) \log (1/|\epsilon|)}{2\pi b}-C<m<\frac{(1-a) \log (1/|\epsilon|)}{2\pi b}+C
$$
so in particular there are at most finitely many values of $m$ and there is a uniform
upper bound on the number of them. We can hence restrict to one fixed value of $m.$
Next recall that from Lemma 7 we have the estimate on the value of $H_\beta^\epsilon$
at each intersection point:

\bea 
H_\beta^\epsilon(p) & \sim & \frac{1}{(\log (1/|\epsilon|))^\gamma}
\int_{|y|<2(\log (1/|\epsilon|))^\gamma} \tilde{H}_\beta(y)dy\\
& + & 
(\log (1/|\epsilon|))^\gamma \int_{|y|>2(\log (1/|\epsilon|))^\gamma}
\frac{\tilde{H}_\beta(y)}{y^2}dy.\\
\eea

By Lemmas 8 and 10 there is at most a uniformly bounded number of
intersection points with $L_{\alpha,n}$ and $L_{\beta,m}^\epsilon$ in $D_1$.
Hence when we estimate the geometric wedge product we can factor out
the contribution from $\beta$ and we get an upper bound of

\bea
\int \left(\sum_p H_\beta^\epsilon \right)d\mu(\beta) & \lesssim & 
 \frac{1}{(\log (1/|\epsilon|))^\gamma}
\int_{|y|<2(\log (1/|\epsilon|))^\gamma} \tilde{H}_\beta(y)dy d\mu(\beta)\\
& + & 
(\log (1/|\epsilon|))^\gamma \int_{|y|>2(\log (1/|\epsilon|))^\gamma}
\frac{\tilde{H}_\beta(y)}{y^2}dyd\mu(\beta).\\
\eea

We collect a few equalities that will be used repeatedly in Lebesgue dominated convergence Theorem.

\begin{lemma}We have the following integral estimates.
\bea
(I) & & 
\frac{1}{(\log (1/|\epsilon|))^\gamma}\int_{|y|<2(\log (1/|\epsilon|))^\gamma}\tilde{H}_\beta(y)dy\\
& = & 
 \frac{1}{(\log (1/|\epsilon|))^\gamma}\int_{|y|<2(\log (1/|\epsilon|))^\gamma}\tilde{H}_\beta(y)
 |y|^{1/\gamma-1}|y|^{1-1/\gamma}dy\\
& \sim & 
\int_{|y|<2(\log (1/|\epsilon|))^\gamma}\tilde{H}_\beta(y)|y|^{1/\gamma-1}\frac{|y|}
{(\log (1/|\epsilon|))^\gamma} \frac{1}{(|y|+1)^{1/\gamma}}dy\\
\eea

\bea
(II)  & & 
(\log (1/|\epsilon|))^\gamma \int_{|y|>
2(\log (1/|\epsilon|))^\gamma} \frac{\tilde{H}_\beta(y)}{y^2}dy \\
& = & \int_{|y|>2(\log (1/|\epsilon|))^\gamma}\tilde{H}_\beta(y)
|y|^{1/\gamma-1}\frac{(\log (1/|\epsilon|))^\gamma}{|y|}|y|^{-1/\gamma}dy\\
\eea

\bea
(III) & & {\mbox {If $U \sim n^\gamma$}}\\
& & \int_{|x-U|<n^{\gamma-1}\log(1/|\eta|)}\frac{\tilde{H}_\alpha(x)}{n^{\gamma-1}\log(1/|\eta|)}dx\\
& \sim & \int_{|x-U|<n^{\gamma-1}\log(1/|\eta|)}\tilde{H}_\alpha(x)|x|^{1/\gamma-1}\frac{dx}{
\log(1/|\eta|)}\\
\eea
\end{lemma}
\end{proof}

We want to show that
$$\int \sum_p H_\beta^\epsilon d\mu(\beta) \rightarrow 0.$$

We use the a priori bound that we found and estimates (I) and (II) in the previous Lemma.

The Lebesgue dominated convergence theorem and Proposition 1 gives the result.

End of the proof of Theorem 8: After integrating with respect to  $\mu$ and using (I) and (II) of Lemma 12 we can use the dominated convergence theorem. 
We estimate the value of $H_{\alpha}$ at one of the intersection points $p\in D_1.$
{}From Lemma 1 we have:

$$
H_{\alpha}(p) =\frac{1}{\pi} \int \tilde{H}_\alpha(x) \frac{\rho^\gamma \sin (\gamma \theta)}{(\rho^\gamma \sin (\gamma \theta))^2+(x-\rho^\gamma \cos (\gamma \theta))^2}dx
$$

Case (i): $|z-\eta|<d|\eta|, |n|<C \log (1/|\eta|).$
By Lemma 8 it follows that $V =\rho^\gamma \sin(\gamma \theta) \sim (\log (1/|\eta|))^\gamma$
and $|U| <\sim (\log (1/|\eta|))^\gamma.$
\bea
H_{\alpha}(p) & \sim &   \int _{|x|< C(\log (1/|\eta|))^\gamma} \frac{ \tilde{H}_\alpha(x)dx }{(\log (1/|\eta|))^\gamma}\\
& + &  \int_{|x|>C(\log (1/|\eta|))^\gamma} \tilde{H}_\alpha(x) \frac{(\log (1/|\eta|))^\gamma}{x^2}dx.\\
\eea

So adding up we get

\bea
\sum_{|n|<\log (1/|\eta|)} h_{\alpha,n}(p_n)
& \sim & 
 \int _{|x|< C(\log (1/|\eta|)^\gamma} \frac{\tilde{H}_\alpha(x)dx }{(\log (1/|\eta|))^{\gamma-1}}\\
& + &  \int_{|x|>C(\log (1/|\eta|))^\gamma} \tilde{H}_\alpha(x) \frac{(\log (1/|\eta|))^{\gamma+1}}{x^2}dx\\
& \sim & 
 \int _{|x|< C(\log (1/|\eta|))^\gamma} \tilde{H}_\alpha(x) |x|^{1/\gamma-1}
 \left(\frac{|x|}{(\log (1/|\eta|))^\gamma}\right)^{1-1/\gamma}dx\\
& + &  \int _{|x|> C(\log (1/|\eta|))^\gamma} \tilde{H}_\alpha(x) |x|^{1/\gamma-1}
 \left(\frac{(\log (1/|\eta|))^\gamma}{|x|}\right)^{1+1/\gamma}dx\\
& & \int \sum_{|n|<C\log (1/|\eta|)} H_{\alpha,n}(p_n)d\mu(\alpha)dx. \\
\eea

Integrating with respect to $\mu$ we get that  $\sum_{|n|<\log (1/|\eta|)} h_{\alpha,n}(p_n)
\rightarrow 0$ as $\epsilon \rightarrow 0$ since $|\eta|<\epsilon.$ We use again the estimates in Lemma 12 and Proposition 1.

Case (ii): $|z-\eta|<d|\eta|, |n|> C\log (1/|\eta|).$
Then by Lemma 8, $n>0$ and we have $U_n \sim n^\gamma, V \sim n^{\gamma-1} \log (1/|\eta|).$
{}From Lemma 9 we have:

\bea
H_{\alpha}(p) & \sim & \int_{|x-U|<n^{\gamma-1}\log (1/|\eta|)}
\frac{\tilde{H}_\alpha(x)}{n^{\gamma-1} \log (1/|\eta|)}dx\\
& + & \int_{|x-U|>n^{\gamma-1} \log (1/|\eta|)} \tilde{H}_\alpha(x) \frac{n^{\gamma-1}\log (1/|\eta|)}
{|x-U|^2}dx.\\
\eea

So,
\bea
H_{\alpha}(p) & \sim & \int_{|x-U|<n^{\gamma-1}\log (1/|\eta|)}
\frac{\tilde{H}_\alpha(x)}{\log (1/|\eta|)}x^{1/\gamma-1}dx\\
& + & \int_{|x-U|>n^{\gamma-1} \log (1/|\eta|)} \tilde{H}_\alpha(x)
 \frac{n^{\gamma-1}\log (1/|\eta|)}
{|x-U|^2}dx.\\
\eea

Therefore,
\bea
\sum_{n>C \log (1/|\eta|)}H_{\alpha,n}(p) & \sim & \sum_{n>C \log (1/|\eta|)}\int_{|x-U_n|<n^{\gamma-1}\log (1/|\eta|)}
\frac{\tilde{H}_\alpha(x)}{\log (1/|\eta|)}x^{1/\gamma-1}dx\\
& + & \sum_{n>C \log (1/|\eta|)}\int_{|x-U_n|>n^{\gamma-1} \log (1/|\eta|)} \tilde{H}_\alpha(x)
 \frac{n^{\gamma-1}\log (1/|\eta|)}
{|x-U_n|^2}dx\\
& = & I+II\\
\eea

We are going to estimate $I$ and $II$ separately.
Note that for a given $x$, the number of integers $n$ for which $U_n-n^{\gamma-1}\log (1/|\eta|)<x<
U_n+n^{\gamma-1}\log (1/|\eta|)$ is bounded above by a multiple of $\log (1/|\eta|).$ It follows that
$I <\sim \int_{(\log(1/|\eta|))^\gamma/C}^\infty \tilde{H}_\alpha(x)x^{1/\gamma-1}dx.$ This
contribution goes to zero as $|\epsilon| \rightarrow 0$ since $|\eta|<|\epsilon|.$

To study II we estimate $U_n$ more precisely. We have
\bea
2n\pi & < & u_n < 2(n+1)\pi\\
\log (1/|\eta|) -2d & <  & v_n  < \log (1/|\eta|)+2d,\\
\eea
and
\bea
(u_n+iv_n)^\gamma & = & u_n^\gamma (1+iv_n/u_n)^\gamma\\
& = & u_n^\gamma +\gamma u_n^{\gamma-1}iv_n+E_n,\\
\eea
with
\bea
|E_n| & \lesssim & u_n^\gamma (v_n/u_n)^2\\
& \sim & n^{\gamma-2} (\log (1/|\eta|))^2.\\
\eea
 
 Hence $|U_n-(2n\pi)^\gamma| << n^{\gamma-1} \log (1/|\eta|).$ We can hence replace $U_n$
 by $(2n\pi)^\gamma$ in $II $ without changing the order of magnitude of the expression.
 We divide $II $ into pieces $II_A,II_B,II_C$. In $II_A$, $x$ is such that $n>C \log (1/|\eta|).$ In
 $II_B$, $n$ has a range of the form $n>x^{1/\gamma}+r(x) \log (1/|\eta|), r(x) \sim 1$ and
 in $II_C$, $C \log (1/|\eta|) <n<x^{1/\gamma} -s(x)\log (1/|\eta|), s(x) \sim 1.$ So

\bea
II & = & II_A+II_B+II_C\\
II_A & = & \int_{x=-\infty}^{C_1 (\log (1/|\eta|))^\gamma}\sum_{n>C \log (1/|\eta|)}
\tilde{H}_\alpha(x) \frac{n^{\gamma-1}\log (1/|\eta|)}
{|x-n^\gamma|^2}dx\\
& \sim & \int_{|x|<C_1 (\log (1/|\eta|))^\gamma}
\tilde{H}_\alpha(x) \frac{\log (1/|\eta|)}
{[10 \log (1/|\eta)|]^\gamma-x}dx\\ 
& + & \int_{x=-\infty}^{-C_1 (\log (1/|\eta|))^\gamma}
\tilde{H}_\alpha(x) \frac{\log (1/|\eta|)}
{[10 \log (1/|\eta)|]^\gamma-x}dx\\
&  \sim & \int_{|x|<C_1 (\log (1/|\eta|))^\gamma}
\tilde{H}_\alpha(x) \frac{1}
{[ \log (1/|\eta)|]^{\gamma-1}}dx\\ 
& + & \int_{x=-\infty}^{-C_1 (\log (1/|\eta|))^\gamma}
\tilde{H}_\alpha(x) \frac{\log (1/|\eta|)}
{|x|}dx\\
&  \sim & \int_{|x|<C_1 (\log (1/|\eta|))^\gamma}
\tilde{H}_\alpha(x)|x|^{1/\gamma-1} \left(\frac{|x|}
{[ \log (1/|\eta|)]^{\gamma}}\right)^{1-1/\gamma}dx\\ 
& + & \int_{x=-\infty}^{-C_1 (\log (1/|\eta|))^\gamma}
\tilde{H}_\alpha(x)|x|^{1/\gamma-1} \left(\frac{(\log (1/|\eta|))^\gamma}
{|x|}\right)^{1/\gamma}dx.\\
\eea
For
\bea
II_B & \sim & \int_{x=C_1 (\log (1/|\eta|))^\gamma}^\infty \sum_{n>x^{1/\gamma}+r(x) \log (1/|\eta|)}
\tilde{H}_\alpha(x) \frac{n^{\gamma-1}\log (1/|\eta|)}
{|x-n^\gamma|^2}dx\\
 & \sim & \int_{x=C_1 (\log (1/|\eta|))^\gamma}^\infty
\tilde{H}_\alpha(x) |x|^{1/\gamma-1}dx,\\
\eea
and
\bea
II_C & \sim & \int_{x=C_2 (\log (1/|\eta|))^\gamma}^\infty \sum_{C \log (1/|\eta|)<n<x^{1/\gamma}-s(x) \log (1/|\eta|)}
\tilde{H}_\alpha(x) \frac{n^{\gamma-1}\log (1/|\eta|)}
{|x-n^\gamma|^2}dx\\
II_C & \sim & \int_{x=C_2 (\log (1/|\eta|))^\gamma}^\infty 
\tilde{H}_\alpha(x) x^{1/\gamma-1}dx.\\
\eea

Integrating with respect to $\mu$ we get
\bea
& & \int II d\mu(\alpha)\\
& \sim & \int_\alpha II_A d\mu(\alpha)
+  \int_\alpha II_B d\mu(\alpha)
+\int_\alpha II_C d\mu(\alpha)\\
&  \lesssim & \int_\alpha \int_{|x|>(\log (1/|\eta|))^\gamma}\tilde{H}_\alpha(x) |x|^{1/\gamma-1}dxd\mu(\alpha)\\
& + &  \int_\alpha \int_{|x|<C_1(\log (1/|\eta|))^\gamma} \tilde{H}_\alpha(x)|x|^{1/\gamma-1}
\left( \frac{|x|}{(\log (1/|\eta|))^\gamma}\right)^{1-1/\gamma} dx d\mu(\alpha),\\
\eea
\noindent which tends to $0$ as $\eta\rightarrow 0.$ (Recall that $|\eta|<\epsilon.$)

Case (iii):
$|w+\lambda \eta|<d|\eta|$. This case is symmetric to cases (i) and (ii), so done.

Case (iv):
$|z|,|w|<c|\epsilon|;\; |z-\eta|, |w+\lambda \eta|>d|\eta|.$
We recall the estimate of $H_{\alpha,n}(p)$ at intersection points from Lemma 11.
The contribution $W$ to the geometric wedge product is:

\bea
\int_\alpha \left[\sum_{|n|>[|1-a|\log (1/|\epsilon|)]/[2\pi b]-C}
\int_{|x|<2|n|^\gamma}\frac{\tilde{H}_\alpha(x)dx}{|n|^\gamma}
+\int_{|x|>2|n|^\gamma}\frac{\tilde{H}_\alpha(x)|n|^\gamma}{x^2}dx \right]d\mu(\alpha).\\
\eea

We divide the first integral into two pieces, so $W=W_A+W_B+W_C.$ We get

\bea
& & W_A \\
& \sim & \int_\alpha \left[
\int_{|x|<2|[|1-a|\log (1/|\epsilon|)]/[2\pi b]-C|^\gamma}\sum_{|n|>[|1-a|\log (1/|\epsilon|)]/[2\pi b]-C}\frac{\tilde{H}_\alpha(x)dx}{|n|^\gamma}
 \right]d\mu(\alpha)\\
& \sim &  \int_\alpha \left[
\int_{|x|<2|[|1-a|\log (1/|\epsilon|)]/[2\pi b]-C|^\gamma}\frac{\tilde{H}_\alpha(x)dx}{(\log (1/|\epsilon|))^{\gamma-1}}
 \right]d\mu(\alpha)\\
& \sim &   \int_\alpha \left[
\int_{|x|<2|[|1-a|\log (1/|\epsilon|)]/[2\pi b]-C|^\gamma}\tilde{H}_\alpha(x)|x|^{1/\gamma-1}
\left(\frac{|x|}{(\log (1/|\epsilon|))^\gamma}\right)^{1-1/\gamma}dx
 \right]d\mu(\alpha)\\
 &\rightarrow & 0\; {\mbox{as}}\; \epsilon \rightarrow 0.\\
 \eea
 
 For $W_B$ we have:
 \bea
 W_B & \sim & \int_\alpha \left[
\int_{|x|>2|[|1-a|\log (1/|\epsilon|)]/[2\pi b]-C|^\gamma}\sum_{(|x|/2)^{1/\gamma}}^\infty \frac{\tilde{H}_\alpha(x)dx}{|n|^\gamma}
 \right]d\mu(\alpha)\\
  W_B & \sim & \int_\alpha \left[
\int_{|x|>2|[|1-a|\log (1/|\epsilon|)]/[2\pi b]-C|^\gamma}\tilde{H}_\alpha(x)|x|^{1/\gamma-1}dx
 \right]d\mu(\alpha)\\
&  \rightarrow & 0,\\
\eea
again by Proposition 1.
\bea
&&W_C\\
 & \sim & 
\int_\alpha \left[
\int_{|x|>2|\frac{|1-a|\log (1/|\epsilon|)}{2\pi b}-C|^\gamma}\sum_{|n|=[|1-a|\log (1/|\epsilon|)]/[2\pi b]-C}^{(|x|/2)^{1/\gamma}}\frac{\tilde{H}_\alpha(x)|n|^\gamma dx}{x^2}
 \right]d\mu(\alpha)\\
 & \sim & 
\int_\alpha \left[
\int_{|x|>2|[|1-a|\log (1/|\epsilon|)]/[2\pi b]-C|^\gamma}\tilde{H}_\alpha(x)|x|^{1/\gamma-1}dx
 \right]d\mu(\alpha),\\
 \eea
 
 \noindent and hence $W_C\rightarrow 0.$
 
Now we have finished the part of the proof of Theorem 8 where we consider
intersection points in $D_1=\{|z|,|w|<c|\epsilon|\}.$

\section{Proof of Theorem 7 for intersection points in $D_2 \subset \Delta^2(0,C|\epsilon|)$
close to the separatrices.}

We split $D_2$ into regions $A'$ and $B$ where $A'$ denotes points close to the separatrices
and $B$ denotes the rest. Then $A'$ has 2 pieces. It suffices to consider one,
$A=\{(z,w); c|\epsilon|<|z|<C|\epsilon|, |w|<r |\epsilon|\}$ where $0<r<c$ depends on the choice of $C.$

We consider intersection points of $L_{\alpha,n}$ and $L^\epsilon_{\beta,m}$ in $A.$
We parametrize $L_\alpha$ with $(u+iv)$ and $L^\epsilon_\beta$ with $u'+iv'$. In $A$ we have for $L_{\alpha,n}:$

\bea
\log (1/|\epsilon|)-C & <  &  v < \log (1/|\epsilon)+C\\
|w| & = & e^{-bu-av}< r |\epsilon|\\
bu+av & > & \log (1/r)+ \log (1/|\epsilon|)\\
u & > & \frac{1-a}{b} \log (1/|\epsilon|)-C\\
n & > & \frac{1-a}{2\pi b} \log (1/|\epsilon|)-C.\\
\eea

For $L^\epsilon_{\beta,m}$ we have in $A:$
\bea
L^\epsilon_{\beta,m}: & & |z'|< C |\epsilon|\\
v' & > & \log (1/|\epsilon|)-C\\
|\epsilon|(1-2c) & < & |w'| <|\epsilon|(1+2c)\\
\log (1/|\epsilon|)-2c & < & bu'+av' < \log (1/|\epsilon|)+2c\\
\frac{1}{b}\log (1/|\epsilon|)-2c/b-av'/b & < & u' < \frac{1}{b} \log (1/|\epsilon|)-av'/b+2c/b.\\
\eea

 The m is estimated later and they depend on which case we are in, $a= 0$ or not.

\begin{lemma}
If $a \neq 0,$ there is an integer $N$ so that for small  $r$, there is at most $N$ intersection points between
any pair $L_{\alpha,n}$ and $L^\epsilon_{\beta,m}.$
\end{lemma}

\begin{proof}
This follows from considering the slopes of the plaques, given by the  forms $\omega, \omega_\epsilon.$
Namely the slope of the $L_{\alpha,n}$ is very small and the slope of
$L^\epsilon_{\beta,m}$ has close to constant larger modulus and close to constant argument on each of $N$ small
squares where there might be an intersection.
\end{proof}

Next we estimate $h_{\alpha,n}$ at an intersection point.

Case (i): $n< \log (1/|\epsilon|):$

In $U,V$ coordinates we have  $V \sim (\log (1/|\epsilon|))^\gamma, |U| \lesssim (\log (1/|\epsilon|))^\gamma.$

Using the expression as a Poisson integral we get:

\bea
h_{\alpha,n}(p) & \sim & \int_{|x|< C(\log (1/|\epsilon|))^\gamma}H_\alpha(x) \frac{1}{(\log (1/|\epsilon|))^\gamma}dx\\ 
& + &  \int_{|x|>C (\log (1/|\epsilon|))^\gamma} H_\alpha(x) \frac{(\log (1/|\epsilon|))^\gamma}{x^2}dx\\
& \sim & \int_{|x|< C(\log (1/|\epsilon|))^\gamma}H_\alpha(x)|x|^{1/\gamma-1} \frac{|x|}{(\log (1/|\epsilon|))^\gamma} |x|^{-1/\gamma}dx\\ 
& + &  \int_{|x|>C (\log (1/|\epsilon|))^\gamma} H_\alpha(x) |x|^{1/\gamma-1}
\frac{(\log (1/|\epsilon|))^\gamma}{|x|}|x|^{-1/\gamma}dx.\\
\eea

Case (ii): $n>\log (1/|\epsilon|)$

Then $U \sim n^\gamma, V \sim n^{\gamma-1} \log (1/|\epsilon|).$ Hence

\bea
h_{\alpha,n}(p) & \sim & \int H_\alpha(x) \frac{n^{\gamma-1}\log (1/|\epsilon|)}
{(n^{\gamma-1}\log (1/|\epsilon|))^2+|x-n^\gamma|^2}dx.\\
\eea

We observe that this integral has already been estimated above. Namely see
Case (ii), integrals I+II. 

So we get 

\bea
\sum_{n> 10 \log (1/|\epsilon|)} h_{\alpha,n}(p) & \lesssim & \int_{|x|>(\log (1/|\epsilon|)^\gamma} H_\alpha(x)|x|^{1/\gamma-1}dx\\
& + & \int_{|x|<C(\log (1/|\epsilon|))^\gamma}H_\alpha(x)|x|^{1/\gamma-1}\left(\frac{|x|}{\log (1/|\epsilon|)}
\right)^\gamma dx.\\
\eea

We estimate next $h^\epsilon_{\beta,m}(p).$ From the above estimates for $u',v'$
we see that $|u'|<\sim v'$ and hence
$V' \sim (v')^\gamma$, $|U'|\lesssim (v')^\gamma$.
We then have:

\bea
h^\epsilon_{\beta,m}(p) & \sim & \int_{|y|< C (v')^\gamma}H_\beta(y) \frac{1}{(v')^\gamma}dy\\
& + & \int_{|y|> C(v')^\gamma } H_\beta(y) \frac{(v')^\gamma}{y^2}dy.\\
h^\epsilon_{\beta,m}(p) & \sim & \int_{|y|< C (v')^\gamma}H_\beta(y) |y|^{1/\gamma-1} \left(\frac{|y|}{(v')^\gamma}\right)^{1-1/\gamma} \frac{1}{v'}dy\\
& + & \int_{|y|> C(v')^\gamma } H_\beta(y)|y|^{1/\gamma-1} \left(\frac{(v')^\gamma}{|y|}\right)^{1+1/\gamma} \frac{1}{v'}dy.\\
\eea

Note that for $a \neq 0,$ we have that 

\bea
\log (1/|\epsilon|)/a-bu'/a -2c/|a| & <  & v' < \log (1/|\epsilon|)/a-bu'/a +2c/|b|\\
 \log (1/|\epsilon|)/a-2m\pi b/a-C & <  & v' < \log (1/|\epsilon|)/a-2m\pi b/a +C,\\
 \eea
 so
 
 \bea
 v' & > & \log (1/|\epsilon|)\\
 \eea
 and
 \bea
 m/a & < & \frac{1}{2\pi b} \log (1/|\epsilon|) (1/a-1)+C.\\
 \eea

Define
\bea
\Sigma & := & \sum_{m/a<\frac{1}{2\pi b} \log (1/|\epsilon|) (1/a-1)+C} h^\epsilon_{\beta,m}. \\
\eea

Then using the above estimates,
\bea
\Sigma &  \lesssim &
\sum_m \int_{|y|< C (  \log (1/|\epsilon|)/a-b2m\pi/a  )^\gamma}H_\beta(y) |y|^{1/\gamma-1} \\
 & * & \left(\frac{|y|}{( \log (1/|\epsilon|)/a-b2m\pi/a)^\gamma}\right)^{1-1/\gamma} 
  *  \frac{1}{ |\log (1/|\epsilon|)/a-b2m\pi/a|}dy\\
& + & \sum_m \int_{|y|> C( \log (1/|\epsilon|)/a-b2m\pi/a)^\gamma } H_\beta(y)|y|^{1/\gamma-1}\\
& * &  \left(\frac{( \log (1/|\epsilon|)/a-b2m\pi/a)^\gamma}{|y|}\right)^{1+1/\gamma} \\
& * & \frac{1}{ \log (1/|\epsilon|)/a-b2m\pi/a}dy\\
& = & I+II.\\
\eea

We study separately I and II.

\bea
I & = & I_A+I_B\\
I_A & = & \int_{|y|<C(\log (1/|\epsilon|))^\gamma}
\sum_{m/a<\frac{1}{2\pi b} \log (1/|\epsilon|) (1/a-1)+C}H_\beta(y) |y|^{1/\gamma-1} \\
 & * & \left(\frac{|y|}{( \log (1/|\epsilon|)/a-b2m\pi/a)^\gamma}\right)^{1-1/\gamma} \\
 & * & \frac{1}{ |\log (1/|\epsilon|)/a-b2m\pi/a|}dy\\
& \sim &  \int_{|y|<C(\log (1/|\epsilon|))^\gamma}H_\beta(y) |y|^{1/\gamma-1}\left(\frac{|y|}{(\log (1/|\epsilon|))^\gamma}\right)^{1-1/\gamma}dy .\\
\eea

We estimate
\bea
I_B & = & \int_{|y|=C(\log (1/|\epsilon|))^\gamma}^{\infty}\sum_{m/a<\frac{\log (1/|\epsilon|)}
{2\pi ab}-(|y|/C)^{1/\gamma}}
H_\beta(y) |y|^{1/\gamma-1} \\
 & * & \left(\frac{|y|}{ (\log (1/|\epsilon|)/a-b2m\pi/a)^\gamma}\right)^{1-1/\gamma} \\
 & * & \frac{1}{ |\log (1/|\epsilon|)/a-b2m\pi/a|}dy\\
 & \sim & 
 \int_{|y|=C(\log (1/|\epsilon|))^\gamma}^{\infty}
H_\beta(y) |y|^{1/\gamma-1} dy.\\
\eea

For II we have
\bea
II & \sim & \sum_{\frac{1}{2\pi b}\log (1/|\epsilon|)(1/a-1)>m/a>\frac{\log (1/|\epsilon|)}
{2\pi ab}-(|y|/C)^{1/\gamma} }\int_{|y|=C(\log (1/|\epsilon|))^\gamma}^{\infty}
H_\beta(y)|y|^{1/\gamma-1}\\
& * &  \left(\frac{( \log (1/|\epsilon|)/a-b2m\pi/a)^\gamma}{|y|}\right)^{1+1/\gamma} \\
& * & \frac{1}{ \log (1/|\epsilon|)/a-b2m\pi/a}dy\\
&  \lesssim & \int_{|y|=C(\log (1/|\epsilon|))^\gamma}^{\infty}
H_\beta(y) |y|^{1/\gamma-1} dy.\\
\eea

With these estimates it follows that Theorem 7 is proved for the region $A$ close
to the separatrices, in the ball $D_2$ provided that $a \neq 0.$

\bigskip

The case $a=0:$

\bigskip

We fix $(\alpha,n)$ and $(\beta,m)$ and investigate intersection points. Note that
since $a=0,$ we need $(\log (1/|\epsilon|)/b-2c/b<u'< (\log (1/|\epsilon|)/b+2c/b.$
Hence there are at most finitely many possible values for $m \sim (\log (1/|\epsilon|))/(2\pi b).$
We proceed as if there is at most one. This will suffice. 
Also note that since we assume that $|z'|<C|\epsilon|$ we also need
$v'> \log (1/|\epsilon|)-C'.$ For every integer $k>0$ we might have an intersection point
$p_{n,k}$  between $L_{\alpha,n}$ and $L^\epsilon_{\beta,m}$ with
$\log (1/|\epsilon|)-C'+k\pi<v' \leq \log (1/|\epsilon|)+(k+1)\pi.$

\bigskip

We estimate $h^\epsilon_{\beta,m}(p_{n,k}).$

\begin{lemma} When $a=0,$ then $\gamma=2.$
\end{lemma}

\begin{proof}
The inequalities $|z|<1, |w|<1$ lead to $u,v>0.$
\end{proof}

We have $u' \sim \log (1/|\epsilon|), v' \sim \log (1/|\epsilon|)+k$
so $U'+iV'= (u')^2-(v')^2+2i u'v' $ Hence if $0<k<C'' \log (1/|\epsilon)$
we have the estimate $|U'| \lesssim (\log (1/|\epsilon))^2 \sim V'$. If $k>C''\log (1/|\epsilon|)$ we have 
$U' \sim -k^2, V' \sim k \log (1/|\epsilon|).$ We consider various cases:

(i) $0<k<C'' \log (1/|\epsilon|):$

\bea
h^\epsilon_{\beta,m}(p_{n,k}) & \sim & \int_{|y|<C (\log (1/|\epsilon|))^2} H_\beta(y) \frac{1}{(\log (1/|\epsilon|))^2}dy\\
& + & \int_{|y|>C(\log (1/|\epsilon|))^2} H_\beta(y) \frac{(\log (1/|\epsilon|))^2}{y^2}dy\\
\eea
Hence
\bea
\sum_{k=0}^{C'' \log (1/|\epsilon)}h^\epsilon_{\beta,m}(p_{n,k}) & \sim & \int_{|y|<C (\log (1/|\epsilon|))^2} H_\beta(y) \frac{1}{(\log (1/|\epsilon|))}dy.\\
\eea

Adding up we get

\bea
& + & \sum_{k=0}^{C'' \log (1/|\epsilon)}\int_{|y|>C(\log (1/|\epsilon|))^2} H_\beta(y) \frac{(\log (1/|\epsilon|))^3}{y^2}dy\\
& \sim & \int_{|y|<C (\log (1/|\epsilon|))^2} H_\beta(y)|y|^{-1/2} \left(\frac{|y|}{(\log (1/|\epsilon|))^2}\right)^{1/2}dy\\
& + & \sum_{k=0}^{C'' \log (1/|\epsilon)}\int_{|y|>C(\log (1/|\epsilon|))^2} H_\beta(y) |y|^{-1/2}
\left(\frac{(\log (1/|\epsilon|))^2}{|y|}\right)^{3/2}dy.\\
\eea

Next we have case (ii).

(ii) $k>C'' \log (1/|\epsilon|):$

Using the location of $p_{n,k}$ we have the estimates 

\bea
h^\epsilon_{\beta,m}(p_{n,k}) & \sim & \int H_\beta (y) \frac{k \log (1/|\epsilon|)}
{(k \log (1/|\epsilon|))^2+(y+k^2)^2}dy\\
& \sim &  \int_{|y+k^2|<k \log (1/|\epsilon|)} H_\beta (y) \frac{1}{k \log (1/|\epsilon|)}dy\\
& + &  \int_{|y+k^2|> k\log (1/|\epsilon|)}H_\beta (y) \frac{k \log (1/|\epsilon|)}
{(y+k^2)^2}dy. \\
\eea
We define $ \Sigma$
\bea
\Sigma & := & \sum_{k>C'' \log (1/|\epsilon|)} h^{\epsilon}_{\beta,m}(p_{n,k}) \\
& = & I_A+I_B +II_A+II_B+II_C.\\
\eea
We have
\bea
I_A & \sim & \int_{y=(-(C'')^2-C'') (\log (1/|\epsilon|))^2}^{(-(C'')^2+C'') \log (1/|\epsilon|)}
\sum_{k=C'' \log (1/|\epsilon|)}^{\sqrt{-y+r(y) \log (1/|\epsilon|)}} H_\beta(y) \frac{1}{k \log (1/|\epsilon|)}dy\\
& & {\mbox{where}}\; r(y) \sim \sqrt{|y|}. \\
\eea
Hence
\bea
I_A & \sim & \int_{y=(-(C'')^2-C'') (\log (1/|\epsilon|))^2}^{(-(C'')^2+C'') \log (1/|\epsilon|)}
 H_\beta(y) \frac{1}{ \log (1/|\epsilon|)}dy\\
 I_A & \sim & \int_{y=(-(C'')^2-C'') (\log (1/|\epsilon|))^2}^{(-(C'')^2+C'') \log (1/|\epsilon|)}
 H_\beta(y) |y|^{-1/2}dy.\\
 \eea
 For  $I_B$ we have
 \bea
 I_B & \sim & \int_{y=-\infty}^{(-(C'')^2-C'') (\log (1/|\epsilon|))^2}
 \sum_{\sqrt{-y-s(y) \log (1/|\epsilon|)}}^{\sqrt{-y+r(y) \log (1/|\epsilon|)}}
 H_\beta (y) \frac{1}{k \log (1/|\epsilon|)}dy\\
 \eea
 with
 \bea
  s(y) &\sim & \sqrt{|y|}. \\
  \eea
  
  Hence
  \bea
 I_B & \sim &  \int_{y=-\infty}^{(-(C'')^2-C'') (\log (1/|\epsilon|))^2} H_\beta(y) \frac{1}{\log (1/|\epsilon|)}dy\\
 I_B & \sim &  \int_{y=-\infty}^{(-(C'')^2-C'') (\log (1/|\epsilon|))^2} H_\beta(y)|y|^{-1/2} 
\left( \frac{|y|}{(\log (1/|\epsilon|))^2}\right)^{1/2}dy.\\
\eea
Next we have for $II_A$
\bea
II_A & \sim & \int_{y=-\infty}^{((-(C'')^2-C''))(\log (1/|\epsilon|))^2}
 \sum_{k=C'' \log (1/|\epsilon|)}^{\sqrt{-y-s(y)\log (1/|\epsilon|)}}H_\beta (y) \frac{k \log (1/|\epsilon|)}
{(y+k^2)^2}dy\\ 
II_A & \sim & \int_{y=-\infty}^{((-(C'')^2-C''))(\log (1/|\epsilon|))^2}
H_\beta (y) |y|^{-1/2}dy\\ 
\eea
and
\bea
II_B & \sim & \int_{y=-\infty}^{((-(C'')^2+C''))(\log (1/|\epsilon|))^2}\sum_{k=\sqrt{-y+r(y)
\log (1/|\epsilon|)}}^\infty H_\beta (y) \frac{k \log (1/|\epsilon|)}
{(y+k^2)^2}dy\\ 
II_B & \sim & \int_{y=-\infty}^{((-(C'')^2+C''))(\log (1/|\epsilon|))^2}
H_\beta (y) |y|^{-1/2}dy.\\ 
\eea

The last term is 
\bea
II_C & \sim & \int_{((-(C'')^2+C''))(\log (1/|\epsilon|))^2}^\infty \sum_{k=C'' \log (1/|\epsilon|)}^\infty
H_\beta (y) \frac{k \log (1/|\epsilon|)}
{(y+k^2)^2}dy\\ 
II_C & \sim & \int_{((-(C'')^2+C''))(\log (1/|\epsilon|))^2}^\infty 
H_\beta (y) \frac{ \log (1/|\epsilon|)}
{y+(C'' \log (1/|\epsilon|))^2}dy\\ 
II_C & \sim & \int_{((-(C'')^2+C''))(\log (1/|\epsilon|))^2}^{(\log (1/|\epsilon))^2} 
H_\beta (y) \frac{1}{ \log (1/|\epsilon|)}dy\\ 
& + & \int_{(\log (1/|\epsilon|))^2}^\infty 
H_\beta (y) \frac{ \log (1/|\epsilon|)}
{y}dy\\ 
II_C & \sim & \int_{((-(C'')^2+C''))(\log (1/|\epsilon|))^2}^{(\log (1/|\epsilon))^2} 
H_\beta (y) |y|^{-1/2}\frac{|y|^{1/2}}{ \log (1/|\epsilon|)}dy\\ 
& + & \int_{(\log (1/|\epsilon|))^2}^\infty 
H_\beta (y) |y|^{-1/2}\frac{|y|^{-1/2}}{ \log (1/|\epsilon|)}dy.\\
\eea

Next we estimate $h_{\alpha,n}(p_{n,k}).$ Note however, that the estimates for the case $a\neq 0$ still applies to $h_{\alpha,n}.$ This condition was not used to estimate $h_\alpha.$ 
Hence we are done with the proof of Theorem 7 for the case of intersection points in $A$

We next consider the set $B$ of points in $\Delta^2(0,C|\epsilon|)$ defined above as consisting of points which are at distance at least
$r|\epsilon|$ from all separatrices.

\section{Proof of Theorem 7 for points in $B$, i.e. points in $D_2$ which are at distance
at least $r|\epsilon|$ from the separatrices.}

\bigskip

We estimate $H_\alpha$ on $L_{\alpha,n}\cap B.$ We can assume $a \neq 1$, otherwise
flip the axes. So,

\bea
r|\epsilon| & <  & |z|< C|\epsilon|\\
r|z| & < & e^{-v} < C |\epsilon|\\
\log (1/|\epsilon|) -C' & < & v < \log (1/|\epsilon|)+C'. \\
\eea
Similarly
\bea
r|\epsilon| & < & |w| < C |\epsilon|\\
r|\epsilon| & < & e^{-bu-av} < C |\epsilon|\\
\log|\epsilon| -C'' & < & -bu-av  < \log |\epsilon| + C''\\
\log (1/|\epsilon|) -C''-av & < & bu < -av+\log (1/|\epsilon|)+C''\\
(1-a) \log (1/|\epsilon|)-C & < & bu < (1-a) \log (1/|\epsilon|)+C\\
\frac{1-a}{b} \log (1/|\epsilon|) -C & < & u < \frac{1-a}{b} \log (1/|\epsilon|)+C.\\
\eea

Using these estimates on $(u,v)$ and similarly for $(u',v')$, Lemma 10 shows that there is an
integer $N$ independent of $\epsilon$
so that if we take  any two plaques of two leaves $L_\alpha, L^\epsilon_{\beta}$, then they
intersect in $B$ in at most $N$ points.  In $U,V$ coordinates,

\bea
(u+iv)^\gamma & = & U +iV,\\
\eea
\noindent hence
\bea
V & \sim & (\log (1/|\epsilon|)^\gamma\\
|U| & < \sim & (\log (1/|\epsilon|))^\gamma.\\
\eea
This gives

\bea
h_{\alpha,n} & \sim & \int H_\alpha (x) \frac{(\log (1/|\epsilon|))^\gamma}
{(\log (1/|\epsilon|))^{2\gamma}+(x-U)^2}dx\\
& \sim & \int_{|x|< C (\log (1/|\epsilon|))^\gamma} H_\alpha (x) \frac{1}{(\log (1/|\epsilon|))^\gamma}dx\\
& + & \int_{|x|> (\log (1/|\epsilon|)^\gamma} H_\alpha (x) \frac{(\log (1/|\epsilon|))^\gamma}
{x^2}dx\\
& \sim & \int_{|x|< C (\log (1/|\epsilon|))^\gamma} H_\alpha (x)(|x|+1)^{1/\gamma-1} \frac{|x|+1}{(\log (1/|\epsilon|))^\gamma} (|x|+1)^{-1/\gamma}dx\\
& + & \int_{|x|> (\log (1/|\epsilon|)^\gamma} H_\alpha (x) (|x|+1)^{1/\gamma-1}
\frac{(\log (1/|\epsilon|))^\gamma}{|x|+1}(|x|+1)^{-1/\gamma}dx.\\
\eea

It follows from these estimate applied to $H_\beta$ as well, that Theorem 7 is
valid for intersection points in $B.$

\section{Theorem 7 for $D_3=\Delta^2(0,\delta)\setminus \Delta^2(0,C|\epsilon|)$}

\bigskip

There are 3 regions to consider: \\

\bea
D_3 & = & R_1 \cup R_2 \cup R_3\\
R_1 & = & \{C|\epsilon|<|z|<\delta, C|\epsilon|<|w|<\delta\}\\
R_2 & = & \{C|\epsilon|<|z|<\delta, |w|<C |\epsilon|\}\\
R_3 & = & \{|z|<C |\epsilon|,C|\epsilon|<|w|<\delta. \}\\
\eea

Note that since we have assumed $a \neq 1,$ the cases of $R_2$ and $R_3$ are
not completely symmetric.  We will leave it to the reader to verify
that the estimates we do later for $R_2$ nevertheless hold for $R_3.$

\section{Theorem 7 for $R_1,$ the diagonal part of $D_3$}

We first outline our approach. Fix parameters $\alpha,\beta$ and corresponding plaques
$L_{\alpha,n}, L^\epsilon_{\beta,m}.$ Next we divide $R_1$ into
dyadic components, rings, $\{R(p)\}$ in the $z-$ direction, $e^{-p-1}<|z|<e^{-p}, C|\epsilon|<|w|<\delta$ Then
we estimate $h_\alpha$ and $h_\beta$ on $L_{\alpha,n}\cap R(p)$ and
$L^\epsilon_{\beta,m}\cap R(p)$ respectively. Next, for fixed
$\alpha,\beta, n,m$ we estimate the values of $p$ where the leaves $L_{\alpha,n},L^\epsilon_{\beta,m}$
might intersect, and the number of intersection points for each such $p.$
Putting this information together we can estimate the contribution from $R_1$ to the
geometric wedge product.

Pick a plaque $L_{\alpha,n}$ and a point $(z,w)$ in $L_{\alpha,n}\cap R(p)$ parametrized by $(u,v).$ Then

\bea
e^{-p-1} & < & |z|=e^{-v} <e^{-p}\\
\log (1/\delta) & < & v  < \log (1/|\epsilon|)-C\\
\log (1/\delta) & < & p  < \log (1/|\epsilon|)-C,\\
\eea
\noindent and
\bea
2n\pi & <  & u < 2(n+1)\pi. \\
\eea

For $w$ we have
\bea
C|\epsilon| & < & |w|<\delta\\
\log (1/\delta) & < & bu+av < \log (1/|\epsilon|)-\log C\\
\frac{\log (1/\delta)}{b}-av/b & < & u < \frac{\log (1/|\epsilon|)-\log C}{b}-av/b.\\
\eea

We divide into  cases depending on whether $a \neq 0$ or $a=0.$\\

First, assume $a \neq 0.$
We choose a constant $0<s<1$ so that $\frac{1}{2}<1+\frac{2sb\pi}{a}<\frac{3}{2}.$\\

Case (i): $a \neq 0, n<sp,$ then\\

\bea
(u+iv)^\gamma & = & U+iV\\
& \sim & U+ip^\gamma, |U| <\sim p^\gamma.\\
\eea
So we have
\bea
H_{\alpha,n} & \sim & \int_{|x|<C p^\gamma} \tilde{H}_\alpha(x)/p^\gamma dx
+\int_{|x|>Cp^\gamma} \tilde{H}_\alpha(x) p^\gamma/x^2dx\\
& \sim & \int_{|x|<C p^\gamma} \tilde{H}_\alpha(x)|x|^{1/\gamma-1} \left(\frac{|x|}{p^\gamma}\right)^{1-1/\gamma}
\frac{1}{p} dx\\
& + & \int_{|x|>Cp^\gamma} \tilde{H}_\alpha(x) |x|^{1/\gamma-1}\left(\frac{p^\gamma}{|x|}\right)^{1+1/\gamma}\frac{1}{p}dx.\\
\eea

Case (ii) $a \neq 0, n>sp$\\

\bea
(u+iv)^\gamma & = & U+iV\\
& \sim & n^\gamma+ipn^{\gamma-1}.\\
\eea
Then
\bea
H_{\alpha,n} & \sim & \int_{|x-n^\gamma|\leq pn^{\gamma-1}}\tilde{H}_\alpha(x)\frac{1}{pn^{\gamma-1}}dx\\
& + & \int_{n^\gamma/2>|x-n^\gamma|>pn^{\gamma-1}}\tilde{H}_\alpha(x)\frac{pn^{\gamma-1}}{|x-n^\gamma|^2}dx\\
& + &\int_{n^\gamma/2<|x-n^\gamma|<2n^{\gamma}}\tilde{H}_\alpha(x)\frac{pn^{\gamma-1}}{n^{2\gamma}}dx\\
& + & \int_{|x-n^\gamma|>2n^\gamma}\tilde{H}_\alpha(x)\frac{pn^{\gamma-1}}{x^2}dx\\
& = & I+II+III+IV\\
\eea

We will usually leave the case $a=0$ to the reader. 

Case (iii) $ a=0$\\

\bea
\gamma & = & 2\\
(u+iv)^2 & = & u^2-v^2+2iuv. \\
\eea
So
\bea
h_{\alpha,n} & = & \int H_{\alpha,n}(x) \frac{uv}{(uv)^2+(x-(u^2-v^2))^2}dx\\
& \sim &  \int H_{\alpha,n}(x) \frac{uv}{(uv)^2+(x+v^2-u^2)^2}dx\\ 
& \sim &  \int _{|x+v^2-u^2|<uv} H_{\alpha,n}(x)\frac{uv}{(uv)^2+(x+v^2-u^2)^2}dx\\ 
& + &  \int_{|x+v^2-u^2|>uv} H_{\alpha,n}(x)\frac{uv}{(uv)^2+(x+v^2-u^2)^2}dx\\ 
& \sim &  \int _{|x+v^2-u^2|<uv} H_{\alpha,n}(x)\frac{1}{uv}dx\\ 
& + &  \int_{|x+v^2-u^2|>uv} H_{\alpha,n}(x)\frac{uv}{(x+v^2-u^2)^2}dx\\ 
\eea
To estimate the integrals we divide into cases.

\bea
u /10& < & v < 10u:\\
& \sim &  \int _{|x|<\sim u^2} H_{\alpha,n}(x)\frac{1}{u^2}dx
 +   \int_{|x|>\sim u^2} H_{\alpha,n}(x)\frac{u^2}{x^2}dx\\ 
& \sim &  \int _{|x|<\sim u^2} H_{\alpha,n}(x)|x|^{-1/2}\left(\frac{|x|}{u^2}\right)^{1/2} \frac{1}{u}dx\\ 
& + &  \int_{|x|>\sim u^2} H_{\alpha,n}(x)|x|^{-1/2}\frac{u^2}{|x|}\frac{1}{|x|^{1/2}}dx.\\ 
\eea

Suppose now

\bea
u & < & v/10\\
h_{\alpha,n}& \sim &  \int _{|x+v^2|<\sim uv} H_{\alpha,n}(x)\frac{1}{uv}dx
 +   \int_{|x+v^2|>\sim uv} H_{\alpha,n}(x)\frac{uv}{(x+v^2-u^2)^2}dx\\ 
& \sim &  \int _{|x+v^2|<\sim uv} H_{\alpha,n}(x)|x|^{-1/2} \frac{1}{u}dx\\ 
& + &  \int_{|x+v^2|>\sim uv} H_{\alpha,n}(x)\frac{uv}{(x+v^2)^2}dx\\  
& \sim &  \int _{|x+v^2|<\sim uv} H_{\alpha,n}(x)|x|^{-1/2} \frac{1}{u}dx\\ 
& + &  \int_{|x+v^2|>\sim uv} H_{\alpha,n}(x)|x|^{-1/2}\frac{uv|x|^{1/2}}{(x+v^2)^2}dx.\\  
\eea

The last case is

\bea
u & > & 10v.\\
h_{\alpha,n}& \sim &  \int _{|x-u^2|<\sim uv} H_{\alpha,n}(x)|x|^{-1/2}\frac{1}{v}dx\\ 
& + &  \int_{|x-u^2|>\sim uv} H_{\alpha,n}(x)|x|^{-1/2}\frac{uv|x|^{1/2}}{(x-u^2)^2}dx.\\ 
\eea

So now we have estimated $h_\alpha$ on $L_{\alpha,n}\cap R(p)$. The analogue estimates
are valid for $H_\beta$ on $L^\epsilon_{\beta,m}.$ The reason is that $e^{-p}>>|\epsilon|.$\\

Our next step is to locate for which $R(p)$ there is an intersection between $L_{\alpha,n}$
and $L^\epsilon_{\beta,m}.$ 

Fix $L_{\alpha,n}$ and $L^\epsilon_{\beta,m}$ and assume $(z,w)\in L_{\alpha,n}
\cap L^\epsilon_{\beta,m}.$ We can write 

\bea
z & = & e^{i(\zeta+(\log |\alpha|)/b}\\
\zeta & = & u+iv\\
2n\pi & < & u < 2(n+1)\pi\\
|z| & = & e^{-v}.\\
\eea

Also $(z,w) = \Phi_\epsilon(z',w'), (z',w')\in L_{\beta,m}.$

\bea
z' & = & e^{i(\zeta'+(\log |\alpha|)/b}\\
\zeta' & = & u'+iv'\\
2m\pi & < & u' < 2(m+1)\pi\\
|z'| & = & e^{-v'}. \\
\eea
So
\bea
z & = & \alpha(\epsilon)+e^{i(\zeta'+(\log |\beta|)/b)}+\epsilon {\mathcal {O}}(z',w')\\
w & = & \beta(\epsilon)+\beta e^{i \lambda(\zeta'+(\log |\beta|)/b)}+\epsilon {\mathcal {O}}(z',w').\\
\eea

Our goal is to locate which $R(p)$ the point $(z,w)$ can belong to. So we need to find $p$ so that
$e^{-p-1}<|z|=e^{-v}<e^{-p},$ i.e. we need to get a good estimate for $v$ in terms of $\alpha,\beta,n,m.$

There are 4 unknowns, $u,v,u',v'.$ However, $u \sim 2n\pi, u' \sim 2m\pi,$ so we only have $v,v'$
left. Also we have two equations for the $z$ and $w$ coordinates respectively.
 (In fact, since these are complex equations, we have $4$ real equations for the two
 real unknowns $v,v'.$)
 
 Before we proceed we show at first that for there to be an intersection, we actually
 must require that $n$ and $m$ are very close.

\begin{lemma}
If $L_{\alpha,n}$ and $L^\epsilon_{\beta,m}$ intersect in $R_1$, it follows that
$|m -n| \leq 1.$
\end{lemma}

\begin{proof}

Recall that
$$
\Phi_\epsilon(z,w)=(\alpha(\epsilon),\beta(\epsilon))+(z,w)+\epsilon {\mathcal O}(z,w).
$$
If $\delta$ is chosen small enough, this implies that
$|\epsilon {\mathcal O}(z,w)|\leq  \sigma |\epsilon|$ for any given $0<\sigma<<1.$

We pick two plaques, $L_{\alpha,n}, L^\epsilon_{\beta,m}$
and consider intersection points in $R_1.$ Let $S>0$ be such that
$ |\epsilon|/S< |\alpha(\epsilon)|-\sigma |\epsilon|,|\beta(\epsilon)|-\sigma |\epsilon|
<|\alpha(\epsilon)|+\sigma |\epsilon|,|\beta(\epsilon)|+\sigma |\epsilon|<S.$
Note that if we increase the constant $C$ used in the definition of $D_4$, we can still use the same $S.$
When the point is in $L_{\alpha,n}$ we have

\bea
z & = & e^{i(\zeta+(\log |\alpha|)/b}\\
|z| & = & e^{-v}\\
\log (1/|\delta|)&  <v & < \log (1/|\epsilon|)-C\\
w & = & \alpha e^{i \lambda (\zeta+(\log |\alpha|)/b)}\\
|w| & = & e^{-bu-av}. \\
\eea

If it is also in $L_{\beta,m},$ then
\bea
 \\
z' & = & e^{i(\zeta'+(\log |\beta|)/b)}\\
|z'| & = & e^{-v'}\\
\log (1/|\delta|) < & v' & < \log (1/|\epsilon|)-C\\
w' & = & \beta e^{i \lambda (\zeta'+(\log |\beta|)/b)}\\
|w'| & = & e^{-bu'-av'}. \\
\eea
Since
\bea
L^\epsilon_{\beta,n} & = & \Phi_\epsilon (L_{\beta,m}),\\
\eea
\noindent the image point can be written
\bea
Z & = & \alpha(\epsilon)+e^{i(\zeta'+(\log |\beta|)/b)}+ \epsilon {\mathcal O}(z',w')\\
W & = & \beta(\epsilon)+ \beta e^{i \lambda (\zeta'+(\log |\beta|)/b)}+ \epsilon {\mathcal O}(z',w')\\
\eea

Consider an intersection point in $R_1$ and set $\zeta'=\zeta+c+id.$ Then
\bea
z & = & Z\\
e^{-v-d}-S|\epsilon| & <  & e^{-v}< e^{-v-d}+S|\epsilon|\\
e^{-d} -Se^v|\epsilon| & <  & 1 < e^{-d}+Se^v|\epsilon|.\\
\eea
So
\bea
S e^v |\epsilon| & < & S (1/(C |\epsilon|)|\epsilon|=S/C <<1.\\
|d| & < & 2Se^v |\epsilon|< 2S/C.\\
\eea

For the other coordinate,
\bea
w & = & W\\
e^{-bu-bc-av-ad}-S|\epsilon| & < & e^{-bu-av} < e^{-bu-bc-av-ad}+S|\epsilon| \\
e^{-bc-ad} -Se^{bu+av}|\epsilon| & < & 1 < e^{-bc-ad}+Se^{bu+av}|\epsilon|.\\
Se^{bu+av} |\epsilon| & < & S/C<<1.\\
|bc+ad| & < & 2Se^{bu+av}|\epsilon|<2S/C\\
|bc| & < & |bc+ad|+ |a||d|\\
& < & 2Se^{bu+av}|\epsilon|+ |a|  2Se^v |\epsilon|.\\
\eea
So
\bea
|c| & < & \frac{1}{|b|} \left(2Se^{bu+av}|\epsilon|+ |a|  2Se^v |\epsilon|\right)\\
& < & 2S \frac{1+|a|}{C|b|}\; {\mbox{and}}\\
|c+id| & < & \frac{2S}{C} \left(1+\frac{1+|a|}{|b|}\right)<<1.\\
\eea

\end{proof}

It is also convenient to show that $\alpha$ and $\beta$ must be very close
if there is an intersection. We estimate first the modulus and next the angle and finally we
combine them.
 
\begin{lemma}
Suppose $L_{\alpha,n}$ intersects $L^\epsilon_{\beta,m}$ in $R_1.$
Then $$|\log (|\beta|/|\alpha|)| \leq 2S|\epsilon| \left[
e^v \left(b+|a|\right)+e^{bu+av}\right].$$
\end{lemma}

\begin{proof}

We have
\bea
z & = & Z\\
e^{i(\zeta+(\log |\alpha|)/b)} &  = & \alpha(\epsilon)+e^{i(\zeta+c+id+(\log |\beta|)/b)}+\epsilon
\mathcal {O}(z',w').\\
\eea
So
\bea
e^{i(\zeta+(\log |\alpha|)/b)}\left[1-e^{ic-d+i(\log (|\beta|/|\alpha|)/b)}\right] & = & \alpha(\epsilon)+\epsilon\mathcal {O}(z',w').\\
\eea
Taking the modulus
\bea
e^{-v}\left|1-e^{ic-d+i(\log (|\beta|/|\alpha|)/b)}\right| & \leq & S|\epsilon|\\
\left|1-e^{ic-d+i(\log |\beta|/|\alpha|)/b)}\right| & \leq & Se^v|\epsilon|<<1.\\
\eea
This gives
\bea
\left|i(c+(\log |\beta|/|\alpha)/b)-d\right| & \leq & 2Se^v|\epsilon|\\
\left|\log (|\beta|/|\alpha|)/b\right| & \leq & 2Se^v|\epsilon|+2Se^{bu+av}|\epsilon|/b+2S(|a|/b) e^v |\epsilon|.\\
\eea

The Lemma follows.

\end{proof}

We remark that the lemma as stated is slightly inaccurate. We only can conclude the estimate modulo
$2\pi.$ However, the parameters $e^{-2\pi b} \leq |\alpha|,|\beta|<1$ so this problem arises when say $|\alpha|$ is close to $1$ and $|\beta|$ is close to $e^{-2\pi b}.$ We ignore this technicality which just means
that $|\alpha|$ and $|\beta|$ get close after we follow the leaf $L_\alpha$ once around $0$ counterclockwise.

\begin{lemma}
Write $\beta/\alpha=|\beta/\alpha|e^{i\theta}.$ If there are intersection points in $R_1$,
$\theta$ is close to $0$ mod $2\pi.$ More precisely:
$$
|\theta| \leq 2Se^{bu+av}|\epsilon|\left[|a|/b +|a|/b+1 \right]  + 2 S|\epsilon|e^v\left[|a|^2/b+b+(|a|+|a|^2/b)\right].
$$
\end{lemma}

\begin{proof}
We again use the parametrization.
\bea
w & = & W\\
\alpha e^{i\lambda(\zeta+(\log |\alpha|/b))} & = & 
\beta(\epsilon)+\beta e^{i\lambda(\zeta+c+id+(\log |\beta|/b))}+\epsilon \mathcal {O}(z',w').\\
\eea
So
\bea
 \beta(\epsilon)+\epsilon\mathcal {O}(z',w') & = & \alpha e^{i\lambda(\zeta+(\log |\alpha|/b))}\left[1-\frac{\beta}{\alpha} e^{i\lambda(c+id+(\log |\beta|/b))}\right] \\
 S e^{bu+av}|\epsilon| & \geq & 
\left|1-\frac{\beta}{\alpha} e^{i\lambda(c+id+(\log |\beta|/b))}
\right|.\\
\eea
Hence
\bea
S e^{bu+av}|\epsilon| & \geq &
\left|1-\frac{\beta}{\alpha} e^{[-bc-ad-(\log|\alpha|/b)]+i[ac-bd+a(\log (|\beta|/|\alpha|))/b]}\right| \\
1>>  S e^{bu+av}|\epsilon| & \geq &
  \left|1-e^{[-bc-ad]+i[\theta+ac-bd+a(\log (|\beta|/|\alpha|))/b]}\right| \\
 2 S e^{bu+av}|\epsilon| & \geq & 
 \left|\theta+ac-bd+a(\log (|\beta|/|\alpha|))/b\right|.\\
 \eea
 Therefore
 \bea
|\theta| & 
  \leq &   |ac|+|bd|+|a||\log (|\beta|/|\alpha|)|/b+
 2 S e^{bu+av}|\epsilon|  \\
&  \leq &  Se^{bu+av}|\epsilon|\left[2|a|/b +2|a|/b+2\right] \\ & +&   S|\epsilon|e^v\left[2|a|^2/b+2b+2(|a|+|a|^2/b)
 \right].\\
\eea
Which gives the estimate.

\end{proof}

Next we  locate more precisely the intersections of $L_{\alpha,n}$ and
$L^\epsilon_{\beta,m}$ in $R_1.$
Let
\bea
z & = & Z.\\
\eea
Then
\bea
e^{i\zeta+i(\log |\alpha|)/b} & = & \alpha(\epsilon)+e^{i\zeta'+i(\log |\beta|)/b}+\epsilon \mathcal{O}(z',w').\\
\eea
We define $\Delta$ by
\bea
\zeta' & = & \zeta+ \Delta.\\
\eea
Then
\bea
e^{i\zeta+i(\log |\alpha|)/b}-e^{i\zeta+i\Delta+i(\log |\beta|)/b} & = & 
\alpha(\epsilon)+\epsilon \mathcal{O}\\
e^{i\zeta+i(\log |\alpha|)/b}\left[1-e^{i\Delta+i(\log (|\beta|/|\alpha|)/b}\right] & = & 
\alpha(\epsilon)+\epsilon \mathcal{O}\\
1-e^{i\Delta+i(\log (|\beta|/|\alpha|)/b}& = & 
e^{-i\zeta-i(\log |\alpha|)/b}\left[\alpha(\epsilon)+\epsilon \mathcal{O}\right].\\
\eea
This gives
\bea
2k\pi+\Delta+ (\log (|\beta|/|\alpha|)/b& = & ie^{-i\zeta-i(\log |\alpha|)/b}\left[\alpha(\epsilon)+\epsilon \mathcal{O}\right]\\
& + & \mathcal{O}( \epsilon e^{-i\zeta})^2.\\
\eea
Using
\bea
w & = & W, \\
\eea
we have
\bea
\alpha e^{i\lambda (\zeta+(\log |\alpha|)/b)} & = & \beta(\epsilon)
+\beta e^{i\lambda(\zeta+\Delta+(\log |\beta|)/b)}+\epsilon \mathcal {O}\\
e^{i\lambda \zeta}\left[\alpha e^{i\lambda (\log |\alpha|)/b)}
-\beta e^{i\lambda(\Delta+(\log |\beta|)/b)}\right] & = & \beta(\epsilon)+\epsilon\mathcal{O}\\
& = & e^{i\lambda \zeta} e^{i\lambda (\log |\alpha|)/b)}\\
& * & \left[\alpha
-\beta e^{i\lambda( ie^{-i\zeta-i(\log |\alpha|)/b}\left[\alpha(\epsilon)+\epsilon \mathcal{O}\right]   )}\right].\\\eea
So
\bea
1-\frac{\beta}{\alpha}e^{i\lambda(\Delta+(\log|\beta|/|\alpha|)/b)} & = &
e^{-i\lambda \zeta}\frac{e^{-i\lambda (\log|\alpha|)/b)}}{\alpha}\left[\beta(\epsilon)+\epsilon \mathcal{O}\right]\\
-\log\left(\frac{\beta}{\alpha}\right)+i\lambda(\Delta+(\log|\beta|/|\alpha|)/b) & \sim &
e^{-i\lambda \zeta}\frac{e^{-i\lambda (\log|\alpha|)/b)}}{\alpha}\left[\beta(\epsilon)+\epsilon \mathcal{O}\right].\\
\eea
We get:
\bea
\Delta +\log (|\beta|/|\alpha|)/b -  \frac{1}{i\lambda}\log\left(\frac{\beta}{\alpha}\right)
& = & \frac{1}{i \lambda}e^{-i\lambda \zeta}\frac{e^{-i\lambda (\log|\alpha|)/b)}}{\alpha}\left[\beta(\epsilon)+\epsilon \mathcal{O}\right] \\
& + & \mathcal{O}(e^{-i\lambda \zeta}\epsilon)^2.\\
\eea

Adding the two expressions with $\Delta:$ 

\begin{lemma}
Suppose that $L_{\alpha,n} \cap L^\epsilon_{\beta,m} \cap R_1 \neq \emptyset.$ Then:
\bea
-\frac{1}{i\lambda}\log\left(\frac{\beta}{\alpha}\right) & = & 
ie^{-i\zeta-i(\log |\alpha|)/b}\left[\alpha(\epsilon)+\epsilon \mathcal{O}\right]
+ \frac{1}{i \lambda}e^{-i\lambda \zeta}\frac{e^{-i\lambda (\log|\alpha|)/b)}}{\alpha}\left[\beta(\epsilon)+\epsilon \mathcal{O}\right]\\
& + &  \mathcal{O}( \epsilon e^{-i\zeta})^2+ \mathcal{O}(e^{-i\lambda \zeta}\epsilon)^2.\\
\eea
\end{lemma}

To continue the search for intersection points of $L_{\alpha,n}, L^\epsilon_{\beta,m}$ in $R_1$,
we divide $R_1$ into $3$ pieces. We let $C_1>1$ be a large constant.

 \bea
 R_{1A} & = & \{C|\epsilon|<|z|,|w|<\delta, C_1|w|\leq|z|\}\\
 R_{1B} & = & \{C|\epsilon|<|z|,|w|<\delta, C_1|z|\leq |w|\}\\
 R_{1C} & = & \{C|\epsilon|<|z|,|w|<\delta, |z|\leq C_1 |w|\leq C_1^2|z|\}\\
 \eea

Here the constant $C_1$ is chosen to work in the slope estimates before Lemma 21.

Observe that $R_{1A}$ and $R_{1B}$ are similar. We will leave it up to the reader to
verify  the estimates for $R_{1B}.$

\bigskip

\section{Theorem 7 for $R_{1A}$, the part of $R_1$ close to the $z-$axis}

We will assume that $a \neq 0$ and leave the verification of the case $a=0$ to the reader.
If $|w|<< |z|,$ then the second term in the expression for $\log (\beta/\alpha)$ in Lemma 18 on the right dominates and we get
\bea
e^{av+bu} |\epsilon| &  \sim &  |(\beta/\alpha) -1|\\
2n\pi & < & u < 2(n+1)\pi\\
av & \sim & \log |(\beta/\alpha)-1|+\log (1/|\epsilon|)-2nb\pi.\\
\eea
So
\bea
|v-\frac{\log |(\beta/\alpha)-1|+\log (1/|\epsilon|)-2nb\pi}{a}| & < & C\\
C|\epsilon| & < & e^{-v}<\delta\\
\log 1/\delta & < & v < \log 1/|\epsilon|-C\\
p & < & v < p+1,\\
\eea
\noindent see the beginning of Section 7.

\begin{lemma}
For intersection points in $R_{1A}$, There is a constant C' such that
$$
\frac{C'|\epsilon|}{\delta} <|\beta-\alpha|<\frac{1}{C'}.
$$
\end{lemma}

\begin{proof}
Since $e^{av+bu}|\epsilon|\sim |\frac{\beta}{\alpha}-1| \sim |\beta-\alpha|$ and
$e^{av+bu} = 1/|w|$ we have $|\beta-\alpha| \sim |\epsilon|/|w|.$ But $C|\epsilon|<|w|
<|z|/C<\delta/C.$ The lemma follows.

\end{proof}

\begin{lemma}
Suppose that $L_{\alpha,n} $ intersects $L^\epsilon_{\beta,m}$ in $R_{1A}$.
Then the intersection points must be  in $R(p)$ for some 
$$|p-\frac{\log |(\beta/\alpha)-1|+\log (1/|\epsilon|)-2nb\pi}{a}|  <  C.$$
For the plaque to enter $R_{1}$ we further need $n$ to satisfy
$$\log |(\beta/\alpha)-1|+\log (1/|\epsilon|)-2nb\pi \in I$$
where $I$ is the interval with endpoints $a \log 1/\delta, a\log (1/|\epsilon|)-aC$
\end{lemma}

Our next step is to verify that there is a uniform bound on the number of intersection
points of $L_{\alpha,n},L^\epsilon_{\beta,m}$ in $R_{1A}.$

In order to study the number of intersections between plaques, we compare
their slopes:

\bigskip

Suppose $(z,w)=(Z,W):=\Phi_\epsilon (z',w')$ is an intersection point of $L_{\alpha,n}$
and $L^\epsilon_{\beta,m}$ in $R_1.$
The slope $S_1$ of $L_{\alpha,n}$ is $\lambda w/z.$ The slope of the perturbed leaf
is $S_2.$ We choose the constant $C_1$ used in the definition of $R_{1A},R_{1B}, R_{1C}$
in the following estimates.

\bea
\Phi_\epsilon'(z',w')(z', \lambda w') & = & (z'+\epsilon \mathcal{O}(z',w'),\lambda w'+\epsilon \mathcal{O}(z',w'))\\
S_2 & = & \frac{\lambda w'+\epsilon \mathcal{O}(z',w')}{z'+\epsilon \mathcal{O}(z',w')}\\
& = & \frac{\lambda W-\lambda\beta(\epsilon)+\epsilon \mathcal{O}(Z,W)}{Z-\alpha(\epsilon)+
\epsilon \mathcal{O}(Z,W)}\\
S_2 & = & \frac{\lambda w-\lambda\beta(\epsilon)+\epsilon \mathcal{O}(z,w)}{z-\alpha(\epsilon)+
\epsilon \mathcal{O}(z,w)}\\
S_2-S_1 & = & \frac{\lambda w-\lambda\beta(\epsilon)+\epsilon \mathcal{O}(z,w)}{z-\alpha(\epsilon)+
\epsilon \mathcal{O}(z,w)}-\lambda w/z\\
& = & \frac{-\lambda \beta(\epsilon)z+\lambda w \alpha(\epsilon)+\epsilon \mathcal{O}(z^2,zw,w^2)}
{z(z-\alpha(\epsilon)+\epsilon \mathcal{O}(z,w))}\\
\frac{1}{C_1}|z|\leq  |w|\leq C_1|z| & : & S_2-S_1\sim  \frac{\lambda}{z^2}\left(w \alpha(\epsilon)-z\beta(\epsilon)\right)\\
|w|>C_1|z| &:& S_2-S_1\sim \frac{\epsilon w}{z^2}\\
|w|<\frac{1}{C_1}|z| & : & S_2-S_1 \sim \frac{\epsilon}{z}\\
\eea

\begin{lemma}
There  is at most a uniformly bounded number of intersection points in 
 $R_{1A}.$ 
\end{lemma}

\begin{proof}
The case of $R_{1A},R_{1B}$ follows from slope estimates. For the case $R_{1C}$, note that
leaves might be tangent when $(w/z)$ is close to $\beta(\epsilon)/\alpha(\epsilon)$. They both have
slope about $\lambda.$ But since we assume that $\lambda \neq \beta'(0)/\alpha'(0)$, this
tangency is at most of order $2.$
\end{proof}

We estimate the contribution to $T \wedge_g T^\epsilon$ from $R_{1A}.$ We assume again that $a \neq 0.$
By Lemma 18, the parameters $\alpha,\beta$ are restricted to the values:
$e^{-2\pi b}<|\alpha|,|\beta|<1, 1/C>|\beta-\alpha|>C|\epsilon|/\delta.$ 
So fix $\alpha,\beta.$ Next, by Lemma 15, we can set $n=m$ to be some integer in the interval given by Lemma 19. The case $n=m\pm 1$ is similar. 
Because of the finiteness of the number of intersection points, see Lemma 20, 
we can set 
$$p=p(n)=\frac{\log |(\beta/\alpha)-1|+\log (1/|\epsilon|)-2nb\pi}{a}$$
\noindent and consider only one intersection point. 
Then we multiply the values of $H_{\alpha,n}$ and $H_{\beta,n}$ using the formulas in Case (i) or (ii)
depending on whether $n<sp$ or $n>sp.$  We then add these products over $n$ and
integrate the result over $d\mu(\alpha)d\mu(\beta).$

Case (i), $n<sp:$

\bea
n & < & s \frac{\log |(\beta/\alpha-1|+\log 1/|\epsilon|-2nb\pi}{a}\\
n(1+\frac{2sb\pi}{a}) & < & s \frac{\log |(\beta/\alpha-1|+\log 1/|\epsilon|}{a}.\;{\mbox{Recall that}}\\
1/2 < 1+\frac{2sb\pi}{a} & < & 3/2.\; {\mbox{We get}}\\
n & < & \frac{s}{1+\frac{2sb\pi}{a}}   \frac{\log |(\beta/\alpha-1|+\log 1/|\epsilon|}{a}=:n(\alpha,\beta,\epsilon).\\
\eea

In this case we have the following estimates at intersection points.

\bea
h_{\alpha,n} 
& \sim & \int_{|x|<C v^\gamma} H_\alpha(x)|x|^{1/\gamma-1} \left(\frac{|x|}{v^\gamma}\right)^{1-1/\gamma}
\frac{1}{v} dx\\
& + & \int_{|x|>Cv^\gamma} H_\alpha(x) |x|^{1/\gamma-1}\left(\frac{v^\gamma}{|x|}\right)^{1+1/\gamma}\frac{1}{v}dx\\
h_{\beta,m} 
& \sim & \int_{|y|<C (v')^\gamma} H_\beta(y)|y|^{1/\gamma-1} \left(\frac{|y|}{(v')^\gamma}\right)^{1-1/\gamma}
\frac{1}{v'} dy\\
& + & \int_{|y|>C(v')^\gamma} H_\beta(y) |y|^{1/\gamma-1}\left(\frac{(v')^\gamma}{|y|}\right)^{1+1/\gamma}\frac{1}{v'}dy\\
\eea

Here we have used that $v$ and $\rho$ are comparable. In fact from the estimate in the beginning
of the section, we see that $u\sim n, n<sp$ so $u<\sim v$, hence $\rho \sim v.$

Here $v,v' \sim \frac{\log |(\beta/\alpha)-1|+\log 1/|\epsilon|-2nb\pi}{a}.$ This allows us to sum over $v$ instead of over $n,$ $\log 1/\delta <v<\log 1/|\epsilon|-C.$

We need to estimate $\sum_{v} h_{\alpha,n}h^\epsilon_{\beta,n}$ and then integrate the answer over
the measure $\mu(\alpha)\mu(\beta)$.\\

Note we will majorize the sum by the product $\sum_{v}h_{\alpha,n} \sum_m h^\epsilon_{\beta,m}$
Then we use the dominated convergence theorem. 

We finally have

\bea
& \sim & \sum_{v=\log 1/\delta}^{\log 1/|\epsilon|-C}
\left[ \int_{|x|<Cv^\gamma} \tilde{H}_\alpha(x) \frac{dx}{v^\gamma}+\int_{|x|>Cv^\gamma}\tilde{H}_\alpha(x) \frac{|v|^\gamma}{|x|^2}dx\right]\\
\eea
We split the integral
\bea
\sum_v h_{\alpha,n}& \sim & \int_{|x|<(\log 1/\delta)^\gamma}\tilde{H}_\alpha(x) \sum_{v=\log 1/\delta}^{\log 1/|\epsilon|-C}
\frac{dx}{v^\gamma}\\
& + & \int_{(\log 1/\delta)^\gamma<|x|<(\log 1/|\epsilon|)^\gamma}
\tilde{H}_\alpha(x) \sum_{v=|x|^{1/\gamma}}^{\log 1/|\epsilon}\frac{dx}{v^\gamma}\\
& + & \int_{(\log 1/\delta)^\gamma<|x|<(\log 1/|\epsilon|)^\gamma}
\tilde{H}_\alpha(x) \frac{1}{x^2}\sum_{v= \log 1/\delta}^{|x|^{1/\gamma}}v^\gamma dx\\
& + & 
\int_{|x|>(\log 1/|\epsilon|)^\gamma}
\tilde{H}_\alpha(x) \frac{1}{x^2}\sum_{v= \log 1/\delta}^{\log 1/|\epsilon}v^\gamma dx.\\
\eea
We estimate the quantities under $\Sigma$ and we get:
\bea
& \sim & \int_{|x|<(\log 1/\delta)^\gamma}\tilde{H}_\alpha(x) 
\frac{1}{(\log 1/\delta)^{\gamma-1}}\\
& + & \int_{(\log 1/\delta)^\gamma<|x|<(\log 1/|\epsilon|)^\gamma}
\tilde{H}_\alpha(x) \frac{1}{|x|^{1-1/\gamma}}dx\\
& + & \int_{(\log 1/\delta)^\gamma<|x|<(\log 1/|\epsilon|)^\gamma}
\tilde{H}_\alpha(x) \frac{1}{x^2} \frac{1}{|x|^{1-1/\gamma}}dx\\
& +& 
\int_{|x|>(\log 1/|\epsilon|)^\gamma}
\tilde{H}_\alpha(x) \frac{1}{x^2}(\log 1/|\epsilon|)^{\gamma+1} dx.\\
\eea

Using Lemma 12 this gives:

\bea
\sum_v h_{\alpha,n}& \sim & \int_{|x|<(\log 1/\delta)^\gamma}\tilde{H}_\alpha(x) |x|^{1/\gamma-1}
\left(\frac{|x|}{(\log 1/\delta)^{\gamma}}\right)^{1-1/\gamma}\\
& + & \int_{(\log 1/\delta)^\gamma<|x|<(\log 1/|\epsilon|)^\gamma}
\tilde{H}_\alpha(x) |x|^{1/\gamma-1}dx\\
& + & \int_{|x|>(\log 1/|\epsilon|)^\gamma}
\tilde{H}_\alpha(x) |x|^{1/\gamma-1} \left(\frac{(\log 1/|\epsilon|)^{\gamma}}{|x|}\right)^{1+1/\gamma} dx\\
& \rightarrow & 0, \delta \rightarrow 0.\\
\eea

Observe that we had to take $\delta$ small.

This finishes the case (i), $n<sp.$ So we have proved:

\begin{lemma}
The contribution to the geometric wedge product from $R_{1A}$
 in case (i),  $a \neq 0, n<sp$ goes to zero when $\delta \rightarrow 0.$
\end{lemma}

We next deal with the case $n>sp.$ 
Recall that:

Case (ii) $a \neq 0, n>sp$\\ We then have:

\bea
(u+iv)^\gamma & = & U+iV\\
& \sim & n^\gamma+ip(n)n^{\gamma-1}.\\
\eea
Then
\bea
H_{\alpha,n} & \sim & \int_{|x-n^\gamma|\leq p(n)n^{\gamma-1}}\tilde{H}_\alpha(x)\frac{1}{p(n)n^{\gamma-1}}dx\\
& + & \int_{n^\gamma/2>|x-n^\gamma|>p(n)n^{\gamma-1}}\tilde{H}_\alpha(x)\frac{p(n)n^{\gamma-1}}{|x-n^\gamma|^2}dx\\
& + &\int_{n^\gamma/2<|x-n^\gamma|<2n^{\gamma}}\tilde{H}_\alpha(x)\frac{p(n)n^{\gamma-1}}{n^{2\gamma}}dx\\
& + & \int_{|x-n^\gamma|>2n^\gamma}\tilde{H}_\alpha(x)\frac{p(n)n^{\gamma-1}}{x^2}dx\\
& = & I+II+III+IV\\
& = & I_n+II_n+III_n+IV_n.\\
\eea

For simplicity of notation we assume $a>0.$ Then we have the following range for $n$ from Lemma 20.
The number $n$ satisfies:

\bea
a \log 1/\delta <  \log |(\beta/\alpha)-1|+\log (1/|\epsilon|)-2nb\pi  & < & a \log 1/|\epsilon|-aC.\\
\eea
This gives
\bea
a \log 1/\delta - \log |(\beta/\alpha)-1|-\log (1/|\epsilon|)& < & -2nb\pi <\\
  - \log |(\beta/\alpha)-1|-\log (1/|\epsilon|) +a \log 1/|\epsilon|-aC. & &\\
  \eea
 So
 \bea
-a \log 1/\delta + \log |(\beta/\alpha)-1|+\log (1/|\epsilon|)& > & 2nb\pi >\\
  \log |(\beta/\alpha)-1|+\log (1/|\epsilon|) -a \log 1/|\epsilon|-aC. & & \\
  \eea
 Hence
 \bea
\frac{ \log |(\beta/\alpha)-1|+(1-a)\log (1/|\epsilon|) -aC}{2b\pi}& < & n <\\
\frac{-a \log 1/\delta + \log |(\beta/\alpha)-1|+\log (1/|\epsilon|)}{2b\pi}.& & \\
\eea

However, $n$ is further restricted because $n>sp$ and $p>\log 1/\delta.$
If we then estimate $IV$ and sum over $n$, we get

\bea
\sum_n IV_n & < \sim & \int_{|x|> (\log 1/\delta)^\gamma} \tilde{H}_\alpha(x)
\frac{1}{x^2} \sum_{n=\log 1/\delta}^{|x|^{1/\gamma}} n^\gamma\\
& < \sim & \int_{|x|> (\log 1/\delta)^\gamma} \tilde{H}_\alpha(x)
|x|^{1/\gamma-1}dx\\
& \rightarrow & 0\\
\eea

Similarly for $\sum_n III_n$ we get to estimate $\sum 1/n^{\gamma} <\sim |x|^{1/\gamma-1}$
which again is fine.

\medskip

Next we handle the terms $II_n.$ For a given $x$, the range of $n $ is on the order of
$2/3 |x|^{1/\gamma}<n< |x|^{1/\gamma}-p(x^{1/\gamma})$ and similarly for $n>|x|^{1/\gamma}.$
Also note that the terms $p(n)  \lesssim |x|^{1/\gamma}$ since $n \sim |x|^{1/\gamma}$ and $p\lesssim n.$
So we sum the expressions $\frac{n^{\gamma-1}}{(x-n^\gamma)^2}$ which integrates
to $\frac{1}{|x-n^\gamma|}$, so inserting the limits of the summation, we get 
a bound of the same form as for $III.$

Finally we sum over the $I_n$. Here we make the rough estimate that
$\log 1/\delta <p(n)<sn.$ So we integrate over $|x-n^\gamma|<sn^\gamma$ but in the
integrand we replace $p(n)$ by $\log 1/\delta.$  With this estimate we get the
integral $\tilde{H}_\alpha(x) \frac{1}{\log( 1/\delta )|x|}<<\tilde{H}_\alpha(x) |x|^{1/\gamma-1}.$ Hence
this also goes to zero with $\delta.$

Hence we have shown the following:

\begin{lemma} The contribution to the geometric wedge product in
the case of $R_{1A}$, case (ii), $a \neq 0, n>sp$ goes to zero when $\delta \rightarrow 0.$
\end{lemma}

\section{Theorem 7 for $R_{1C},$ the diagonal part of $R_1$}

We are in the set $\{C|\epsilon|<|z|,|w|<\delta, |z| \sim |w|\}.$

On $L_{\alpha,n},$ we have the following estimate for $u,v.$

\bea
2n\pi & < & u <2(n+1)\pi\\
|v-\frac{2n}{1-a}| & < & C'',\\
\log 1/\delta & < & v < \log (1/|\epsilon|)-C. \\
\eea
In the $U,V$ coordinates,
\bea
(u+iv)^\gamma & = & U+iV\\
V & \sim & |n|^\gamma \\
|U| & <\sim & |n|^\gamma.\\
\eea
So at intersection points
\bea
h_{\alpha,n} & \sim & \int \tilde{H}_\alpha (x) \frac{n^\gamma}{n^{2\gamma}+(x-U)^2}dx\\
& \sim & \int_{|x|<2n^\gamma}\tilde{H}_\alpha(x) \frac{dx}{n^\gamma}
+\int_{|x|>2n^\gamma}\tilde{H}_\alpha(x) \frac{n^\gamma dx}{x^2}.\\
\eea

Adding up the contributions
\bea
\sum_n h_{\alpha,n} & \sim &
\int_{|x|<(\log (1/|\delta))^\gamma} \tilde{H}_\alpha(x) \left( \sum_{n= \log 1/\delta}^\infty \frac{1}{n^\gamma}\right) dx\\
& + &  \int_{|x|>(\log (1/|\delta))^\gamma} \tilde{H}_\alpha(x) \left( \sum_{n= x^{1/\gamma}}
 \frac{1}{n^\gamma}\right) dx\\
& + & \int_{|x|>(\log (1/|\delta))^\gamma} \tilde{H}_\alpha(x) \left( \sum_{n= \log 1/\delta}^{x^{1/\gamma}} \frac{n^\gamma}{x^2} \right) dx.\\
\eea

After estimating the sums we get
\bea
\sum_n h_{\alpha,n} & \sim &
 \int_{|x|<(\log (1/|\delta))^\gamma} \tilde{H}_\alpha(x) \frac{1}{(\log (1/\delta))^\gamma} dx\\
& + & \int_{|x|>(\log (1/|\delta))^\gamma} \tilde{H}_\alpha(x) 
\frac{1}{(x^{1/\gamma})^{\gamma-1}}dx\\
& + & \int_{|x|>(\log (1/|\delta))^\gamma} \tilde{H}_\alpha(x) 
\frac{(x^{1/\gamma})^{\gamma+1}}{x^2}dx\\
\eea

So,

\bea
\sum_n h_{\alpha,n} & \sim & \int_{|x|>(\log (1/\delta))^\gamma}\tilde{H}_\alpha(x)|x|^{1/\gamma-1}dx\\
& + & 
\int_{|x|<(\log (1/\delta))^\gamma}\tilde{H}_\alpha(x)|x|^{1/\gamma-1}
\left(\frac{|x|}{(\log (1/\delta))^\gamma}\right)^{1-1/\gamma}dx\\
\eea

This is arbitrarily small as long as $\delta$ is chosen small enough.

\section{Theorem 7 for $R_2$, the part of $D_3$ close to the $z-$ axis}

This case is divided in two subcases depending on whether one is close
to one of the indicatrices ($R_{2A}$) or not ($R_{2B}$).

\section{Theorem 7 for $R_{2A}$ close to an indicatrix}

Again we assume that $a \neq 0.$
There are two indicatrices, $w=0$ and $w$ close to $\beta(\epsilon).$ By symmetry it suffices
to do one of them. We choose to estimate close to the indicatrix $w=0.$
So we set $R_{2A}=\{C|\epsilon|<|z|<\delta, |w|<s|\epsilon|\}$ for some small constant $s>0.$
Let $L_{\beta,m}^\epsilon$ and $L_{\alpha,n}$ be plaques intersecting at $(z,w)$ in $R_{2A}$
for parameters $(u',v'), (u,v).$

Since the point $(z,w)$ is about distance $|\beta'(0)||\epsilon|$ away from the
indicatrix for the perturbed lamination, we get ($w'= \beta(\epsilon)+\beta
e^{i \lambda (u'+(\log |\beta|/b)+iv')}+ \cdots$). This gives:

\bea
2m\pi & < & u' <2(m+1)\pi\\
C_1 & < & av'+2mb\pi +\log|\epsilon| < C_2.\\
\eea
We also have
\bea
C|\epsilon|<|z| =e^{-v} & = & |z'|=|\alpha(\epsilon)+e^{i(u'+\log |\beta|/b)-v'}+ \cdots|,\\
\eea
hence
\bea
C_3 & < &  v-v' <C_4\\
C_4 & < & av+2mb\pi +\log|\epsilon| < C_5.\\
2n\pi & < & u <2(n+1)\pi.\\
\eea
 Using
 \bea
|w| & < & s|\epsilon|\\
\eea
we get
\bea
e^{-bu-av} & < & s |\epsilon|\\
\log (1/s) & < & av+2nb\pi +\log |\epsilon| \\
2(n-m)b \pi & = &( av+2nb\pi+\log |\epsilon|)-(av+2m\pi b+\log |\epsilon|)\\
& > & \log (1/s)-C_1.\\
\eea

These calculations show that for the given plaques, the pairs $(u,v), (u',v')$
belong to rectangles of uniformly bounded size. Hence the number of intersection
points can easily be estimated by using slope estimates for the plaques.
We get a uniformly bounded number of intersection points.

We divide this into cases $I,II,III$.\\

For $I$, we have $1/C \log (1/|\epsilon|)<2mb\pi+\log|\epsilon|<C \log ( 1/|\epsilon|)$.\\

For $II$ we have $2mb\pi+\log|\epsilon|<1/C \log ( 1/|\epsilon|)$.\\

For $III$ we have $2mb\pi+\log|\epsilon|>C \log ( 1/|\epsilon|)$.
 We note however, that in case $III$, $v'$ must be very large in comparison with $\log 1/|\epsilon|.$
 This implies that $|z'| << |\epsilon|$ hence there are no intersection points in this case.
 So we are left with the two cases $R_{2AI}, R_{2AII}.$

\section{Theorem 7 for $R_{2AI}$ close to an indicatrix.}

It follows in this case that $v,v' \sim \log (1/|\epsilon|).$ Hence

\bea
u'+iv' & \sim & 2m\pi+i \log (1/|\epsilon|)\\
\eea
and
\bea
U'+iV' & \sim & U'+i (\log (1/|\epsilon|))^\gamma.\\
\eea
In particular
\bea
|U'| & < \sim & (\log (1/|\epsilon|))^\gamma.\\
\eea

Using the Poisson integral we estimate
\bea
h^\epsilon_{\beta,m} & \sim & \int \tilde{H}_\beta(y) \frac{(\log(1/|\epsilon|))^\gamma}
{(\log (1/|\epsilon|))^{2\gamma} +(x-U')^2}dy\\
& \sim & \int_{|y|<2(\log (1/|\epsilon|))^\gamma} \tilde{H}_\beta(y) \frac{1}{(\log(1/|\epsilon|))^\gamma}dy\\
& + & \int _{|y|>2(\log (1/|\epsilon|))^\gamma} \tilde{H}_\beta(y) \frac{(\log(1/|\epsilon|))^\gamma}
 {y^2 }dy.\\
 \eea
 Adding up
 \bea
 \sum_{m\in I} h^\epsilon_{\beta,m} & \sim &
 \int_{|y|<2(\log (1/|\epsilon|))^\gamma} \tilde{H}_\beta(y) |y|^{1/\gamma-1}
 \left(\frac{|y|}{(\log(1/|\epsilon|))^\gamma}\right)^{1-1/\gamma}dy\\
& + & \int _{|y|>2(\log (1/|\epsilon|))^\gamma} \tilde{H}_\beta(y) |y|^{1/\gamma-1}
\left(\frac{(\log(1/|\epsilon|))^\gamma}{|y|}\right)^{1/\gamma+1}dy.\\
\eea

Next we estimate $h_{\alpha,n}$. There are two cases to consider:\\

a): $n<C \log (1/|\epsilon|)$\\

b): $n>C \log (1/|\epsilon|)$\\

The contribution for case a) is:\\

Case $R_{2AIa}:$

Recall that we have $n>m-C_6.$ Hence we have that $|n|< C \log (1/|\epsilon|).$ This means
that we can write $u+iv \sim 2n\pi+i (\log (1/|\epsilon|)) $. Hence the estimates work as for
$h^\epsilon_{\beta,m}.$

\bea
\sum_{|n|<C\log (1/|\epsilon|)} h_{\alpha,n} & \sim &
 \int_{|x|<2(\log (1/|\epsilon|))^\gamma} \tilde{H}_\alpha(x) |x|^{1/\gamma-1}
 \left(\frac{|x|}{(\log(1/|\epsilon|))^\gamma}\right)^{1-1/\gamma}dx\\
& + & \int _{|x|>2(\log (1/|\epsilon|))^\gamma} \tilde{H}_\alpha(x) |x|^{1/\gamma-1}
\left(\frac{(\log(1/|\epsilon|))^\gamma}{|x|}\right)^{1/\gamma+1}dx.\\
\eea

Case $R_{2AIb}$. We have

\bea
u+iv & \sim & n+i \log (1/|\epsilon|)\\
U+ i V & \sim & n^\gamma +in^{\gamma-1}\log (1/|\epsilon|),\;{\mbox{and}}\\
h_{\alpha,n} & \sim & \int \tilde{H}_\alpha(x) \frac{n^{\gamma-1} \log (1/|\epsilon|)}
{(n^{\gamma-1}\log(1/|\epsilon|))^2+(x-n^\gamma)^2}dx.\\
\eea

This integral has already been estimated. See the calculations for the set $D_1$ in the
region where $|z-\eta|<d|\eta|$, case (ii) where $n>10 \log (1/|\eta|).$ It follows that the contributions from that region goes to zero with $\epsilon.$

\section{Theorem 7 for $R_{2AII}$ close to an indicatrix.}

We restrict for simplicity to the case $a>0$. 
We can divide into three cases:

\medskip

\noindent a)  $n>m>v,v'$\\
b) $ n>v,v'>m$\\
c) $ v,v'>n>m$\\

\section{Theorem 7 for $R_{2AIIa}$ close to an indicatrix.}

We have 
\bea
(u+iv)^\gamma & = & U+iV\\
& \sim & n^\gamma+ivn^{\gamma-1}\\
(u'+iv')^\gamma & = & U'+iV'\\
& \sim & m^\gamma+iv'm^{\gamma-1}\\
& \sim & (\log 1/|\epsilon|)^\gamma+ iv' (\log (1/|\epsilon|))^{\gamma-1}\\
\log 1/\delta & < & v' < \log 1/|\epsilon|.\\
\eea

We now estimate

\bea
H_\beta & \sim & \int \tilde{H}_\beta (y) \frac{v' (\log (1/|\epsilon|))^{\gamma-1}}{[v' (\log (1/|\epsilon|))^{\gamma-1}]^2+(y-m^\gamma)^2}dy.\\
\eea

We divide the integral and estimate each term.

\bea
H_\beta & \sim & \int_{|y-m^\gamma|<c v' (\log 1/|\epsilon|)^{\gamma-1}}
 \tilde{H}_\beta (y) \frac{1}{v' (\log (1/|\epsilon|))^{\gamma-1}}dy\\
& + & \int_ {(\log 1/|\epsilon|)^\gamma/2>|y-(\log 1/|\epsilon|)^\gamma|>c v' (\log 1/|\epsilon|)^{\gamma-1}}\tilde{H}_\beta (y) \frac{v' (\log (1/|\epsilon|))^{\gamma-1}}{(y-(\log 1/|\epsilon|)^\gamma)^2}dy\\
& + & \int_ {|y-(\log 1/|\epsilon|)^\gamma|>(\log 1/|\epsilon|)^{\gamma}/2}\tilde{H}_\beta (y) \frac{v' (\log (1/|\epsilon|))^{\gamma-1}}{(y-(\log 1/|\epsilon|)^\gamma)^2}dy.\\
\eea
So
\bea
H_\beta & \sim & \int_{|y-m^\gamma|<c v' (\log 1/|\epsilon|)^{\gamma-1}}
 \tilde{H}_\beta (y) \frac{y^{1/\gamma-1}}{v' }dy\\
& + & \int_ {|y-(\log 1/|\epsilon|)^\gamma|>c v' (\log 1/|\epsilon|)^{\gamma-1}}\tilde{H}_\beta (y) \frac{v' (\log (1/|\epsilon|))^{\gamma-1}}{(y-(\log 1/|\epsilon|)^\gamma)^2}dy\\
& = & H_{\beta_{1,v'}}+H_{\beta_{2,v'}}.\\
\eea
For $H_\alpha$ we have
\bea
H_\alpha & \sim & \int \tilde{H}_\alpha (x) \frac{v n^{\gamma-1}}{[v n^{\gamma-1}]^2+  (x-n^\gamma)^2}dx\\
& \sim & \int_{|x-n^\gamma|<c v n^{\gamma-1}}
 \tilde{H}_\alpha (x) \frac{1}{v n^{\gamma-1}}dx\\
& + & \int_ {|x-n^\gamma|>c v n^{\gamma-1}}\tilde{H}_\alpha (x) \frac{v n^{\gamma-1}}{(x-n^\gamma)^2}dx.\\
\eea

To sum up over the intersection points, we note at first that for a given plaque
$L_{\beta,m}$ there is a finite range of $v'$ and $v-v'$ is bounded, so we can
assume that there is one intersection point with $L_{\alpha,n}$ for each $n>m$.
Hence we sum first over the plaques $L_{\alpha,n}$, $m<n<\infty.$

\bea
\sum_n \int_{|x-n^\gamma|<c v n^{\gamma-1}}
 \tilde{H}_\alpha (x) \frac{1}{v n^{\gamma-1}}dx & \sim & 
 \sum_{n=x^{1/\gamma}-v}^{n=x^{1/\gamma}+v} \int \\
 & \sim &  \int_{x>m^\gamma} \tilde{H}_\alpha(x) |x|^{1/\gamma-1}dx.\\
 \eea
 
 The other contribution is 
\bea
 \sum_n\int_ {|x-n^\gamma|>c v n^{\gamma-1}}\tilde{H}_\alpha (x) \frac{v n^{\gamma-1}}{(x-n^\gamma)^2}dx\\
 & \sim & \int_{x>m^\gamma-cvm^{\gamma-1}}\tilde{H}_\alpha(x)|x|^{1/\gamma-1}\\
& + & \int_{x<m^\gamma-cvm^{\gamma-1}}\tilde{H}_\alpha(x)\frac{v}{|x-m^\gamma|}dx,\\
 \eea
 
 so we conclude:
 
 \bea
 \sum_{n>m}H_\alpha & \sim & \int_{x>m^\gamma} \tilde{H}(x)|x|^{1/\gamma-1}+
 \int_{x<m^\gamma-cvm^{\gamma-1}}\tilde{H}_\alpha(x)\frac{v}{|x-m^\gamma|}dx\\
 & < \sim & \int_{|x|>m^\gamma/2}\tilde{H}_\alpha |x|^{1/\gamma-1}+
 \int_{|x|<m^\gamma/2}\tilde{H}_\alpha(x)|x|^{1/\gamma-1}\left(\frac{|x|}{m^\gamma}\right)^{1-1/\gamma}dx.\\
 \eea

In this case $m$ will have approximately the range $(\log 1/|\epsilon|)/2<m<\log 1/|\epsilon|$, hence 
we have

\bea
\sum_{n>m}H_\alpha & <\sim & 
  \int_{|x|>(\log 1/|\epsilon|)^\gamma}\tilde{H}_\alpha |x|^{1/\gamma-1}\\
  & + & 
 \int_{|x|< (\log 1/|\epsilon|)^\gamma}\tilde{H}_\alpha(x)|x|^{1/\gamma-1}\left(\frac{|x|}{(\log 1/|\epsilon|)^\gamma}\right)^{1-1/\gamma}dx.\\
 \eea

Next we sum $H_\beta$ over $m$ or equivalently over $v', \log 1/\delta < v'< (\log1/|\epsilon|)/2.$
 We integrate first over $H_{\beta_{1,v'}}$. For a given $y$, the range of $v'$ is
 in the interval with endpoints $(1\pm c) \frac{y-(\log 1/|\epsilon|)^\gamma}{(\log 1/|\epsilon|)^{\gamma-1}}.$ This part is bounded by
 
 \bea
 \int_{|y-(\log 1/|\epsilon|)^\gamma|<(\log 1/|\epsilon|)^\gamma/2}
 \tilde{H}_\beta(y)|y|^{1/\gamma-1} dy & \rightarrow & 0.\\
 \eea
 The second part is bounded by
  \bea
& & \int_{|y|<2(\log 1/|\epsilon|)^\gamma}\tilde{H}_\beta (y) |y|^{1/\gamma-1}
\left(\frac{|y|}{(\log 1/|\epsilon|)^\gamma}\right)^{1-1/\gamma}dy\\
 & + & \int_{|y|>2(\log 1/|\epsilon|)^\gamma}\tilde{H}_\beta(y) |y|^{1/\gamma-1}
\left( \frac{(\log 1/|\epsilon|)^\gamma}{|y|}\right)^{1+1/\gamma}dy.\\
 \eea
 
Again the contribution goes to zero by Proposition 1 and Lemma 12.

\section{Theorem 7 for $R_{2AIIb}$ close to an indicatrix.}

In this case $n>v,v'>m.$
First we recall the estimates for $H_\alpha$ which are the same as in the case $R_{2AIIa}.$
We have

\bea
(u+iv)^\gamma & = & U+iV\\
& \sim & n^\gamma+ivn^{\gamma-1}\;{\mbox{with}}\\
\log 1/\delta & < & v,v' < \log 1/|\epsilon|.\\
\eea
So
\bea
H_\alpha & \sim &   \int \tilde{H}_\alpha (x) \frac{v n^{\gamma-1}}{[v n^{\gamma-1}]^2+(x-n^\gamma)^2}dx\\
& \sim & \int_{|x-n^\gamma|<c v n^{\gamma-1}}
 \tilde{H}_\alpha (x) \frac{1}{v n^{\gamma-1}}dx\\
& + & \int_ {|x-n^\gamma|>c v n^{\gamma-1}}\tilde{H}_\alpha (x) \frac{v n^{\gamma-1}}{(x-n^\gamma)^2}dx.\\
\eea

Next we estimate $H_\beta.$ We have

\bea
(u'+iv')^\gamma & = & U'+iV'\\
(\log 1/|\epsilon|)/2 & < & v' < \log 1/|\epsilon|\\
m+v' &  = & \log 1/|\epsilon|\\
V'  & \sim & (\log 1/|\epsilon|)^\gamma\\
|U'| & < \sim & (\log 1/|\epsilon|)^\gamma.\\
\eea
Hence
\bea
H_\beta & \sim & \int_{|y|<2 (\log 1/|\epsilon|)^\gamma} \tilde{H}_\beta \frac{1}{(\log 1/|\epsilon|)^\gamma}dy\\
& + & \int_{|y|>2(\log 1/|\epsilon|)^\gamma}\tilde{H}_\beta \frac{(\log 1/|\epsilon|)^\gamma}{y^2}dy.\\
\eea

Next we estimate the contribution to the geometric wedge product.
So fix $\alpha,\beta.$ Next fix a plaque $L_{\beta,m}$, $ v,v' \sim \log 1/|\epsilon|-m.$
Next we consider the contribution from $H_\alpha$ for all $n>v.$
This is the same estimate as in the previous section, so goes to zero
when $\epsilon \rightarrow 0.$ To sum up over $m,$ notice that we have about $\log 1/|\epsilon|$
terms of the same order of magnitude. From this we get that the contribution goes to zero when
$\epsilon \rightarrow 0.$

To estimate the geometric wedge product, we sum independently over
$n,m$ throwing out the condition that $n>m$. We get as in the previous
section that the contribution goes to zero.

\section{Theorem 7 for $R_{2AIIc}$ close to an indicatrix.}

Here we deal with the case when

$v,v'>n>m.$ In this case the same formula as in the last section applies to both
$H_\alpha$ and $H_\beta.$ We have:

\bea
H_\alpha & \sim & \int_{|x|<2 (\log 1/|\epsilon|)^\gamma} \tilde{H}_\alpha \frac{1}{(\log 1/|\epsilon|)^\gamma}dx\\
& + & \int_{|x|>2(\log 1/|\epsilon|)^\gamma}\tilde{H}_\alpha \frac{(\log 1/|\epsilon|)^\gamma}{x^2}dx,\;{\mbox{and}}\\
H_\beta & \sim & \int_{|y|<2 (\log 1/|\epsilon|)^\gamma} \tilde{H}_\beta \frac{1}{(\log 1/|\epsilon|)^\gamma}dy\\
& + & \int_{|y|>2(\log 1/|\epsilon|)^\gamma}\tilde{H}_\beta \frac{(\log 1/|\epsilon|)^\gamma}{y^2}dy.\\
\eea

So again the contribution goes to zero.

\section{Theorem 7 for $R_{2B}$ away from the indicatrices.}

At an intersection point $p=(z,w)$ of $L_{\alpha,n}, L^\epsilon_{\beta,m}$ we have

\bea
s|\epsilon| & < & |w|< C|\epsilon|\\
s|\epsilon| & < & |w-\beta(\epsilon)| <C |\epsilon|.\\
\eea
So
\bea
\log |\epsilon|-C & < & -av-bu < \log |\epsilon|+C\\
\log |\epsilon|-C & < & -av'-bu' < \log |\epsilon|+C.\\
\eea
This gives:
\bea
-C & < & v-v' < C\\
-C & < & n-m < C\\
\log (1/\delta) & < & v,v' < \log (1/|\epsilon|)-C\\
-C \log (1/|\epsilon|) & < & u,u',n,m < C \log (1/|\epsilon|).\\
\eea

Given $(\alpha,\beta,n,m)$ we need to estimate the values of $v,v'$ corresponding to
an intersection, as well as the number of intersections. The following is immediate.
There is no dependence on $\alpha,\beta$.

\begin{lemma}
At intersection points of $L_{\alpha,n},L^\epsilon_{\beta,m}$ in $R_{2B}$ away from
the indicatrices, we have 
$$
-2nb\pi/a+1/a \log (1/|\epsilon|)-C<v,v'<-2nb\pi/a+1/a \log (1/|\epsilon|)+C.
$$
\end{lemma}

It follows that intersection points are localized in bounded rectangles. To
show finiteness of number of intersection points for given plaques,
we use slope estimates.

We divide the estimates in two cases, (i) if $v,v' \sim \log (1/|\epsilon|)$ and
(ii) if $\log (1/\delta)<v,v' < 1/C \log (1/|\epsilon|).$

\section{Theorem 7 for $R_{2Bi}$ when $v \sim \log (1/|\epsilon|)$}

Recall that this means that for a large constant $A$,
$\frac{1}{A} \log \frac{1}{|\epsilon|}<v< A \log \frac{1}{|\epsilon|}.$
The estimates for $h_{\alpha,n}$ and $h^\epsilon_{\beta,m}$ are similiar.
We have 
\bea
U+iV & = & (u+iv)^\gamma\\
& \sim & U+i (\log (1/|\epsilon|))^\gamma\\
|U| & <\sim & (\log (1/|\epsilon|))^\gamma,\\
\eea
and at intersection points
\bea
h_{\alpha,n} & \sim & \int \tilde{H}_\alpha(x) \frac{(\log (1/|\epsilon|))^\gamma}{(\log (1/|\epsilon|))^{2\gamma}+(x-U)^2}dx\\
& \sim & 
\int_{|x|<C(\log (1/|\epsilon|))^\gamma}\tilde{H}_\alpha(x) \frac{1}{(\log (1/|\epsilon|))^\gamma}dx\\
& + & \int {|x|>C(\log (1/|\epsilon|))^\gamma} \tilde{H}_\alpha(x) \frac{(\log (1/|\epsilon|))^\gamma}
{x^2}dx\\
& \sim & 
\int_{|x|<C(\log (1/|\epsilon|))^\gamma}\tilde{H}_\alpha|x|^{1/\gamma-1}
 \left(\frac{|x|}{(\log (1/|\epsilon|))^\gamma}\right)^{1-1/\gamma} \frac{1}{\log (1/|\epsilon|)}dx\\
& + & \int_{|x|>C(\log (1/|\epsilon|))^\gamma} \tilde{H}_\alpha |x|^{1/\gamma-1}
\left(\frac{(\log (1/|\epsilon|))^\gamma}{|x|}\right)^{1+1/\gamma}\frac{1}{\log (1/|\epsilon|)}dx.\\
\eea
We estimate the total contribution,
\bea
\sum_n h_{\alpha,n} & \sim &
\int_{|x|<C(\log (1/|\epsilon|))^\gamma}\tilde{H}_\alpha|x|^{1/\gamma-1}
 \left(\frac{|x|}{(\log (1/|\epsilon|))^\gamma}\right)^{1-1/\gamma} dx\\
& + & \int_{|x|>C(\log (1/|\epsilon|))^\gamma} \tilde{H}_\alpha |x|^{1/\gamma-1}
\left(\frac{(\log (1/|\epsilon|))^\gamma}{|x|}\right)^{1+1/\gamma}dx\\
\eea
\noindent which will converge to $0$ by Proposition 1.

\section{Theorem 7 for $R_{2Bii}$ when  $v <\frac{1}{A} \log (1/|\epsilon|)$}

In this case we have $u,u',n,m \sim \log (1/|\epsilon|).$
The estimates for $h_{\alpha,n}, h^\epsilon_{\beta,m}$ are similar. In the following
$0<d<<1.$ More precisely, $d$ will be close to  $|a|A$, see the 4th inequality below.

\bea
(1-d)\log (1/|\epsilon|) & < & 2nb \pi < (1+d) \log (1/|\epsilon|)\\
\log |\epsilon|-C & < & -av-bu < \log |\epsilon|+C\\
 \log |\epsilon|+2nb\pi-C & < & -av < \log |\epsilon|+2bn\pi + C\\
 -d\log(1/ |\epsilon| )-C & < & -av < d \log (1/|\epsilon|)+C.\\
 \eea
 In $U,V$ coordinates:
 \bea
 U+iV & = & (u+iv)^\gamma\\
 & \sim & (\log (1/|\epsilon|))^\gamma + i (\log (1/|\epsilon|))^{\gamma-1} v.\\
 \eea
 This gives
 \bea
 h_{\alpha,n} & \sim & \int \tilde{H}_\alpha(x) \frac{ (\log (1/|\epsilon|))^{\gamma-1} v}
 { ((\log (1/|\epsilon|))^{\gamma-1} v)^2+(x-U)^2}dx.\\
\eea

When we sum up over $h_{\alpha,n}, h^\epsilon_{\beta,m}$ we can take for simplicity $n=m$
and $v=v'$ since $|n-m|, |v-v'|$ are uniformly bounded in $R_{2B}$ as stated above.
The product of contributions is estimated by

\bea
h_{\alpha,n}h^\epsilon_{\beta,m}
& \sim &
 \int \tilde{H}_\alpha(x) \frac{ (\log (1/|\epsilon|))^{\gamma-1} v}
 { ((\log (1/|\epsilon|))^{\gamma-1} v)^2+(x-U)^2}dx\\
 & * &  \int \tilde{H}_\beta(y) \frac{ (\log (1/|\epsilon|))^{\gamma-1} v}
 { ((\log (1/|\epsilon|))^{\gamma-1} v)^2+(y-U)^2}dy\\
 \eea
 So
 \bea
h_{\alpha,n}h^\epsilon_{\beta,m}& \sim & 
 [ \int_{|x-U|<(\log (1/|\epsilon|))^{\gamma-1}|v|}
   \tilde{H}_\alpha(x) \frac{1}{ (\log (1/|\epsilon|))^{\gamma-1} v}dx\\
& + &   \int_{|x-U|>(\log (1/|\epsilon|))^{\gamma-1}|v|}\tilde{H}\alpha(x) \frac{ (\log (1/|\epsilon|))^{\gamma-1} v}
 { (x-U)^2}dx]\\
& * & [ \int _{|y-U|<(\log (1/|\epsilon|))^{\gamma-1}|v|}\tilde{H}_\beta(y) \frac{1}{ (\log (1/|\epsilon|))^{\gamma-1} v}dy\\
& + &  \int _{|y-U|>(\log (1/|\epsilon|))^{\gamma-1}|v|} \tilde{H}_\beta(y) \frac{ (\log (1/|\epsilon|))^{\gamma-1} v}
 { (y-U)^2}dy]\\
 &  = & [I+II][III+IV].\\
\eea

There are 4 cases to sum over: $(I,III)$, $(II,III)$, $(II,IV)$ and $(I,IV)$.
The case $(I,IV)$ is similar to $(II,III)$ so we can skip it without any loss.

\section{Theorem 7 for $R_{2Bii(I,III)}$}
 We have
\bea
h_{\alpha,n}h^\epsilon_{\beta,m}
& \sim &
  \int_{|x-U|<(\log (1/|\epsilon|))^{\gamma-1}|v|}
   \tilde{H}_\alpha(x) \frac{1}{ (\log (1/|\epsilon|))^{\gamma-1} v}dx\\
& * &  \int _{|y-U|<(\log (1/|\epsilon|))^{\gamma-1}|v|}\tilde{H}_\beta(y) \frac{1}{ (\log (1/|\epsilon|))^{\gamma-1} v}dy\\
& \lesssim & \frac{1}{v^2}  \int_{|x-U|<1/C (\log (1/|\epsilon|))^{\gamma}}
   \tilde{H}_\alpha(x) \frac{1}{ (\log (1/|\epsilon|))^{\gamma-1} }dx\\
& * &  \int _{|y-U|<1/C(\log (1/|\epsilon|))^{\gamma}}\tilde{H}_\beta(y) \frac{1}{ (\log (1/|\epsilon|))^{\gamma-1} }dy\\
& \sim & \frac{1}{v^2}  \int_{|x-U|<1/C (\log (1/|\epsilon|))^{\gamma}}
   \tilde{H}_\alpha(x) |x|^{1/\gamma-1}dx\\
& * &  \int _{|y-U|<1/C(\log (1/|\epsilon|))^{\gamma}}\tilde{H}_\beta(y) |y|^{1/\gamma-1}dy.\\
\eea
Since
\bea
\log (1/\delta) & <  & v < 1/A \log (1/|\epsilon|)\\
\eea
we get
\bea
\sum h_{\alpha,n}h^\epsilon_{\beta,m} & \lesssim &
\int_{\log (1/\delta)}^{d \log (1/|\epsilon|)}\frac{1}{v^2}
 \int_{|x-U|<1/C (\log (1/|\epsilon|))^{\gamma}}
   \tilde{H}_\alpha(x) |x|^{1/\gamma-1}dx\\
& * &  \int _{|y-U|<1/C(\log (1/|\epsilon|))^{\gamma}}\tilde{H}_\beta(y) |y|^{1/\gamma-1}dy,\;{\mbox{or}}\\
\sum h_{\alpha,n}h^\epsilon_{\beta,m} & \lesssim &
 \frac{1}{\log (1/\delta)}\int_{|x-U|<1/C (\log (1/|\epsilon|))^{\gamma}}
   \tilde{H}_\alpha(x) |x|^{1/\gamma-1}dx\\
& * &  \int _{|y-U|<1/C(\log (1/|\epsilon|))^{\gamma}}\tilde{H}_\beta(y) |y|^{1/\gamma-1}dy.\\
\eea
Finally
\bea
\sum h_{\alpha,n}h^\epsilon_{\beta,m} & \lesssim &
 \frac{1}{\log (1/\delta)}\int_{|x-(\log (1/|\epsilon|))^\gamma|<1/C (\log (1/|\epsilon|))^{\gamma}}
   \tilde{H}_\alpha(x) |x|^{1/\gamma-1}dx\\
& * &  \int _{|y-(\log (1/|\epsilon|))^\gamma|<1/C (\log (1/|\epsilon|))^{\gamma}}\tilde{H}_\beta(y) |y|^{1/\gamma-1}dy.\\
\eea

This contribution goes to zero when $\epsilon \rightarrow 0.$

\section{Theorem 7 for $R_{2Bii(II,III)}$}

We estimate
\bea
h_{\alpha,n}h^\epsilon_{\beta,m}
& \sim &
  \int_{|x-U|>(\log (1/|\epsilon|))^{\gamma-1}|v|}\tilde{H}\alpha(x) \frac{ (\log (1/|\epsilon|))^{\gamma-1} v}
 { (x-U)^2}dx\\
& * &  \int _{|y-U|<(\log (1/|\epsilon|))^{\gamma-1}|v|}\tilde{H}_\beta(y) \frac{1}{ (\log (1/|\epsilon|))^{\gamma-1} v}dy.\\
\eea

Here $\log (1/\delta) < v< d \log (1/|\epsilon|),0<d<<1$ and
$ -av = \log |\epsilon|+ 2bn \pi+{\mathcal {O}}(1).$  Also we can take $n=m.$ When we sum 
over $n$, $v$ runs through  $\log (1/\delta) < v< d \log (1/|\epsilon|)$.
Hence the contribution to the geometric wedge product is

\bea
\sum _{n,m} h_{\alpha,n}h^\epsilon_{\beta,m}
& \sim &
\sum_{v=\log (1/\delta)}^{d \log (1/|\epsilon|)}
  \int_{|x-U|>(\log (1/|\epsilon|))^{\gamma-1}|v|}\tilde{H}_\alpha(x) \frac{ 1}{(x-U)^2}dx\\
& * &  \int _{|y-U|<(\log (1/|\epsilon|))^{\gamma-1}|v|}\tilde{H}_\beta(y)dy\\
& \sim & 
\sum_{v=\log (1/\delta)}^{d \log (1/|\epsilon|)}
  \int_{|x-(\log (1/|\epsilon|))^\gamma|>(\log (1/|\epsilon|))^{\gamma-1}|v|}\tilde{H}_\alpha(x) \frac{ 1}{(x-(\log (1/|\epsilon|)^\gamma)^2}dx\\
& * &  \int _{|y-(\log (1/|\epsilon|))^\gamma|<(\log (1/|\epsilon|))^{\gamma-1}|v|}\tilde{H}_\beta(y)dy\\
\eea

We introduce a counting function, $N(x,y)$, which tells us for a given $(x,y)$
for how many terms of the sum $(x,y)$ is in the domain of integration for the above integrals:

\bea
|x-(\log (1/|\epsilon|))^\gamma| & > & (\log (1/|\epsilon|))^{\gamma-1}|v|\\
|y-(\log (1/|\epsilon|))^\gamma| & < & (\log (1/|\epsilon|))^{\gamma-1}|v|.\\
\eea

We divide the above domain.

\bea
P_1 & = & \{|x-(\log (1/|\epsilon|))^\gamma|>d (\log (1/|\epsilon|))^\gamma,\\
&  & 
|y-(\log (1/|\epsilon|))^\gamma|<\log (1/\delta) (\log (1/|\epsilon|))^{\gamma-1}\}\\
\eea
in which case
\bea
N_1(x,y) & \sim & d \log (1/|\epsilon|).\\
P_2 & = & \{|x-(\log (1/|\epsilon|))^\gamma|>d (\log (1/|\epsilon|))^\gamma,\\
&  & 
\log (1/\delta) (\log (1/|\epsilon|))^{\gamma-1}<|y-(\log (1/|\epsilon|))^\gamma|<
d (\log (1/|\epsilon|))^{\gamma}\},\\
\eea
for $P_2$
\bea
N_2(x,y) & \sim & \frac{d( \log (1/|\epsilon|))^\gamma-|y-( \log (1/|\epsilon|))^\gamma|}
{( \log (1/|\epsilon|))^{\gamma-1}}.\\
P_3 & = & \{\log (1/\delta) (\log (1/|\epsilon|))^{\gamma-1}<|x-(\log (1/|\epsilon|))^\gamma|<
d (\log (1/|\epsilon|))^{\gamma},\\
& & \log (1/\delta) (\log (1/|\epsilon|))^{\gamma-1}<|y-(\log (1/|\epsilon|))^\gamma|<
d (\log (1/|\epsilon|))^{\gamma}\},\\
\eea
for $P_3$
\bea
N_3(x,y) & \sim & \frac{|x-(\log (1/|\epsilon|))^\gamma|-|y-(\log (1/|\epsilon|))^\gamma|}
{( \log (1/|\epsilon|))^{\gamma-1}}\\
\eea
 when positive
 \bea
N_3(x,y) & \sim & \frac{|x-y|}
{( \log (1/|\epsilon|))^{\gamma-1}} .\\
\eea

\section{Theorem 7 for $R_{2Bii(II,III)P_1}$}

This gives the estimate for the product

\bea
\sum _{n,m} h_{\alpha,n}h^\epsilon_{\beta,m}& \sim &
d \log (1/|\epsilon|) \int_{P_1} \frac{\tilde{H}_\alpha(x)\tilde{H}_\beta(y)}{(x-(\log (1/|\epsilon|))^\gamma)^2}
dxdy\\
& \sim &  \int_{P_1} \frac{\tilde{H}_\alpha(x)|x|^{1/\gamma-1}|x|^{1-1/\gamma}\tilde{H}_\beta(y)|y|^{1/\gamma-1}}{(x-(\log (1/|\epsilon|))^\gamma)^2 ((\log (1/|\epsilon|)^\gamma)^{1/\gamma-1}}\log (1/|\epsilon|)\\
& \sim &  \int_{P_1} \frac{\tilde{H}_\alpha(x)|x|^{1/\gamma-1}|x|^{1-1/\gamma}\tilde{H}_\beta(y)|y|^{1/\gamma-1}}{|x-(\log (1/|\epsilon|))^\gamma|^{1-1/\gamma} ((\log (1/|\epsilon|)^\gamma)^{2/\gamma}}\log (1/|\epsilon|)\\
&  \lesssim & \int_{P_1} \tilde{H}_\alpha(x)|x|^{1/\gamma-1}\tilde{H}_\beta(y)|y|^{1/\gamma-1}
\frac{1}{\log (1/|\epsilon|)}\\
& \rightarrow & 0\\
\eea
\noindent when $\epsilon \rightarrow 0.$

\section{Theorem 7 for $R_{2Bii(II,III)P_2}$}

We get the estimate

\bea
\sum _{n,m} h_{\alpha,n}h^\epsilon_{\beta,m}& \sim &
 \int_{P_2} \frac{\tilde{H}_\alpha(x)\tilde{H}_\beta(y)}{(x-(\log (1/|\epsilon|))^\gamma)^2}
\frac{d (\log (1/|\epsilon|))^\gamma-|y-(\log (1/|\epsilon|))^\gamma|}{(\log (1/|\epsilon|))^{\gamma-1}}
dxdy\\
& \sim &
 \int_{P_2} \frac{\tilde{H}_\alpha(x)|x|^{1/\gamma-1}\tilde{H}_\beta(y)
 |y|^{1/\gamma-1}}{(x-(\log (1/|\epsilon|))^\gamma)^2}|x|^{1-1/\gamma}|y|^{1-1/\gamma}\\
 & * & 
\frac{d (\log (1/|\epsilon|))^\gamma-|y-(\log (1/|\epsilon|))^\gamma|}{(\log (1/|\epsilon|))^{\gamma-1}}
dxdy.\\
\eea

Using the definition of $P_2,$
\bea
\sum _{n,m} h_{\alpha,n}h^\epsilon_{\beta,m}
& \sim &
 \int_{P_2} \frac{\tilde{H}_\alpha(x)|x|^{1/\gamma-1}\tilde{H}_\beta(y)
 |y|^{1/\gamma-1}}{(x-(\log (1/|\epsilon|))^\gamma)^2}|x|^{1-1/\gamma}\\
 & * & 
(d (\log (1/|\epsilon|))^\gamma-|y-(\log (1/|\epsilon|))^\gamma|)
dxdy\\
& \sim &
 \int_{P_2} \tilde{H}_\alpha(x)|x|^{1/\gamma-1}\tilde{H}_\beta(y)
 |y|^{1/\gamma-1} \frac{|x|}{|x-(\log (1/|\epsilon|))^\gamma|}|x|^{-1/\gamma}\\
 & * & 
\frac{(d (\log (1/|\epsilon|))^\gamma-|y-(\log (1/|\epsilon|))^\gamma|)}
{|x-(\log (1/|\epsilon|))^\gamma|}dxdy\\
& <\sim & 
\int_{P_2} \tilde{H}_\alpha(x)|x|^{1/\gamma-1}\tilde{H}_\beta(y)
 |y|^{1/\gamma-1}\frac{1}{\log (1/|\epsilon|)}dxdy\\
& \rightarrow & 0\\
\eea
\noindent as $\epsilon \rightarrow 0.$

\section{Theorem 7 for $R_{2Bii(II,III)P_3}$}

We estimate, using the definition of $P_3.$

\bea
\sum _{n,m} h_{\alpha,n}h^\epsilon_{\beta,m}& \sim &
 \int_{P_3} \frac{\tilde{H}_\alpha(x)\tilde{H}_\beta(y)}{(x-(\log (1/|\epsilon|))^\gamma)^2}
\frac{|x-y|}{(\log (1/|\epsilon|))^{\gamma-1}}
dxdy\\
& \sim &
 \int_{P_3} \frac{\tilde{H}_\alpha(x)|x|^{1/\gamma-1}\tilde{H}_\beta(y)|y|^{1/\gamma-1}}{(x-(\log (1/|\epsilon|))^\gamma)^2 (((\log (1/|\epsilon|))^{\gamma})^{1/\gamma-1})^2}\\
& * & \frac{|x-y|}{(\log (1/|\epsilon|))^{\gamma-1}}
dxdy\\
& \sim &
 \int_{P_3} \frac{\tilde{H}_\alpha(x)|x|^{1/\gamma-1}\tilde{H}_\beta(y)|y|^{1/\gamma-1}}{(x-(\log (1/|\epsilon|))^\gamma)^2 (\log (1/|\epsilon|))^{2-2\gamma}}\\
 & * & \frac{|x-y|}{(\log (1/|\epsilon|))^{\gamma-1}}dxdy\\
 & \sim &
 \int_{P_3} \tilde{H}_\alpha(x)|x|^{1/\gamma-1}\tilde{H}_\beta(y)|y|^{1/\gamma-1}\\
& * &  
 \frac{|x-y|}{(x-(\log (1/|\epsilon|))^\gamma)^2 (\log (1/|\epsilon|))^{1-\gamma}} dxdy.\\
 \eea
 Hence
 \bea
 \sum _{n,m} h_{\alpha,n}h^\epsilon_{\beta,m} & <\sim &
 \int_{P_3} \tilde{H}_\alpha(x)|x|^{1/\gamma-1}\tilde{H}_\beta(y)|y|^{1/\gamma-1}\\
& * &  
 \frac{1}{|x-(\log (1/|\epsilon|))^\gamma| (\log (1/|\epsilon|))^{1-\gamma}} dxdy\\
   & <\sim &
 \int_{P_3} \tilde{H}_\alpha(x)|x|^{1/\gamma-1}\tilde{H}_\beta(y)|y|^{1/\gamma-1}\\
& * &  
 \frac{1}{\log (1/\delta) (\log (1/|\epsilon|))^{\gamma-1}(\log (1/|\epsilon|))^{1-\gamma}} dxdy\\
& \sim & \frac{1}{\log (1/\delta)} \int_{R_3} \tilde{H}_\alpha(x)|x|^{1/\gamma-1}\tilde{H}_\beta(y)|y|^{1/\gamma-1}dxdy\\
& \rightarrow & 0\\
\eea
\noindent because we can choose $\delta$ small enough.

\section{Theorem 7 for $R_{2Bii(II,IV)}$}

Recall from Lemma 24 that:
$$
-2nb\pi/a+1/a \log (1/|\epsilon|)-C<v,v'<-2nb\pi/a+1/a \log (1/|\epsilon|)+C.
$$

We estimate the contribution

\bea
h_{\alpha,n}h^\epsilon_{\beta,m}
& \sim & 
   \int_{|x-U|>(\log (1/|\epsilon|))^{\gamma-1}|v|}\tilde{H}_\alpha(x) \frac{ (\log (1/|\epsilon|))^{\gamma-1} v}
 { (x-U)^2}dx\\
& * &
 \int _{|y-U|>(\log (1/|\epsilon|))^{\gamma-1}|v|} \tilde{H}_\beta(y) \frac{ (\log (1/|\epsilon|))^{\gamma-1} v}
 { (y-U)^2}dy\\
 & \sim & 
   \int_{|x-(\log (1/|\epsilon|))^\gamma|>(\log (1/|\epsilon|))^{\gamma-1}|v|}\tilde{H}_\alpha(x) \frac{ (\log (1/|\epsilon|))^{\gamma-1} v}
 { (x-(\log (1/|\epsilon|))^\gamma)^2}dx\\
& * &
 \int _{|y-(\log (1/|\epsilon|))^\gamma|>(\log (1/|\epsilon|))^{\gamma-1}|v|} \tilde{H}_\beta(y) \frac{ (\log (1/|\epsilon|))^{\gamma-1} v}
 { (y-(\log (1/|\epsilon|))^\gamma)^2}dy\\
\eea

Note that when we sum over $n$, $v$ depends linearly on $n$
and as seen above, ranges from $\log 1/\delta$ to $d \log (1/|\epsilon|),$ $0<d<<1.$

Hence we need to estimate the expression $I(\alpha,\beta)$for given $(\alpha,\beta)$:

\bea
I(\alpha,\beta) & : = & \sum_{k=\log 1/\delta}^{d \log (1/|\epsilon|)}
  \int_{|x-(\log (1/|\epsilon|))^\gamma|>(\log (1/|\epsilon|))^{\gamma-1}k}\tilde{H}_\alpha(x) \frac{ (\log (1/|\epsilon|))^{\gamma-1} k}
 { (x-(\log (1/|\epsilon|))^\gamma)^2}dx\\
& * &
 \int _{|y-(\log (1/|\epsilon|))^\gamma|>(\log (1/|\epsilon|))^{\gamma-1}k} \tilde{H}_\beta(y) \frac{ (\log (1/|\epsilon|))^{\gamma-1} k}
 { (y-(\log (1/|\epsilon|))^\gamma)^2}dy\\
 \eea
 
 We introduce the integrals 
 \bea 
 I_{j,\alpha} & := &  
 \int_{(\log (1/|\epsilon|))^{\gamma-1}j<|x-(\log (1/|\epsilon|))^\gamma|<(\log (1/|\epsilon|))^{\gamma-1}(j+1)}\tilde{H}_\alpha(x) \frac{ (\log (1/|\epsilon|))^{\gamma-1} }
 { (x-(\log (1/|\epsilon|))^\gamma)^2}dx\\
 & \sim & 
  \int_{(\log (1/|\epsilon|))^{\gamma-1}j<|x-(\log (1/|\epsilon|))^\gamma|<(\log (1/|\epsilon|))^{\gamma-1}(j+1)}\tilde{H}_\alpha(x) \frac{1}{j^2(\log (1/|\epsilon|))^{\gamma-1} }dx\\
  & \sim & 
  \frac{1}{j^2}  \int_{(\log (1/|\epsilon|))^{\gamma-1}j<|x-(\log (1/|\epsilon|))^\gamma|<(\log (1/|\epsilon|))^{\gamma-1}(j+1)}\tilde{H}_\alpha(x) |x|^{1/\gamma-1}dx\\
  & = & 
\frac{1}{j^2}  \hat{I}_{j,\alpha},\\
\eea
and
\bea
I_{\infty,\alpha} & := & \int_{|x-(\log (1/|\epsilon|))^\gamma|>d(\log (1/|\epsilon|))^{\gamma}}\tilde{H}_\alpha(x) \frac{ (\log (1/|\epsilon|))^{\gamma-1} }
 { (x-(\log (1/|\epsilon|))^\gamma)^2}dx\\
 & \sim & 
  \int_{d(\log (1/|\epsilon|))^\gamma <|x-(\log (1/|\epsilon|))^\gamma|<Cd(\log (1/|\epsilon|))^{\gamma}}\tilde{H}_\alpha(x) \frac{1}{ (\log (1/|\epsilon|))^{\gamma+1} }dx\\
 & + &  \int_{|x-(\log (1/|\epsilon|))^\gamma|>C(\log (1/|\epsilon|))^{\gamma}}\tilde{H}_\alpha(x) \frac{ (\log (1/|\epsilon|))^{\gamma-1} }{x^2}dx\\
 & < \sim & \frac{1}{(\log (1/|\epsilon|))^2}
\int \tilde{H}_\alpha(x)|x|^{1/\gamma-1}dx\\
 & = & I^1_{\infty,\alpha}\\
\eea
 
\noindent and similarly for $\beta.$
 
 We get:
 
 \bea
 I(\alpha,\beta) & = & \sum_{k=\log 1/\delta}^{d \log (1/|\epsilon|)}
 \left[k\left( \left( \sum_{j=k}^{d \log (1/|\epsilon|)}I_{j,\alpha}\right)+I_{\infty,\alpha}\right)\right]
  \left[k\left( \left( \sum_{i=k}^{d \log (1/|\epsilon|)}I_{i,\beta}\right)+I_{\infty,\beta}\right)\right]\\
& \sim &  
  \sum_{k=\log 1/\delta}^{d \log (1/|\epsilon|)}k^2
 \left[\left( \left( \sum_{j=k}^{d \log (1/|\epsilon|)}\frac{\hat{I}_{j,\alpha}}{j^2}\right)+I_{\infty,\alpha}\right)\right]
  \left[\left( \left( \sum_{i=k}^{d \log (1/|\epsilon|)}\frac{\hat{I}_{i,\beta}}{i^2}\right)+I_{\infty,\beta}\right)\right]\\
 & = & 
  \sum_{k=\log 1/\delta}^{d \log (1/|\epsilon|)}k^2
 \left[\sum_{j=k}^{d \log (1/|\epsilon|)}\frac{\hat{I}_{j,\alpha}}{j^2}\right]
  \left[ \sum_{i=k}^{d \log (1/|\epsilon|)}\frac{\hat{I}_{i,\beta}}{i^2}\right]\\
& + &  
   \sum_{k=\log 1/\delta}^{d \log (1/|\epsilon|)}k^2
 I_{\infty,\alpha}
  \left[\sum_{i=k}^{d \log (1/|\epsilon|)}\frac{\hat{I}_{i,\beta}}{i^2}\right]
  +    \sum_{k=\log 1/\delta}^{d \log (1/|\epsilon|)}k^2
 \left[ \sum_{j=k}^{d \log (1/|\epsilon|)}\frac{\hat{I}_{j,\alpha}}{j^2}\right]
  I_{\infty,\beta}\\
& + & 
   \sum_{k=\log 1/\delta}^{d \log (1/|\epsilon|)}k^2
 I_{\infty,\alpha}
 I_{\infty,\beta}\\
 & = & I+II+III+IV\\
\eea 
 
 Here $II$ and $III$ are symmetric. It suffices to estimate $II.$\\

 We estimate $IV$ first. Since $\sum k^2 \sim (\log (1/|\epsilon|))^3,$ this is immediately small when multiplied with $I_{\infty,\alpha}, I_{\infty,\beta}.$
 For $II$, we get:
 
 \bea
  II & = & 
  \sum_{k=\log 1/\delta}^{d \log (1/|\epsilon|)}k^2
 I_{\infty,\alpha}
  \left[\sum_{i=k}^{d \log (1/|\epsilon|)}\frac{\hat{I}_{i,\beta}}{i^2}\right]\\
  & < & 
 I_{\infty,\alpha} \sum_{k=\log 1/\delta}^{d \log (1/|\epsilon|)}
  \left[\sum_{i=k}^{d \log (1/|\epsilon|)}\hat{I}_{i,\beta}\right]\\
  & < & \frac{1}{(\log (1/|\epsilon|))}\int \tilde{H}_\alpha(x) |x|^{1/\gamma-1}dx\int \tilde{H}_\beta(y) |y|^{1/\gamma-1}dy\\
  & \rightarrow & 0\\
 \eea
 
 Finally we estimate $I$.
 
 \bea
  I & = &  \sum_{k=\log 1/\delta}^{d \log (1/|\epsilon|)}k^2
 \left[\sum_{j=k}^{d \log (1/|\epsilon|)}\frac{\hat{I}_{j,\alpha}}{j^2}\right]
  \left[ \sum_{i=k}^{d \log (1/|\epsilon|)}\frac{\hat{I}_{i,\beta}}{i^2}\right]\\
& < &  
  \sum_{k=\log 1/\delta}^{d \log (1/|\epsilon|)}\frac{1}{k^2}
 \left[\sum_{j=k}^{d \log (1/|\epsilon|)}\hat{I}_{j,\alpha}\right]
  \left[ \sum_{i=k}^{d \log (1/|\epsilon|)}\hat{I}_{i,\beta}\right]\\
\eea
  
  We can make this as small as we wish by choosing $\delta$ small.

\section{Proof of Theorem 4}

\begin{proof}

We use the approach in  \cite{FS2005}.

Let $T$ be a positive harmonic current directed by $\mathcal F$. We want to show that
$\int T \wedge T=0.$ Let $T_\epsilon =(\Phi_\epsilon)_*T$ and define
$T_\epsilon^\delta$ as the average of $T_\epsilon$ using a small neighborhood
of identity in $U(3).$ 
Then since $T_\epsilon \rightharpoonup T,$ we have
$\int T \wedge T = \lim_{\epsilon \rightarrow 0} \int T \wedge T_\epsilon.$
On the other hand $T_\epsilon^\delta=
\omega+\partial S_\epsilon^\delta+ \overline{\partial S}_\epsilon^\delta
+i\partial \overline{\partial}u_\epsilon^\delta$ and $S_\epsilon^\delta \rightarrow S_\epsilon$ in $L^2.$
So $\int T \wedge T_\epsilon =\lim_{|\delta|,|\delta' |\rightarrow 0,
|\delta|,
|\delta'|<< \epsilon} \int T_\epsilon^\delta \wedge T^{\delta'}.$ Hence as in \cite{FS2005}
it is enough to show that

$$
\lim_{\delta,\delta',\epsilon \rightarrow 0, |\delta|,|\delta'|<<|\epsilon|}
\int T^\delta_\epsilon \wedge T^{\delta'} =0.
$$

 We can compute the geometric intersection $T^\delta_\epsilon \wedge T^{\delta'}$
and it is enough to estimate $T_\epsilon \wedge_g T$. Recall that if $\phi$ is a test function supported
in $B,$ then we define
$$
<T_\epsilon \wedge_g T, \phi> =
 \int \sum_{J^\epsilon_{\alpha,\beta}}\phi(p)  H_\alpha(p)H^\epsilon_\beta(p)d\mu(\alpha)
d\mu(\beta).
$$
\noindent where $J^\epsilon_{\alpha,\beta}$ consists of intersection points
of $\Delta_\alpha$ and $\Delta^\epsilon_\beta.$ 
The following Lemma is proved in \cite{FS2005}.

\begin{lemma}
We have that $\int T \wedge T_\epsilon= \int T \wedge_g T_\epsilon.$ 
The same holds for $T^\delta, T^{\delta'}_\epsilon.$
\end{lemma}

$$
<T_\epsilon \wedge_g T, \phi> \leq
C\|\phi\|_\infty \int \sum_{J^\epsilon_{\alpha,\beta}} H_\alpha(p)H^\epsilon_\beta(p)d\mu(\alpha)
d\mu(\beta),
$$

We know that the number of points in $J^\epsilon_{\alpha,\beta}$ is bounded by a fixed constant independent of $\epsilon.$ For $p$ out of a fixed neighborhood of the singularities the integral converges to zero. This is the  case considered in
\cite{FS2005}. So it is enough to show that for $\delta>0$ small enough

$$
J_\epsilon (\delta):= \int \sum_{J^\epsilon_{\alpha,\beta} } H_\alpha(p)H^\epsilon_\beta(p)
d\mu(\alpha) d\mu(\beta)$$
\noindent is arbitrarily small. This is precisely the content of Theorem 7, since all estimates
are valid after composition by automorphisms in a small neighborhood of $U(3).$

Consequently if $T_1,T_2$ are two such currents then $\int \frac{T_1+T_2}{2} \wedge
\frac{T_1+T_2}{2}=0.$ Hence $\int T_1 \wedge T_2=0,$ therefore
$T_1,T_2$ are proportional.

\end{proof}

We give a dynamical consequence of the uniqueness of the harmonic current for ${\mathcal F}
\in {\mathcal H}(d),$ here ${\mathcal H}(d)$ is the Zariski open set of foliations of degree $d,$ introduced in Theorem 2. 

\begin{corollary}
Let ${\mathcal F} \in {\mathcal H}(d)$. Let $\phi:\Delta \rightarrow L$ be the universal covering of a leaf $L.$ Let $\tau_r: = \frac{\phi_* [\log^+ \frac{r}{|z|}\Delta_r]}
{\|\phi_*[\log^+ \frac{r}{|z|}\Delta_r]\|}.$ Then $\lim_{r \rightarrow 1}\tau_r = T,$ where $T$ is the unique harmonic current directed by ${\mathcal F}.$
\end{corollary}

Here $\Delta_r$ denotes the disc of center $0$ and radius $r$. The Corollary which is a consequence of paragraph 5 in \cite{FS2005} says that the normalized images of $[\log^+\frac{r}{|z|} \Delta_r]$ converge to $T.$ This is similar to the pointwise ergodic theorem, since we are averaging on an orbit.

Recall that the limit set of a leaf $L$ is defined as 
$\lim(L) = \cap_n \overline{L \setminus K_n},$ where $K_n \subset K_{n+1}$ is an exhaustion of $L$
by compact sets. One of the main questions in foliation theory is to describe the limit set
of a foliation ${\mathcal F}$: $\lim({\mathcal F}):= \overline{\cup_{L \in \mathcal F}\lim(L)}.$ Corollary 2 implies in particular that for ${\mathcal F } \in {\mathcal H}(d),$ for every leaf $L \in \mathcal F$,
$\lim(L)$ contains $ \mbox{supp}(T).$ Indeed as shown in \cite{FS2005},

$$
\|\Phi_* \left[ \log^+ \frac{r}{|z|} \Delta_r\right] \| \rightarrow \infty
$$

\noindent as $r \rightarrow 1.$ Hence supp$(T) \subset \overline{L \setminus K_n}$ for every $n.$

\begin{corollary}
The map $\lambda \rightarrow T_\lambda$ is continuous from ${\mathcal H}(d)$ with values in the positive harmonic currents of mass one. Let ${\mathcal F}_\lambda$ be a holomorphic family of
foliations in ${\mathcal H}(d).$ Let $(T_\lambda)$ be the associated currents. If a hyperbolic point
$p_0\in {\mbox{Supp}}(T_{\lambda_0}),$ then the perturbed hyperbolic point $p_\lambda$ belongs to
Supp$(T_\lambda)$.
\end{corollary}

\begin{proof}
Assume ${\mathcal F}_{\lambda_n} \rightarrow {\mathcal F}_{\lambda_0}$ in ${\mathcal H}(d).$
Let $(T_{\lambda_n})$ be the normalized positive harmonic currents associated to ${\mathcal F}_{\lambda_n}$. Since $\|T_{\lambda_n}\|=1,$ the sequence $(T_{\lambda_n})$ has cluster points. It is clear that any cluster point $S$ is positive harmonic and directed by ${\mathcal F}_{\lambda_0}.$ So $S=T_{\lambda_0}$ by uniqueness. Assume  the support of $T_{\lambda_0}$ intersects a ball 
$B(p_0,r)$ where $p_0$ is a hyperbolic singular point of ${\mathcal F}_{\lambda_0}$ and the ball is contained in the common domain of linearization of $p_\lambda\in {\mbox{Sing}}
({\mathcal F}_{\lambda}), p_\lambda \rightarrow p_0,$ $p_\lambda$ hyperbolic. 

{} From our local study of positive harmonic currents near a hyperbolic singular point $p_0
\in {\mbox{Supp}}(T_{\lambda_0}).$ Since $T_{\lambda} \rightarrow T_{\lambda_0}$, $T_{\lambda}$ gives mass to $B(p_0,r)$, applying again the local study for $T_{\lambda}$
we get that $p_\lambda \in {\mbox{Supp}}(T_{\lambda}).$

\end{proof}

 \begin{remark}
Let $f$ be a holomorphic endomorphism of ${\mathbb P}^2.$ Let ${\mathcal F}$ be a foliation with
only hyperbolic singularities. Then $f^*{\mathcal F}$ is a foliaton and its singularities are not necessarily
hyperbolic. However there is only one positive harmonic current of mass $1$, directed by $f^*{\mathcal F}.$
Indeed let $T$ be any such current. We will show that $\int T \wedge T =0$ which implies
the uniqueness. Observe that $f_*T$ is a current directed by ${\mathcal F}.$ Hence $\int f_*T
\wedge f_*T=0.$ Since $f^*$ is a finite covering of degree $d^2$ we have 
$$
\int T \wedge T\leq \int f^* \left[f_*T \wedge f_* T\right]=d^2\int f_*T \wedge f_*T=0.
$$
\end{remark}

\section{Measure associated to a harmonic current}

Let ${\mathcal F} \in {\mathcal H}(d)$ be a holomorphic foliation as in Theorem 2. We know that there
is a unique positive harmonic current $T$ of mass one directed by ${\mathcal F}.$

We are going to associate to $T$ a conformal, measurable metric along leaves that we will denote by $g_T$ and also a positive $\underline{{\mbox{finite}}}$ measure $\mu_T$ which is invariant under the harmonic flow associated also to $T.$ The metric $g_T$ and the measure $\mu_T$ where first considered by S. Frankel, in the non singular case \cite{F1995} he proved in that case a version of Proposition 2 and Proposition 3.

On a flow box $B$ disjoint from $E={\mbox{Sing}}({\mathcal F}),$ the current $T$ can be written

$$
T=\int h_\alpha [V_\alpha]d\mu(\alpha)
$$

\noindent where $h_\alpha$ are positive harmonic functions and $\mu$ is a positive measure on a transversal $A.$ The $[V_\alpha]$ are the currents of integration on plaques. On $B, \partial T=
\tau \wedge T$ with $\tau=\frac{\partial h_\alpha}{h_\alpha},\; \mu$ almost everywhere. Observe that
$\tau$ is independent of the choice of $h_\alpha:$ if we replace $h_\alpha$ by $c_\alpha h_\alpha,
c_\alpha \in \mathbb R^+$ then $\tau$ is unchanged.

\medskip

We define the metric $g_T$ on leaves by $g_T= \frac{i}{2}\tau \otimes \overline{\tau}.$ Along the plaque $V_\alpha$ with a choice of coordinate $(z_\alpha)$ we have

$$
g_T = \frac{i}{2} \left| \frac{\partial h_\alpha}{\partial z_\alpha}\right|^2 \frac{1}{h_\alpha^2}dz_\alpha
\otimes d \overline{z}_\alpha\;\;\;\;\;\;\;(1)
$$

Define ${\mathcal C}_T=\{(\alpha,z);\frac{\partial h_\alpha}{\partial z}(\alpha,z)=0\}$ it's the critical set of the "metric" $g_T.$ We also define the current of bidegree $(2,2)$, $\mu_T,$ which we identify with a measure

$$
\mu_T:= i\tau \wedge \overline{\tau}\wedge T.
$$

In local coordinates in a flow box $B$, we have:

$$
\mu_T= \int d\nu(\alpha) \int_{[V_\alpha]}\left|\frac{\partial h_\alpha}{\partial z_\alpha}\right|^2\frac{1}{h_\alpha}(idz_\alpha \wedge d\overline{z_\alpha}).\;\;\;\;\;(2)
$$

\begin{proposition}
Let ${\mathcal F} \in {\mathcal H}(d).$ The metric $g_T$ has constant negative curvature out of the set $\mathcal{C}_T$ where the metric vanishes.
\end{proposition}

\begin{proof}
Since the current $T$ is unique, every measurable set of leaves $\mathcal A$ has zero or full measure with respect to $\|T\|.$ Define ${\mathcal N}_g:=
\{{\mbox{leaves on which}}\; g_T\; {\mbox{vanishes identically}}\}.$ Since $h_\alpha$ is measurable,
then ${\mathcal N}_g$ is measurable. So ${\mathcal N}_g$ is of zero or full measure. But if 
${\mathcal N}_g$ is of full measure, $\partial T=0$ and by conjugation $\overline{\partial}T=0,$ 
hence $T$ is closed. A foliation ${\mathcal F}$ in ${\mathcal H}(d)$ admits no positive closed current directed by ${\mathcal F}$ since all singularities are hyperbolic. So ${\mathcal N}_g$ is of zero $\|T\|$ measure.

{}From (1) it is clear that the metric is conformal. On a flow box $B,$ the curvature $\kappa(g)$ has the following expression out of $C_T.$ The curvature is given by
$$
\kappa(g)=-\frac{1}{4} \frac{\Delta \log g}{g}=\frac{1}{2} \frac{\Delta \log h_\alpha}{\left|\frac{\partial h_\alpha}{\partial z_\alpha}\right|^2 \frac{1}{h_\alpha^2}}.
$$

So 

$$
\kappa(g_T)=  \frac{h_\alpha^2}{|h_{\alpha,z}|^2} \left( \frac{\partial}{\partial \overline{z}}
\left( \frac{h_{\alpha,z}}{h_\alpha}\right)\right).
$$

Since $h_\alpha$ is harmonic we get $\kappa(g_T)=-1.$

\end{proof}
Because of the nature of the singularities, the leaves are uniformized by the unit disc $\Delta.$ Let $g$ denote the Poincar\'e metric on leaves. We choose a normalization so that the curvature $\kappa(g)$ of $g$ on leaves is $-1.$ 

\begin{proposition}
Let $T$ be the harmonic current associated to ${\mathcal F}\in {\mathcal H}(d).$ If $g_T$ is the associated metric on leaves, then $g_T\leq g.$
\end{proposition}

\begin{proof}
We have normalized the metric $g_T$ so that on each leaf $L_\alpha,$  $g_T$ has curvature $-1$
on $L_\alpha \setminus \mathcal{C}(T).$ The Ahlfors' Schwarz lemma, applied to the abstract Riemann surface $L_\alpha \setminus \mathcal{C}_T$ implies that $g_T \leq g.$
\end{proof}

 We will denote by $\Phi_\alpha:\Delta \rightarrow \mathcal L_\alpha,$ the uniformizing map from 
 $\Delta$ to $\mathcal L_\alpha.$ When we fix a transversal $A$ in a flow box we can choose for each 
 $\alpha \in A$ a uniformizing map $\Phi_\alpha(0)=\alpha,$ then $\Phi_\alpha$ vary measurably. We will denote by $\Gamma_\alpha$ the group of deck transformations for the map $\Phi_\alpha.$

We want to define a vector field $\chi$ on ${\mathcal F}$ associated to the current $T.$ The vector field will be defined as the metric $g_T$ only $\|T\|$ a.e. On $L_\alpha, \chi_\alpha$ is collinear with
the gradient field of $h_\alpha.$ We define $\chi_\alpha$ on a flow box with local coordinates
$z_\alpha=x_\alpha+iy_\alpha$ by

$$
\chi_\alpha:= c \frac{h_\alpha}{|h_z|^2}(h_{x_\alpha},h_{y_\alpha}).
$$

We choose the constant  $c$ so that $g_T(\chi_\alpha, \chi_\alpha)=1.$ 
The vector field $\chi_\alpha$ is independent of the choice of $h.$ It blows up at every point of 
$\mathcal{C}_T.$ Which means that the integral curves of $\chi_\alpha$ approach these points at infinite speed. So we have to take out these trajectories in order to have a well defined flow. Observe that the set of these trajectories is of $\mu_T$ measure zero.  It is clear that the integral curves of $\chi_\alpha$ are along the level sets of the harmonic conjugates of $h_\alpha$ such that $f_\alpha= h_\alpha+
iv_\alpha$ is holomorphic.

\begin{theorem}
Let $T$ be the positive harmonic current associated to ${\mathcal F}\in {\mathcal H}(d).$ Then the measure $\mu_T$ is $\underline{{\mbox{finite}}}.$  Moreover, if ${\mathcal {F}}_\lambda$ is a holomorphic family of foliations in 
${\mathcal {H}}(d)$, $\lambda \in \Delta(\lambda_0,r),$ then the mass of $\mu_{T_\lambda}$ near
 hyperbolic singularities is uniformly small in a fixed neighborhood of the singularities.
\end{theorem}

\begin{proof}
For a flow box $B$ away from the singularities, it is clear that $\mu_T$ has finite mass. Indeed the functions $h_\alpha$ are positive harmonic, and by Harnack $\frac{h_\alpha}{|\partial h_\alpha|}\leq c,$ hence $\mu_T$ has finite mass in $B.$ It is enough to show that $\mu_T$ has finite mass in a flow box $B_i$ near a hyperbolic singularity given by $\omega=zdw-\lambda wdz, \lambda=a+ib,b \neq 0.$ We use the parametrization 
$$
\psi_\alpha(\zeta)=(e^{i(\zeta+(\log |\alpha|)/b)}, \alpha e^{i\lambda (\zeta+(\log |\alpha|)/b)})
$$

\noindent by a sector near  the hyperbolic singularity. Since $\psi_\alpha^* h_\alpha=H_\alpha$ is a positive harmonic function and $\mu$ a.e. $H_\alpha(\zeta) \rightarrow 0$ when $\Im \zeta \rightarrow
 +\infty,$ then again by Harnack 
 $\psi^*_\alpha(\tau)$ is bounded. The total mass of $\mu_T$ in $B_i$ satisfies

$$
\int_{B_i}\mu_T \leq
\int_{D(w_0,r) \times S_\lambda}i\psi^*_\alpha(\tau) \wedge \psi^*_\alpha(\overline{\tau})
\wedge \psi^*_\alpha [V_\alpha]H_\alpha d\mu(\alpha)
$$

\noindent $\psi^*_\alpha[V_\alpha]$ is a graph in the flow box. It is of bounded area and $\int_{D(w_0,r)}H_\alpha d\mu(\alpha)$ defines a bounded harmonic function. So the mass $\mu_T$ is bounded near the origin.

\medskip

 Basically the slicing of $\mu_T$ along the leaves gives the area measure on leaves associated to the metric $g_T.$
Let $T_\lambda$ be the current associated to ${\mathcal F}_\lambda,$ and let $\mu^\lambda$ denote the corresponding measure on a transversal. The linearizations associated to a holomorphically varying hyperbolic singularity vary holomorphically. Then $\int H^\lambda_\alpha d\mu^\lambda (\alpha) \rightarrow 0$ when $\Im \zeta \rightarrow +\infty,$ uniformly when $\lambda$ is near
$\lambda_0.$ (We don't say that $H_\alpha^\lambda$ vary holomorphically.) So the mass of 
$\mu_{T_\lambda}$ is uniformly small in a fixed neighborhood of the singularities if $\lambda$ is close enough to $\lambda_0.$

\end{proof}

\begin{theorem}
Let $\lambda \rightarrow \mathcal F_\lambda$ be a holomorphic family of foliations in $\mathcal H(d),$ parametrized by a disc $\Delta.$ Then $\lambda \rightarrow \mu_{\lambda}$ is a continuous family of measures.
\end{theorem}

\begin{proof}
Let $(T_\lambda)$ be the family of the positive harmonic currents directed by $\mathcal F_\lambda$. Recall that $\mu_{T_\lambda}=i\tau_\lambda \wedge \overline{\tau}_\lambda \wedge T_\lambda.$ 

\medskip

Fix a flow box $B$ for $\mathcal F_{\lambda_0}$ away from the singularities. We can consider 
$(\phi_\lambda)$ local biholomorphisms straightening ${\mathcal F}_\lambda$ in $B,$ when $\lambda \rightarrow \lambda_0.$ We know that the currents $S_\lambda:=(\phi_\lambda)_* {T_\lambda}$ depend continuously on $\lambda.$ We can write in $B$,

$$
S_\lambda =\int [w=\alpha] h_\alpha^\lambda (z) d \mu_\lambda(\alpha)
$$

\noindent where $\mu_\lambda$ is the measure on a fixed transversal $(z=z_0).$ We can assume that 
$h_\alpha^{\lambda}(z_0) =1$ for all $\alpha,\lambda.$ 

\medskip

Since $S_\lambda \rightarrow S_{\lambda_0}$ then for every $z$ we have $h_\alpha^\lambda(z) \mu_\lambda(\alpha) \rightarrow h_\alpha^{\lambda_0}\mu_{\lambda_0 }(\alpha)$
weakly when $\lambda \rightarrow \lambda_0.$

\medskip

The $(h^\lambda_\alpha)^2$ also vary slowly, by Harnack, so we also get that $\lambda \rightarrow
 (h_\alpha^\lambda(z))^2 \mu_\lambda(\alpha)$ is continuous for every $z.$ Define

$$
U_\lambda:= \int [w=\alpha] (h_\alpha^\lambda)^2 (z) d \mu_\lambda(\alpha).
$$

The family of positive currents $U_\lambda$ is also continuous because $(h_\alpha^\lambda)^2$ is uniformly bounded.
It follows that $\lambda \rightarrow i \partial \overline{\partial} U_\lambda$ is continuous i.e.

$$
\lambda \rightarrow \int |h^\lambda_{\alpha, z}|^2 [w=\alpha] d \mu_\lambda(\alpha).
$$

Using again Harnack inequalities for $\frac{1}{h_\alpha^{\lambda^2}},$ we find that 
 $\lambda \rightarrow \mu_{T_\lambda}$ is continuous in $B.$

\medskip

We have seen in Theorem 9 that $\mu_{T_\lambda}$ has uniformly small mass near the singularities. Hence $\lambda \rightarrow \mu_{T_\lambda}$ is continuous.

\end{proof}

Let $|g_T^\alpha|$ denote the measure induced by the metric $g_T$ on the leaf $L_\alpha.$ We will omit $\alpha,$ most frequently. 

\medskip

We will say that a set $E$ is invariant if up to a set of $\mu_T$ measure zero, it is a union
of orbits of $\chi.$ For a measurable set $E$ we denote by $E_\alpha$ the intersection 
$E \cap L_\alpha.$

\begin{theorem}
Either there is an invariant set $E$ for $\chi$ such that for $\|T\|$ almost every leaf $L_\alpha,$
$|g_T|(E_\alpha)>0$ and $|g_T|(E^c_\alpha)>0$ or the measure
$\mu_T$ is ergodic.
\end{theorem}

\begin{proof}
Fix a countable family $(B_i)$ of flow boxes such that $\cup_i B_i=
\mathbb P^2 \setminus ({\mbox{Sing}}(\mathcal F).$ Let $E$ be an invariant set for $\chi$ such that
$\mu_T(E)>0.$ Define $E_i=\{\alpha; |g_T| (L_\alpha \cap E \cap B_i)=0\}.$
$\mathcal E := \cap_i E_i$ is measurable. It is a union of leaves. Since the current $T$ is unique
and $\mu_T({\mathcal E})=0$, then $\|T\|$ almost every leaf is in ${\mathcal E}^c.$

\medskip

For $L_\alpha \in {\mathcal E}^c, |g_T|(E_\alpha)>0.$ We can do a similar construction
for $E^c$ if $\mu_T(E^c)>0.$ We then get a set of $\|T\|$ full measure of leaves such that
$|g_T|(E_\alpha) >0$ and $|g_T|(E^c_\alpha)>0.$

\end{proof}

\end{document}